\def\paperversion{2}
\ifnum\paperversion=2
\documentclass[opre, nonblindrev]{informs4}
\SingleSpacedXI

\renewcommand{\theARTICLETOP}{}
\renewcommand{\theARTICLERULE}{\vspace{12pt}}
\renewcommand{\theRRHFirstLine}{}
\renewcommand{\theRRHSecondLine}{}
\renewcommand{\theLRHFirstLine}{}
\renewcommand{\theLRHSecondLine}{}
\renewcommand{\theECLRHFirstLine}{}%
\renewcommand{\theECLRHSecondLine}{}%
\renewcommand{\theECRRHFirstLine}{}%
\renewcommand{\theECRRHSecondLine}{}%
\renewcommand{\endproof}{\endTrivlist\addvspace{2ex}}
\makeatletter
\renewcommand\section{\@startsection {section}{1}{\z@}{-13pt plus -6pt minus -3pt}{4pt}%
  {\fs.13.15.\bfseries\RAGG}}%
\renewcommand\subsection{\@startsection{subsection}{2}{\z@}{-13pt plus -6pt minus -3pt}{4pt}%
  {\TEN\bfseries\RAGG}}%
\makeatother
\let\FigureCaptionFontStyle\undefined
\def\FigureCaptionFontStyle{\mdseries}

\usepackage[scaled=.98,sups,osf]{XCharter}%

\fi
\ifnum\paperversion=1
\documentclass[opre, blindrev]{informs4}
\OneAndAHalfSpacedXI
\fi
\usepackage[numbers]{natbib}
 \def\bibfont{\small}%
\TheoremsNumberedThrough   
\ECRepeatTheorems
\theoremstyle{TH}%
\newtheorem{condition}{Condition}
\usepackage{endnotes}
\DeclareMathAlphabet{\mathsfit}{T1}{\sfdefault}{\mddefault}{\sldefault}
\SetMathAlphabet{\mathsfit}{bold}{T1}{\sfdefault}{\bfdefault}{\sldefault}
\DeclareMathAlphabet{\mathcal}{OMS}{cmsy}{m}{n}
\usepackage{anyfontsize}
\let\footnote=\endnote
\usepackage{booktabs}
\usepackage{bm}
\usepackage{subfigure}

\usepackage[flushleft]{threeparttable}

\DeclareMathAlphabet{\mathsfit}{T1}{\sfdefault}{\mddefault}{\sldefault}
\SetMathAlphabet{\mathsfit}{bold}{T1}{\sfdefault}{\bfdefault}{\sldefault}
\DeclareMathAlphabet{\mathcal}{OMS}{cmsy}{m}{n}
\usepackage[table, svgnames, dvipsnames, clash]{xcolor}
\usepackage[hypertexnames=false]{hyperref}
\hypersetup{%
  colorlinks=true,%
  linkcolor={blue!66!black},%
  linkbordercolor=red,%
  citecolor={blue!66!black},
}
\usepackage{algorithm}
\usepackage{algorithmic}
\usepackage{eqparbox}

\usepackage{etoolbox}  %
\makeatletter
\patchcmd{\algorithmic}{\addtolength{\ALC@tlm}{\leftmargin} }{\addtolength{\ALC@tlm}{\leftmargin}}{}{}
\makeatother

\usepackage{mathtools}
\usepackage{booktabs,multirow}
\usepackage{enumitem}
\setlist[enumerate,1]{label=\normalfont{(\Roman*)},leftmargin=*}

\newcommand*{\QED}{%
\leavevmode\unskip\penalty9999 \hbox{}\nobreak\hfill
    \quad\hbox{$\square$}%
}
\newcommand*{\QEG}{%
\leavevmode\unskip\penalty9999 \hbox{}\nobreak\hfill
    \quad\hbox{$\clubsuit$}%
}
\newcommand*{\QDEF}{%
\leavevmode\unskip\penalty9999 \hbox{}\nobreak\hfill
    \quad\hbox{$\diamondsuit$}%
}

\newcommand{\cP}{\mathcal{P}}

\newcommand{\bP}{\mathbb{P}}
\newcommand{\bR}{\mathbb{R}}
\newcommand{\bQ}{\mathbb{Q}}

\newcommand{\cN}{\mathcal{N}}
\newcommand{\ox}{\overline{x}}
\newcommand{\hx}{\hat{x}}
\newcommand{\W}{\mathcal{W}}
\newcommand{\esssup}{\operatorname{ess}\sup}
\newcommand{\orho}{\overline{\rho}}
\newcommand{\hP}{\widehat{\mathbb{P}}}
\newcommand{\hF}{\widehat{F}}
\newcommand{\proj}{\mathrm{Proj}}
\newcommand{\Reg}{\epsilon}
\newcommand{\fh}{\mathfrak{h}}
\newcommand{\Primal}{{}}
\newcommand{\Dual}{\mathrm{D}}
\newcommand{\trans}{^{\mathrm T}}
\newcommand{\frakM}{\mathfrak{M}}
\newcommand{\diff}{\,\mathrm{d}}

\newcommand{\Var}{\mathbb{V}\mathrm{ar}}

\newcommand{\cO}{\mathcal{O}}
\newcommand{\tO}{\widetilde{\mathcal{O}}}

\newcommand{\cZ}{\mathcal{Z}}
\newcommand{\cM}{\mathcal{M}}
\newcommand{\bE}{\mathbb{E}}
\newcommand{\prox}{\mathrm{Prox}}
\newcommand{\nout}{n}

\renewcommand{\emph}[1]{\textit{#1}}

\DeclarePairedDelimiterX{\inp}[2]{\langle}{\rangle}{#1, #2}
\definecolor{longhorn}{rgb}{0.8, 0.33, 0.0}

\newcommand{\Jie}[1]{{\textcolor{black}{#1}}}

\TITLE{Sinkhorn Distributionally Robust Optimization}
\ARTICLEAUTHORS{%
\AUTHOR{Jie Wang}
\AFF{School of Industrial and Systems Engineering, Georgia Institute of Technology, Atlanta, GA 30332, \EMAIL{jwang3163@gatech.edu}}
\AUTHOR{Rui Gao}
\AFF{Department of Information, Risk, and Operations Management, University of Texas at Austin, Austin, TX 78712, \EMAIL{rui.gao@mccombs.utexas.edu}}
\AUTHOR{Yao Xie}
\AFF{School of Industrial and Systems Engineering, Georgia Institute of Technology, Atlanta, GA 30332, \EMAIL{yao.xie@isye.gatech.edu}}
} %

\ABSTRACT{
We study distributionally robust optimization with Sinkhorn distance---a variant of Wasserstein distance based on entropic regularization.
We derive a convex programming dual reformulation for general nominal distributions, transport costs, and loss functions.
\Jie{
To solve the dual reformulation, we develop a stochastic mirror descent algorithm with biased subgradient estimators and derive its computational complexity guarantees.
}
Finally, we provide numerical examples using synthetic and real data to demonstrate its superior performance.
}

\RUNAUTHOR{Wang, Gao, and Xie}
\RUNTITLE{Sinkhorn DRO}

\begin{document}
\maketitle

\section{Introduction}
Decision-making problems under uncertainty arise in various fields such as operations research, machine learning, engineering, and economics. In these scenarios, uncertainty in the data arises from factors like measurement error, limited sample size, contamination, anomalies, or model misspecification. 
Addressing this uncertainty is crucial to obtain reliable and robust solutions. 
In recent years, Distributionally Robust Optimization (DRO) has emerged as a promising data-driven approach to tackle these challenges.
DRO aims to find a minimax robust optimal decision that minimizes the expected loss under the most adverse distribution within a predefined set of relevant distributions, known as an ambiguity set. This approach provides a principled framework to handle uncertainty and obtain solutions that are resilient to distributional variations. It goes beyond the traditional sample average approximation (SAA) method used in stochastic programming and offers improved out-of-sample performance.
For a comprehensive overview of DRO, we refer interested readers to the recent survey by \Jie{\cite{daneilsurvey}}. 

At the core of distributionally robust optimization lies the crucial task of selecting an appropriate ambiguity set. 
An ambiguity set should strike a balance between computational \Jie{efficiency} and practical interpretability while being rich enough to encompass relevant distributions and avoiding unnecessary ones that may lead to overly conservative decisions.
In the literature, various formulations of DRO have been proposed, among which the ambiguity set based on Wasserstein distance has gained significant attention in recent years \citep{wozabal2012framework, Mohajerin18, blanchet2019quantifying, gao2016distributionally}.
The Wasserstein distance incorporates the geometry of the sample space, making it suitable for comparing distributions with non-overlapping supports and hedging against data perturbations \citep{gao2016distributionally}. 
The Wasserstein ambiguity set has received substantial theoretical attention, with provable performance guarantees \citep{shafieezadeh2019regularization, blanchet2019robust, blanchet2019confidence, blanchet2021statistical, gao2020finitesample}. 
Empirical success has also been demonstrated across a wide range of applications, including operations research \citep{blanchet2018distributionallyportfolio, cherukuri2019cooperative, singh2020tight, nguyen2021robustifying, SINGH2021121}, machine learning \citep{Shafieezadeh15, chen2019selecting, Luo19, blanchet2019multivariate, nguyen2020robust}, stochastic control \citep{yang2017convex, YangWDRO20, wang2021reliable}, and more.

However, the current Wasserstein DRO framework has its limitations. First, \Jie{
the computational efficiency of Wasserstein DRO is achieved under somewhat stringent conditions, as its dual formulation involves a subproblem that requires the global supremum of some regularized loss function over the sample space.
Let $\min_{\theta\in\Theta}\max_{\bP\in\frakM}~\bE_{z\sim \bP}[f_{\theta}(z)]$ denote a typical Wasserstein DRO formulation, where the loss function $f_{\theta}(z)$ is convex in $\theta$ belonging to a closed and convex feasible region $\Theta$, and the ambiguity set $\frakM$ is centered around a nominal distribution $\hP$ and contains distributions supported on a space $\cZ$. 
Table~\ref{Tab:compare:DRO} summarizes the known cases where solving Wasserstein DRO is computationally efficient.
}
One general approach to solving it is to use a finite and discrete grid of scenarios to approximate the entire sample space.
This involves solving the formulation restricted to the approximated sample space \citep{pflug2007ambiguity, chen2021decomposition, liu2021discrete}, but suffers from the curse of dimensionality. 
Simplified convex reformulations are known when the loss function can be expressed as a pointwise maximum of finitely many concave functions \citep{esfahani2018data, gao2016distributionally, shafieezadeh2023new}, or when the loss is the generalized linear model \citep{Shafieezadeh15, yu2022fast, selvi2022wasserstein, shafieezadeh2023new}. 
In addition, efficient first-order algorithms have been developed for Wasserstein DRO with strongly convex transport cost, smooth loss functions, and sufficiently small radius (or equivalently, sufficiently large Lagrangian multiplier) so that the involved subproblem becomes strongly convex \citep{sinha2018certifiable, blanchet2021optimal}. 
However, beyond these conditions on the loss function and the transport cost, solving Wasserstein DRO becomes a computationally challenging task. 
Second, from a \emph{modeling} perspective, in data-driven Wasserstein DRO, where the nominal distribution is finitely supported (usually the empirical distribution), the worst-case distribution is shown to be a discrete distribution \cite{gao2016distributionally} (which is unique when the regularized loss function has a unique maximizer).
This is the case even though the underlying true distribution in many practical applications may be continuous. 
Consequently, concerns arise regarding whether Wasserstein DRO hedges the right family of distributions and whether it induces overly conservative solutions.

\begin{table}[ht] 
\small
\caption{\Jie{Known cases of Wasserstein DRO where it is computationally efficient to solve}
}
\label{Tab:compare:DRO}
		\centering
			\setlength{\tabcolsep}{2pt}\renewcommand{\arraystretch}{1.2}
			\begin{threeparttable}
			\begin{tabular}{c| c| c |c| c}
				\toprule  
				\multicolumn{1}{c|}{References}&	\multicolumn{1}{c|}{Loss function $f_{\theta}(z)$} &	\multicolumn{1}{c|}{ Transport cost} & 	\multicolumn{1}{c|}{Nominal distribution $\hP$} & Support $\cZ$ \\
				\hline
				\multirow{2}{2cm}{
				\citep{pflug2007ambiguity, chen2021decomposition, liu2021discrete}
				}	&  
				\multirow{2}{2.6cm}{\centering General} &  
				\multirow{2}{2cm}{\centering General} &  
				\multirow{2}{3.7cm}{\centering General} & 
				\multirow{2}{3.7cm}{\centering Discrete and finite set} \\
				& &   &   & \\
				\cline{1-5}
                    \multirow{2}{2cm}{\centering 
				\citep{esfahani2018data, gao2016distributionally, shafieezadeh2023new}
				}&
				\multirow{2}{2.6cm}{\centering Piecewise concave in $z$} & \multirow{2}{2cm}{\centering General}  & \multirow{2}{3.7cm}{\centering Empirical distribution} &\multirow{2}{3.7cm}{\centering General} \\
				& &   &   & \\
				\cline{1-5}
				\multirow{2}{2cm}{\centering 
				\citep{Shafieezadeh15, yu2022fast, selvi2022wasserstein, shafieezadeh2023new}
				}&
				\multirow{2}{2.8cm}{\centering Generalized linear model in $(z, \theta)$} & 
				\multirow{2}{2cm}{\centering General} & \multirow{2}{3.7cm}{\centering General} &
				\multirow{2}{3.5cm}{\centering Whole Euclidean space\footnote{
				Here the references essentially assume the numerical part of the probability vector is supported on the whole Euclidean space, such as the numerical features 				is supported on the entire Euclidean space.
				}} \\
				& &   &  & \\
				\cline{1-5}
    \multirow{2}{2cm}{\centering 
				\citep{sinha2018certifiable, blanchet2021optimal}
				}& 
				\multirow{2}{3.4cm}{\centering $z\mapsto f_{\theta}(z) - \lambda^*c(x,z)$ is strongly concave\footnote{
				\citet{sinha2018certifiable} approximately solves the Wasserstein DRO by penalizing the Wasserstein ball constraint with fixed Lagrangian multiplier $\lambda^*$.
				Here the assumption of loss function holds for $\hP$-almost every $x$.
				}} &  
				\multirow{2}{2.5cm}{\centering Strongly convex function\footnote{
We say a transport cost $c(\cdot,\cdot)$ is strongly convex if $c(x,y)=u(x-y)$ for a strongly convex function $u(\cdot)$.
                }} & 
				\multirow{2}{3.7cm}{\centering General} & 
				\multirow{2}{3.7cm}{\centering General}\\
				& &  &   & \\
				\bottomrule
				\end{tabular}%
		\end{threeparttable}
	\end{table}

To address the aforementioned concerns while retaining the advantages of Wasserstein DRO, we propose a novel approach called Sinkhorn DRO. 
Sinkhorn DRO leverages the Sinkhorn distance \citep{cuturi2013sinkhorn}, which hedges against distributions that are close to a given nominal distribution in Sinkhorn distance. The Sinkhorn distance can be viewed as a smoothed version of the Wasserstein distance and is defined as the minimum transport cost between two distributions associated with an optimal transport problem with entropic regularization (see Definition \ref{EntW:Wasserstein} in Section \ref{Sec:model}).
To the best of our knowledge, this paper is the first to explore the DRO formulation using the Sinkhorn distance. Our work makes several contributions, which are summarized below:
\begin{enumerate}
    \item 
We derive a strong duality reformulation for Sinkhorn DRO (Theorem~\ref{Theorem:strong:duality}) in a highly general setting, where the loss function, transport cost, nominal distribution, and probability support are allowed to be arbitrary.
The dual objective of Sinkhorn DRO smooths the dual objective of Wasserstein DRO, where the level of smoothness is controlled by the entropic regularization parameter (Remark~\ref{Remark:Wasserstein:DRO}).
    \item
Our duality proof yields an insightful characterization of the worst-case distribution~(Remark~\ref{Remark:worst:Sinkhorn:DRO}). 
Unlike Wasserstein DRO, where the worst-case distribution is typically discrete and finitely supported, the worst-case distribution in Sinkhorn DRO is absolutely continuous with respect to a pre-specified reference measure, such as the Lebesgue or counting measure.
This characteristic of Sinkhorn DRO highlights its flexibility as a modeling choice and provides a more realistic representation of uncertainty that better aligns with the underlying true distribution in practical scenarios.

    \item
\Jie{The dual reformulation of Sinkhorn DRO can be viewed as a conditional stochastic optimization~\citep{hu2020sample, Yifan20, hu2021biasvar} involving an expectation (with respect to observed samples) of nonlinear transformation of a conditional expectation~(with respect to a conditional distribution).
In our work, we introduce and analyze an efficient stochastic mirror descent algorithm with biased subgradient estimators to solve this problem (Section~\ref{Sec:first:order}).
We quantify the computational cost using the number of generated samples from the outer expectation and the number of generated samples from the inner expectation.
Our algorithm achieves both complexities of $\tO(\delta^{-2})$ for a fixed entropic regularization parameter $\Reg$.%
}
    \item
To validate the effectiveness and efficiency of the proposed Sinkhorn DRO model, we conduct a series of experiments in Section~\ref{Sec:numerical}, including the newsvendor problem, mean-risk portfolio optimization, and multi-class adversarial classification.
Using synthetic and real datasets, we compare the Sinkhorn DRO model against benchmarks such as SAA, Wasserstein DRO, and KL-divergence DRO.
The results demonstrate that the Sinkhorn DRO model consistently outperforms the benchmarks in terms of out-of-sample performance and computational speed.
\Jie{
We also provide a comprehensive set of experiment studies to show that there exists a large number of parameter choices under which Sinkhorn DRO outperforms Wasserstein DRO.
}
\end{enumerate}

\subsection*{Related Literature}\label{Sec:literature}

In the following, we first compare our work with the four most closely related papers that appear recently.

\Jie{
\citet{feng2018model} studied the Wasserstein DRO formulation with an additional differential entropy constraint on the optimal transport mapping, which
is closely related to our Sinkhorn DRO formulation. } They derived a weak dual formulation and characterized the worst-case distribution under the assumption that strong duality holds. 
It is important to note that such an assumption cannot be taken for granted for the considered infinite-dimensional problem. 
\Jie{Instead, we provided a rigorous proof of strong duality for our Sinkhorn DRO formulation.}
Moreover, their results heavily depend on the assumption that the nominal distribution $\hP$ is absolutely continuous with respect to the Lebesgue measure.
This limits the applicability of their formulation in data-driven settings where $\hP$ is discrete.
Since the initial submission of our work, \citet{azizian2022regularization} have presented a duality result similar to ours, but with different assumptions. Their results apply to more general regularization beyond entropic regularization, but they assume a continuous loss function and a compact probability space under the Slater condition. 
\citet{song2022efficient} have recently explored the application of Sinkhorn DRO in reinforcement learning. Their duality proof rely on the boundedness of the loss function and the discreteness of the probability support. 
These three papers do not present numerical algorithms to solve the dual formulation.
\citet[Section~3.2]{blanchet2020semi} solved a log-sum-exp approximation of the Wasserstein DRO dual formulation. This smooth approximation can be viewed as a special case of the dual reformulation of our Sinkhorn DRO model.
However, their study did not specifically explore the primal form of Sinkhorn DRO. 
Their algorithm employed unbiased subgradient estimators, even though the second-order moment could be unbounded. The paper did not provide explicit theoretical convergence guarantees for their algorithm. Additionally, numerical comparisons detailed in Appendix \ref{Appendix:compare:opt:alg} suggest that our proposed algorithm outperforms theirs in terms of empirical performance.

Next, we review papers on several related topics. 

\textit{On DRO models.}
In the literature on DRO, there are two main approaches to constructing ambiguity sets.
The first approach involves defining ambiguity sets based on descriptive statistics, such as support information \citep{bertsimas2006persistence}, moment conditions \citep{scarf1958min, Delage10, Goh10, Zymler13, wiesemann2014distributionally, Chen19,Bertsimas19}, shape constraints \citep{popescu2005semidefinite,van2015generalized}, marginal distributions \citep{frechet1960tableaux,natarajan2009persistency,agrawal2012price,doan2012complexity}, etc.
The second approach, which has gained popularity in recent years, involves considering distributions within a pre-specified statistical distance from a nominal distribution.
Commonly used statistical distances in the literature include $\phi$-divergence \citep{hu2013kullback, Ben13,wang2016likelihood,bayraksan2015data, Duchi21}, Wasserstein distance \citep{pflug2007ambiguity,wozabal2012framework,Mohajerin18,zhao2018data,blanchet2019quantifying, gao2016distributionally,chen2018data, xie2019distributionally}, and maximum mean discrepancy \citep{staib2019distributionally, Kernelzhu}.
Our proposed Sinkhorn DRO can be seen as a variant of the Wasserstein DRO.
In the literature on Wasserstein DRO,
researchers have also explored the regularization effects and statistical inference of the approach. 
In particular, it has been shown that Wasserstein DRO is asymptotically equivalent to a statistical learning problem with variation regularization \citep{gao2020wasserstein,blanchet2019robust,shafieezadeh2019regularization}. When the radius is chosen properly, the worst-case loss of Wasserstein DRO serves as an upper confidence bound on the true loss~\citep{blanchet2019robust,blanchet2019confidence,gao2020finitesample,blanchet2021statistical}.
Variants of Wasserstein DRO have been proposed by combining it with other information, such as moment information
\citep{Wang18} or marginal distributions \citep{eckstein2020robust} to enhance its modeling capabilities.

\textit{On Sinkhorn distance.}
Sinkhorn distance \citep{cuturi2013sinkhorn} was proposed to improve the computational complexity of Wasserstein distance, by regularizing the original mass transportation problem with relative entropy penalty on the transport plan.
It has been demonstrated to be beneficial because of lower computational cost in various applications, including domain adaptations~\citep{courty2014domain2, courty2016optimal,courty2017joint}, generative modeling~\citep{Aude18a,petzka2018on,luise2018differential,patrini2020sinkhorn}, dimension reduction~\citep{lin2020projection2, wang2020kerneltwosample, Riemannianhuang21}, etc.
In particular, this distance can be computed from its dual form by optimizing two blocks of decision variables alternatively, which only requires simple matrix-vector products and therefore significantly improves the computation speed~\citep{Computational19, mensch2020online, lin2022efficiency, Altschuler17}.
Such an approach first arises in economics and survey statistics~\citep{kruithof1937telefoonverkeersrekening, yule1912methods, deming1940least, bacharach1965estimating}, and its convergence analysis is attributed to the mathematician Sinkhorn~\citep{sinkhorn1964relationship}, which gives the name of Sinkhorn distance.
\Jie{Computing the Sinkhorn distance between a discrete distribution and an arbitrary distribution can be reformulated as a stochastic optimization problem with a log-sum-exp–type loss function~\citep[Section~5.4]{Computational19}.
For Sinkhorn DRO, the dual objective takes the form of an expectation involving the logarithm of a conditional expectation of an exponential function.
When the inner expectation is over a discrete distribution, the problem retains a similar structure and can be effectively addressed using standard stochastic optimization techniques.
In the case of general (non-discrete) distributions, the formulation becomes more challenging to solve, primarily due to the difficulty of obtaining unbiased (sub)gradient estimators.}

\textit{On algorithms for solving DRO models.}
In the introduction, we have elaborated on the literature that proposes efficient optimization algorithms for solving the Wasserstein DRO dual formulation~\citep{zhao2018data, chen2021decomposition, liu2021discrete,sinha2018certifiable,esfahani2018data, gao2016distributionally,Shafieezadeh15, yu2022fast, selvi2022wasserstein,blanchet2021optimal,shafieezadeh2023new}, in which the computational efficiency is limited to a certain class of loss functions, transport costs, and nominal distributions.
To solve the $\phi$-divergence DRO, one common approach is to employ sample average approximation (SAA) to approximate the dual formulation. However, SAA requires storing the entire set of samples, making it inefficient in terms of storage usage.
An alternative approach is to use first-order stochastic subgradient algorithms, which are more storage-efficient. These algorithms have the advantage of complexity that can be independent of the sample size of the nominal distribution \citep{levy2020large, namkoong2016stochastic, qi2022stochastic}.
Our derived dual reformulation of Sinkhorn DRO can be seen as an instance of the CSO problem \cite{hu2020sample, Yifan20, hu2021biasvar}.
In this context, we have developed stochastic mirror descent algorithms with biased subgradient oracles.
Notably, our proof can be adjusted to show that the proposed algorithm achieves near-optimal complexity for general CSO problems with both smooth and nonsmooth loss functions, marking an improvement over the state-of-the-art~\citep[Theorem~3.2]{Yifan20} that is sub-optimal for nonsmooth loss functions.

\medskip

The rest of the paper is organized as follows. In Section~\ref{Sec:model}, we describe the main formulation for the Sinkhorn DRO model. In Section~\ref{Sec:Reformulation:DRO}, we develop its strong dual reformulation.
In Section~\ref{Sec:first:order}, we propose a first-order optimization algorithm that solves the reformulation efficiently.
We report several numerical results in Section~\ref{Sec:numerical}, and conclude the paper in
Section~\ref{Sec:conclusion}.
All omitted proofs can be found in Appendices.

\section{Model Setup}\label{Sec:model}

\noindent\textit{Notation.}
Assume the logarithm function $\log$ is taken with base $e$. For a positive integer $N$, we write $[N]$ for $\{1,2,\ldots,N\}$.
For a measurable set $\cZ$, denote by $\cM(\cZ)$ the set of measures (not necessarily probability measures) on $\cZ$, and $\cP(\cZ)$ the set of probability measures on $\cZ$.
Given a probability distribution $\bP$ and a measure $\mu$, we denote $\mathrm{supp}\, \bP$ the support of $\bP$, and write $\bP\ll \mu$ if $\bP$ is absolutely continuous with respect to $\mu$.
Given a measure $\mu\in\cM(\cZ)$ and a measurable variable $f:~\cZ\to\mathbb{R}$, we write $\bE_{z\sim\mu}[f(z)]$ for $\int f(z)\diff\mu(z)$.
For a given element $x$, denote by $\delta_x$ the one-point probability distribution supported on $\{x\}$.
Denote $\bP\otimes\bQ$ as the product measure of two probability measures $\bP$ and $\bQ$.
Denote by $\proj_{1\#}\gamma$ and $\proj_{2\#}\gamma$ the first and the second marginal distributions of $\gamma$, respectively.
For a function $\omega:~\Theta\to\mathbb{R}$, we say it is $\kappa$-strongly convex with respect to norm $\|\cdot\|$ if $\inp{\theta' - \theta}{\nabla \omega(\theta') - \nabla\omega(\theta)}\ge \kappa\|\theta' - \theta\|^2, \forall \theta,\theta'\in\Theta$.

We first review the definition of Sinkhorn distance.
\begin{definition}[Sinkhorn Distance]
\label{EntW:Wasserstein}
Let $\cZ$ be a measurable set.
Consider distributions $\bP,\bQ\in\cP(\cZ)$, and let $\mu,\nu\in\cM(\cZ)$ be two reference measures such that $\bP\ll\mu$, $\bQ\ll\nu$.
For regularization parameter $\epsilon \ge 0$, the \emph{Sinkhorn distance} between two distributions $\bP$ and $\bQ$ is defined as
\[\label{Eq:Sinkhorn:primal}
\mathcal{W}_{\Reg}(\bP, \bQ) \:= \inf_{
\substack{
\gamma\in\Gamma(\bP, \bQ)
}}~ \left\{\mathbb{E}_{(x,y)\sim\gamma}[c(x,y)] + \Reg H(\gamma\mid\mu\otimes\nu)\right\},
\]
where $\Gamma(\bP,\bQ)$ denotes the set of joint distributions whose first and second marginal distributions are $\bP$ and $\bQ$ respectively, $c(x,y)$ denotes the transport cost, and $H(\gamma\mid\mu\otimes\nu)$ denotes the relative entropy of $\gamma$ with respect to the product measure $\mu\otimes\nu$:
\[
H(\gamma\mid\mu\otimes\nu)\:= \bE_{(x,y)\sim\gamma}\left[\log\left( 
\frac{\diff\gamma(x,y)}{\diff\mu(x)\diff\nu(y)}
\right)\right],
\]
where $\frac{\diff\gamma(x,y)}{\diff\mu(x)\diff\nu(y)}$ stands for the density ratio of $\gamma$ with respect to $\mu\otimes\nu$ evaluated at $(x,y)$.
\QDEF
\end{definition}

\begin{remark}[Variants of Sinkhorn Distance]
Sinkhorn distance in Definition~\ref{EntW:Wasserstein} is based on general reference measures $\mu$ and $\nu$.
Special forms of distance have been investigated in the literature. For instance, the entropic regularized optimal transport distance $\mathcal{W}_{\Reg}^{\texttt{Ent}}(\bP, \bQ)$ \citep[Equation~(2)]{cuturi2013sinkhorn} chooses $\mu$ and $\nu$ as the Lebesgue measure when the corresponding $\bP$ and $\bQ$ are continuous, or counting measures if $\bP$ and $\bQ$ are discrete.
For given $\bP$ and $\bQ$, 
one can check the two distances above are equivalent up to a constant:
\begin{align*}
\mathcal{W}_{\Reg}^{\texttt{Ent}}(\bP, \bQ) &=
\mathcal{W}_{\Reg}(\bP, \bQ) + \bE_{(x,y)\sim\gamma}\left[\log\left(
\frac{\diff\mu(x)\diff\nu(y)}{\diff x\diff y}
\right)\right]\\
&=\mathcal{W}_{\Reg}(\bP, \bQ)  + 
\bE_{x\sim\bP}\left[ 
\log\left(
\frac{\diff\mu(x)}{\diff x}
\right)
\right]
+\bE_{y\sim\bQ}\left[ \log\left(
\frac{\diff\nu(y)}{\diff y}
\right)\right].
\end{align*}
Another variant is to chose $\mu$ and $\nu$ to be $\bP,\bQ$, 
respectively \citep[Section~2]{aude2016stochastic}.
A hard-constrained variant of the relative entropy regularization has been discussed in \citep[Definition~1]{cuturi2013sinkhorn} and \citep{bai2020information}:
\[
\mathcal{W}_{R}^{\texttt{Info}}(\bP, \bQ)
\:= \inf_{\gamma\in\Gamma(\bP, \bQ)}\left\{
\mathbb{E}_{(X,Y)\sim\gamma}[c(X,Y)]:\ H(\gamma\mid\bP\otimes\bQ)\le R
\right\},
\]
where $R\ge0$ quantifies the upper bound for the relative entropy between distributions $\gamma$ and $\bP\otimes\bQ$.
\QEG
\end{remark}

\begin{remark}[Choice of Reference Measures]\label{remark:choice:reference:measure}
We discuss below our choice of the two reference measures $\mu$ and $\nu$ in Definition \ref{EntW:Wasserstein}.
For the reference measure $\mu$, observe from the definition of relative entropy and the law of probability, we can see that the regularization term in $\W_{\Reg}(\hP, \bP)$ can be written as
\begin{align*}
H(\gamma\mid\mu\otimes\nu)&\:= \bE_{(x,y)\sim\gamma}\left[\log\left( 
\frac{\diff\gamma(x,y)}{\diff\hP(x)\diff\nu(y)}
\right) +
\log\left( 
\frac{\hP(x)}{\diff\mu(x)}
\right)\right]\\
&=\bE_{(x,y)\sim\gamma}\left[\log\left( 
\frac{\diff\gamma(x,y)}{\diff\hP(x)\diff\nu(y)}
\right)\right] + \bE_{x\sim\hP}\left[\log\left( 
\frac{\hP(x)}{\diff\mu(x)}
\right)\right].
\end{align*}
Therefore, any choice of the reference measure $\mu$ satisfying $\hP\ll\mu$ is equivalent up to a constant. 
For simplicity, in the sequel we will take $\mu=\hP$.
For the reference measure $\nu$, observe that the worst-case solution $\bP$ in \eqref{Inf-W-E} should satisfy that $\bP\ll\nu$ since otherwise the entropic regularization in Definition~\ref{EntW:Wasserstein} is undefined.
As a consequence, we can
choose $\nu$ which the underlying true distribution is absolutely continuous with respect to and is easy to sample from.
For example, if we believe the underlying distribution is continuous, then we can choose $\nu$ to be the Lebesgue measure or Gaussian measure, or if we believe the underlying distribution is discrete, we can choose $\nu$ to be a counting measure.
We refer to \citep[Section~3.6]{pichler2021mathematical} for the construction of a general reference measure.
\QEG
\end{remark}

In this paper, we study the Sinkhorn DRO model.
Given a loss function $f$, a nominal distribution $\hP$ and the Sinkhorn radius $\rho$, the primal form of the worst-case expectation problem of Sinkhorn DRO is given by
\begin{equation}
\tag{\texttt{Primal}}\label{Inf-W-E}
\begin{aligned}
&V_{\Primal}
:=\sup_{\mathbb{P}\in\mathbb{B}_{\rho,\Reg}(\hP)}~
\mathbb{E}_{z\sim \mathbb{P}}[f(z)],\\
\end{aligned}
\end{equation}
where $\mathbb{B}_{\rho,\Reg}(\hP):= \big\{
\mathbb{P}:\,\mathcal{W}_{\Reg}(\hP, \mathbb{P})\le \rho
\big\}$ is the Sinkhorn ball of the radius $\rho$ centered at the nominal distribution $\hP$.
Due to the convex entropic regularization in $\mathcal{W}_{\Reg}(\hP, \mathbb{P})$~\citep{Cover06}, the Sinkhorn distance $\mathcal{W}_{\Reg}(\hP, \mathbb{P})$ is convex in $\bP$, i.e., it holds that
$
\mathcal{W}_{\Reg}(\hP, \lambda\bP_1 + (1-\lambda)\bP_2)\le 
\lambda \mathcal{W}_{\Reg}(\hP, \bP_1) + (1-\lambda)\mathcal{W}_{\Reg}(\hP, \bP_2)
$
for all probability distributions $\bP_1$ and $\bP_2$ and all $0\le \lambda\le 1$.
Therefore, the Sinkhorn ball is a convex set, and the problem \eqref{Inf-W-E} is an (infinite-dimensional) convex program.

Our goal for the rest of the paper is to derive the dual reformulation and efficient algorithms to solve the Sinkhorn DRO model.

\section{Strong Duality Reformulation}\label{Sec:Reformulation:DRO}
Problem~\eqref{Inf-W-E} is an infinite-dimensional optimization problem over probability distributions. To obtain a more tractable form, in this section, we derive a strong duality result for \eqref{Inf-W-E}.
Our main goal is to derive the strong dual program
\begin{equation}\label{Eq:dual:SDRO:general}
  V_{\Dual} := \inf_{\lambda\ge0}\;\bigg\{\lambda\rho + \lambda\Reg\; \mathbb{E}_{x\sim\hP}\Big[ \log\mathbb{E}_{z\sim\nu}\big[
   e^{
  (f(z) - \lambda c(x,z))/(\lambda\Reg)}
  \big] \Big]\bigg\},
\end{equation}
where the dual variable $\lambda$ corresponds to the Sinkhorn ball constraint in \eqref{Inf-W-E}, and by convention, we define the dual objective evaluated at $\lambda=0$ as the limit of the objective values with $\lambda\downarrow0$, which equals the essential supremum of the objective function with respect to the measure $\nu$.
Or equivalently, by defining the constant
\begin{equation}\label{Eq:def:orho:general:unified}
\orho := \rho+\Reg\,\mathbb{E}_{x\sim\hP}\Big[\log\mathbb{E}_{z\sim\nu}\big[ e^{-c(x,z)/\Reg}\big]\Big],
\end{equation}
and the kernel probability distribution
\begin{equation}
    \label{Eq:bQ:x:Gibbs}
\diff\bQ_{x,\Reg}(z) := \frac{e^{-c(x,z)/\Reg}}{
\bE_{u\sim\nu}\left[e^{-c(x,u)/\Reg}\right]
}\diff\nu(z),
\end{equation}
we have
\begin{equation}
\begin{aligned}\label{Eq:dual:SDRO:general:unified}
V_{\Dual} =
\inf_{\lambda\ge0}\;\bigg\{
\lambda\orho+\lambda\Reg\,\mathbb{E}_{x\sim\hP}\Big[\log
\mathbb{E}_{z\sim\bQ_{x,\Reg}}\big[ 
e^{f(z)/(\lambda\Reg)}
\big]\Big]
\bigg\}.
\end{aligned}
\tag{\texttt{Dual}}
\end{equation}

The rest of this section is organized as follows.
In Section~\ref{sec:main:result}, we summarize our main results on the strong duality reformulation of Sinkhorn DRO.
Next, we provide detailed examples in Section~\ref{Sec:example:reformulation} and discussions in Section~\ref{Sec:ent:duality}.
In Section~\ref{Proof:Theorem:strong:duality:main}, we provide a proof sketch of our main results.

\subsection{Main Theorem}
\label{sec:main:result}

To make the above primal \eqref{Inf-W-E} and dual \eqref{Eq:dual:SDRO:general:unified} problems well-defined, we introduce the following assumptions on the transport cost $c$, the reference measure $\nu$, and the loss function $f$.
\begin{assumption}\label{Assumption:distance:measure:function}
\begin{enumerate}
    \item\label{Assumption:distance:measure:function:1} 
    \Jie{
    The cost function $c(x,z)$ is $\hP\otimes\nu$-measurable and satisfies
    }
    $\nu\{z:~0\le c(x,z)<\infty\}=1$ for $\hP$-almost every $x$; %
    \item\label{Assumption:distance:measure:function:2} $
    \bE_{z\sim\nu}\left[e^{-c(x,z)/\Reg}\right]<\infty$ for $\hP$-almost every $x$;
    \item\label{Assumption:distance:measure:function:3} $\cZ$ is a measurable space, and the function $f:~\cZ\to\mathbb{R}\cup\{\infty\}$ is measurable.

    \item\label{Assumption:distance:measure:function:4}
    For every joint distribution $\gamma$ on $\cZ\times\cZ$ with first marginal distribution $\hP$, it has a regular conditional distribution $\gamma_x$ given the value of the first marginal equals $x$.
\end{enumerate}
\end{assumption}

Assumption~\ref{Assumption:distance:measure:function}\ref{Assumption:distance:measure:function:1} implies that $0\le c(x,y)<\infty$ for $\hP\otimes\nu$-almost every $(x,y)$.
By \citep[Proposition~4.1]{Computational19}, the Sinkhorn distance has an equivalent formulation
\[
\mathcal{W}_{\Reg}(\hP, \bP)=
\min_{\gamma\in\Gamma(\hP, \bP)}~
\int \log\left( 
\frac{\diff\gamma}{\diff\mathcal{K}}(x,y)
\right)\diff\gamma(x,y),\quad 
\text{where }\diff\mathcal{K}(x,y) = e^{-c(x,y)/\Reg}\diff\hP(x)\diff\nu(y).
\]
Therefore Assumption~\ref{Assumption:distance:measure:function}\ref{Assumption:distance:measure:function:1} ensures that the reference measure $\mathcal{K}$ is well-defined.
Assumption~\ref{Assumption:distance:measure:function}\ref{Assumption:distance:measure:function:2} ensures the optimal transport mapping $\gamma_*$ for Sinkhorn distance $\mathcal{W}_{\Reg}(\hP, \bP)$ exists with density value $\frac{\diff\gamma_*(x,y)}{\diff\hP(x)\diff\nu(y)}\propto e^{-c(x,y)/\Reg}$.
Hence, Assumptions~\ref{Assumption:distance:measure:function}\ref{Assumption:distance:measure:function:1} and \ref{Assumption:distance:measure:function}\ref{Assumption:distance:measure:function:2} together ensure the Sinkhorn distance is well-defined. 
Assumption~\ref{Assumption:distance:measure:function}\ref{Assumption:distance:measure:function:3} ensures the expected loss $\mathbb{E}_{z\sim\bP}[f(z)]$ to be well-defined 
for any distribution $\bP$.
Assumption~\ref{Assumption:distance:measure:function}\ref{Assumption:distance:measure:function:4} ensures the joint distribution $\gamma$ can be written as $\diff\gamma(x,z) = \diff\hP(x)\diff\gamma_x(z)$ and the law of total expectation holds; we refer to \cite[Chapter 5]{kallenberg1997foundations} for the concept of the regular conditional distribution.
\Jie{Such an assumption is very mild; for instance, it holds if $\cZ$ is a Polish space~\citep{blackwell1963non}.}

To distinguish the cases $V_D<\infty$ and $V_D=\infty$, we introduce the light-tail condition on $f$ in Condition~\ref{Assumption:light:tailed:f}.
In Appendix \ref{app:sufficient condition}, we present sufficient conditions for Condition~\ref{Assumption:light:tailed:f} that are easy to verify.
\begin{condition}\label{Assumption:light:tailed:f}
There exists $\lambda>0$ such that $\mathbb{E}_{z\sim\bQ_{x,\Reg}}[ e^{f(z)/(\lambda\Reg)}]<\infty$ for $\hP$-almost every $x$.
\end{condition}

In the following, we provide the main results of the strong duality reformulation.
\begin{theorem}[Strong Duality]\label{Theorem:strong:duality}
Let $\hP\in\cP(\cZ)$, and assume Assumption \ref{Assumption:distance:measure:function} holds. Then the following holds:
\begin{enumerate}
\item\label{Theorem:strong:duality:I}
The primal problem \eqref{Inf-W-E} is feasible if and only if $\orho\ge0$;
\item\label{Theorem:strong:duality:II}
  Whenever $\orho\ge0$, it holds that $V=V_D$.
\item\label{Theorem:strong:duality:III}
If, in addition, Condition \ref{Assumption:light:tailed:f} holds and $\orho>0$, it holds that
$V_{\Primal} = V_{\Dual} <\infty$; otherwise $V_{\Primal} = V_{\Dual} = \infty$.
\Jie{
\item\label{Theorem:strong:duality:IV}
Assume in addition that Condition \ref{Assumption:light:tailed:f} holds and $\orho>0$.
Define the event $A:=\{z:~f(z)=\esssup\limits_{\nu}(f)\}$ with $\esssup\limits_{\nu}(f):=\inf\{t:~\nu\{f(z)>t\}=0\}$.
The dual minimizer $\lambda^*=0$ if and only if $\esssup\limits_{\nu}(f)<\infty$ and $\orho\ge \Reg\bE_{x\sim\hP}[\log(1/\bP_{z\sim\bQ_{x,\Reg}}\{A\})]$.
}
\end{enumerate}
\end{theorem}

We remark that if $\orho<0$,  by convention, $V=-\infty$ and $V_D=-\infty$ as well by Lemma \ref{Lemma:finite:dual} in Appendix~\ref{Sec:appendix:proof}. Therefore, we have $V=V_D$ as long as Assumption \ref{Assumption:distance:measure:function} holds.

\subsection{Discussions}\label{Sec:ent:duality}
In the following, we make several remarks regarding the strong duality result.

\begin{remark}[Comparison with Wasserstein DRO]\label{Remark:Wasserstein:DRO}
As the regularization parameter $\Reg\to0$, the dual objective of the Sinkhorn DRO \eqref{Eq:dual:SDRO:general:unified} converges to (see Appendix~\ref{Sec:proof:sec:ent:duality} for details)
\[
\lambda\rho + \mathbb{E}_{x\sim\hP}\Big[
\sup_{z\in\mathrm{supp}\,\nu}
\big\{f(z) - \lambda c(x,z)\big\}\Big],
\]
which essentially follows from the fact that the log-sum-exp function is a smooth approximation of the supremum.
Particularly, when $\mathrm{supp}\,\nu=\cZ$, the dual objective of the Sinkhorn DRO converges to the dual objective of the Wasserstein DRO \citep[Theorem~1]{gao2016distributionally}.
  The main computational difficulty in Wasserstein DRO is solving the maximization problem inside the expectation above.
  All results in Table \ref{Tab:compare:DRO} ensure this inner maximization can be efficiently solved.
  \Jie{
  On the one hand, as Sinkhorn DRO does not need to solve this maximization, it does not need stringent assumptions on $f(\cdot)$ and thus enables efficient implementation for a larger class of loss functions for a fixed regularization parameter $\Reg$~(see detailed discussion in Section~\ref{Sec:first:order}).
  On the other hand, when the stringent assumptions on $f(\cdot)$ are satisfied, (data-driven) Wasserstein DRO typically admits a finite-dimensional convex reformulation, which can be solved more efficiently than Sinkhorn DRO (see our numerical comparison of CPU times in Appendix~\ref{Sec:baseline:run}). The key reason is that, in these special cases, Wasserstein DRO yields a finitely supported worst-case distribution, whereas Sinkhorn DRO always results in a worst-case distribution supported over the entire sample space.
  }

We also remark that Sinkorn DRO and Wasserstein DRO result in different conditions for finite worst-case values. 
From Condition~\ref{Assumption:light:tailed:f} we see that Sinkhorn DRO is finite if and only if under a light-tail condition on $f$, while Wasserstein DRO is finite if and only if the loss function satisfies a growth condition \citep[Theorem~1 and Proposition~2]{gao2016distributionally}: $f(z)\le L_fc(z,z_0)+M, \forall z\in\cZ$ for some constants $L_f,M>0$ and $z_0\in\cZ$. 
\QEG
\end{remark}

\begin{remark}[Worst-case Distribution]\label{Remark:worst:Sinkhorn:DRO}
\Jie{
Assume $\orho>0$ and Condition~\ref{Assumption:light:tailed:f} holds, and there exists (actually, Lemma~\ref{Lemma:lambda:positive} ensures its uniqueness) an optimal Lagrangian multiplier $\lambda^*>0$ in \eqref{Eq:dual:SDRO:general:unified}.
}
As we will demonstrate in the proof of Theorem~\ref{Theorem:strong:duality}, the worst-case distribution maps every $x\in\mathrm{supp}\,\hP$ to a (conditional) distribution \Jie{$\gamma^*_x$ that solves a strictly convex program~(i.e., Problem~\eqref{Eq:v:x:lambda}),} whose density function (with respect to $\nu$) is
\[
  \alpha_x\cdot
  \exp\Big(\big( 
  f(z) - \lambda^\ast c(x,z) \big)/(\lambda^\ast\Reg)\Big),
\]
where $\alpha_x:=\left(
\bE_{z\sim\nu}\left[e^{(f(z) - \lambda^{\ast}c(x,z))/(\lambda^{\ast}\Reg)}\right]
  \right)^{-1}$ is a normalizing constant to ensure the conditional distribution well-defined.
  \Jie{The uniqueness of $\gamma^*_x, x\in\mathrm{supp}\,\hP$ ensures the uniqueness of the worst-case distribution $\bP_{\ast}$}, whose density becomes 
\[
\frac{\diff\bP_{\ast}(z)}{\diff\nu(z)}=
\bE_{x\sim\hP}\left[
\alpha_x\cdot \exp\Big(\big( 
  f(z) - \lambda^\ast c(x,z) \big)/(\lambda^\ast\Reg)\Big)
\right].
\]
As such, the worst-case distribution shares the same support as the measure $\nu$.

Particularly, when $\hP$ is the empirical distribution $\frac{1}{n}\sum_{i=1}^n \delta_{\hx_i}$ and $\nu$ is any continuous distribution on $\mathbb{R}^d$, the worst-case distribution $\bP_\ast$ is supported on the entire $\mathbb{R}^d$.
In contrast, the worst-case distribution for Wasserstein DRO is supported on at most $n+1$ points \cite{gao2016distributionally}.
In Fig.~\ref{fig:DRO:transport} we visualize the worst-case distributions from Wasserstein/Sinkhorn DRO models.
The loss function and transport cost used in this plot follow the setup described in Example~\ref{Example:DRO:LR}.
The Wasserstein ball radius, Sinkhorn ball radius, and entropic regularization value are fine-tuned to ensure that the optimal dual multipliers for all instances equal 5.
Notably, the support points of the worst-case distributions from the Wasserstein DRO model correspond to the modes of the continuous worst-case distributions from the Sinkhorn DRO model.
\begin{figure}[H]
    \centering
    \includegraphics[width=0.24\textwidth]{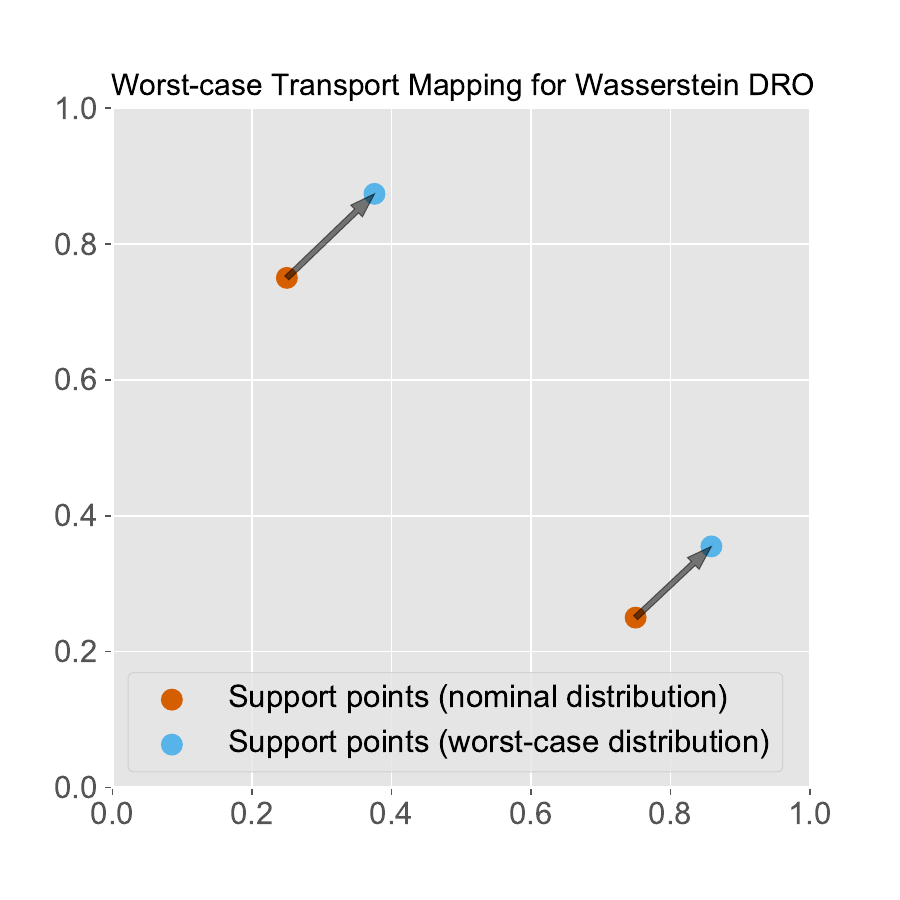}
    \includegraphics[width=0.24\textwidth]{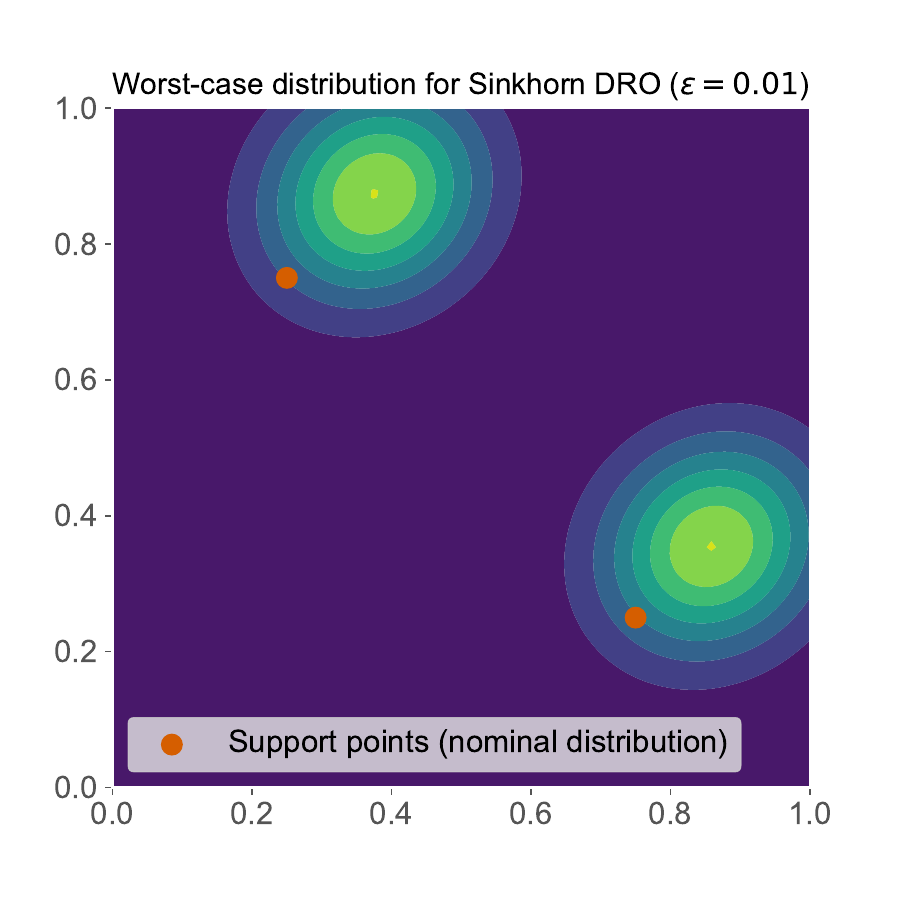}
    \includegraphics[width=0.24\textwidth]{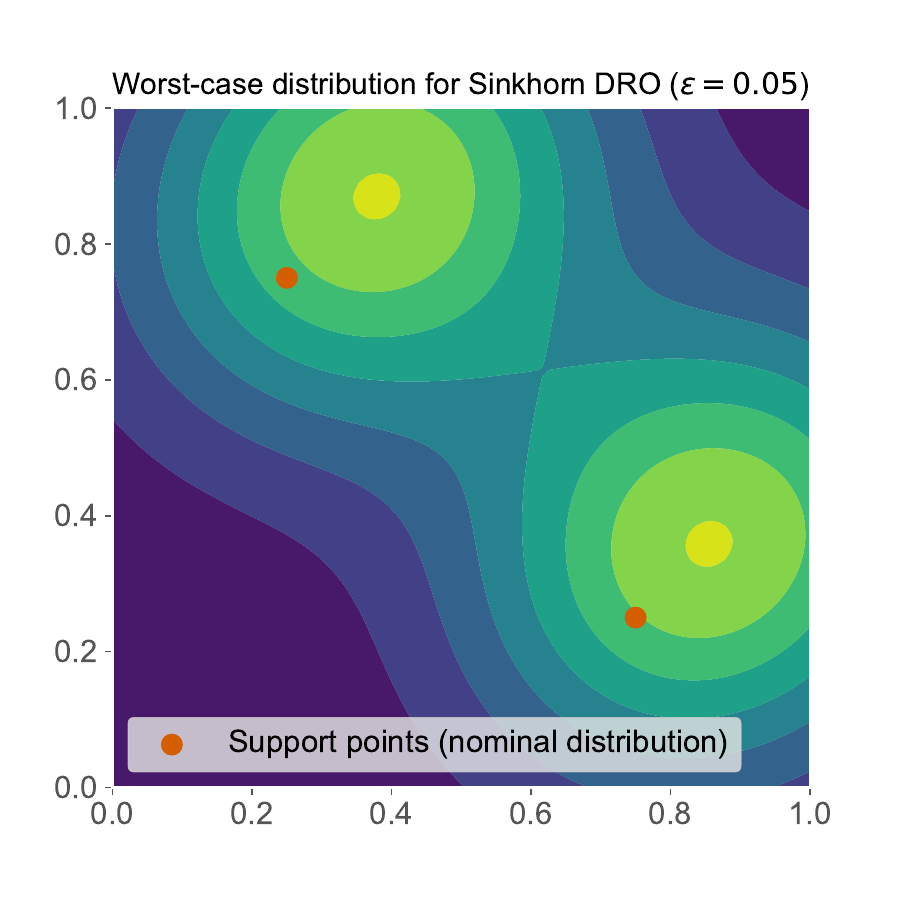}
    \includegraphics[width=0.24\textwidth]{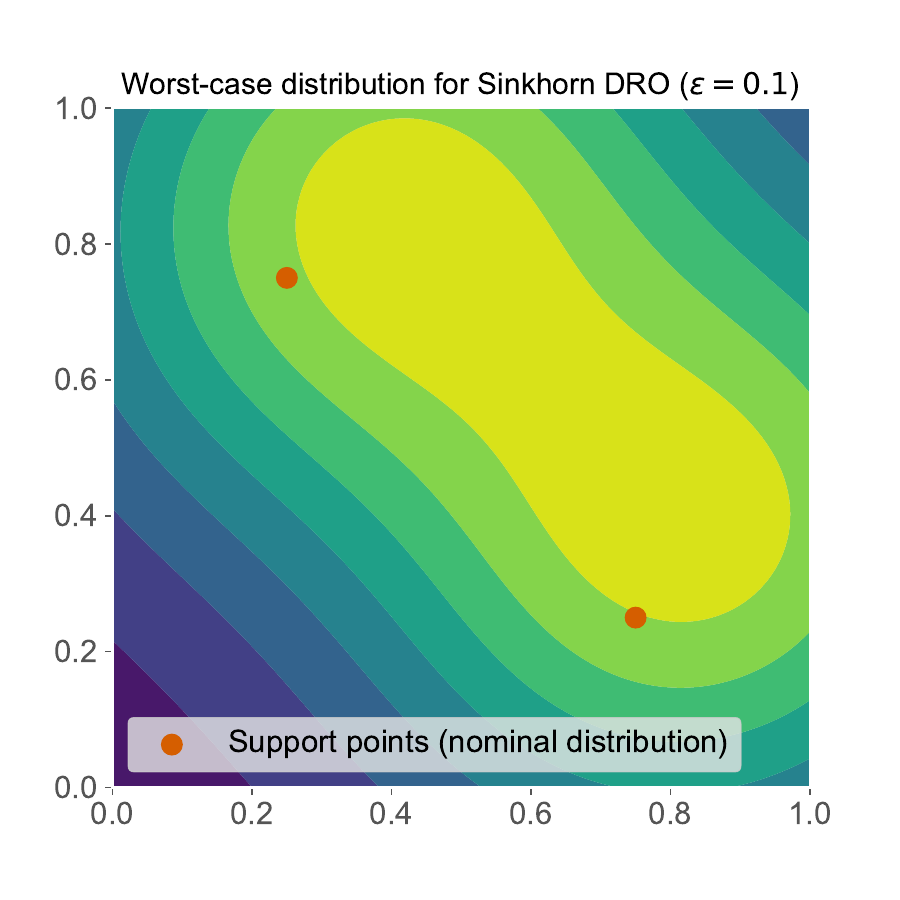}
    \caption{
    Visualization of worst-case distributions from Wasserstein DRO (left plot) and Sinkhorn DRO models (right three plots) with varying choices of $\Reg$.%
    }
    \label{fig:DRO:transport}
\end{figure}
The above demonstrates another difference~(perhaps advantage) of Sinkhorn DRO compared with Wasserstein DRO.
Indeed, for many practical problems, the underlying distribution is modeled as a continuous distribution. The worst-case distribution for Wasserstein DRO is often finitely supported, raising the concern of whether it hedges against the wrong family of distributions and thus results in suboptimal solutions. 
The numerical results in Section \ref{Sec:numerical} demonstrate some empirical advantages of Sinkhorn DRO.
\QEG

\end{remark}

\begin{remark}[Connection with KL-divergence DRO]\label{Remark:KL:DRO}
Using Jensen's inequality, we can see that the dual objective function of the Sinkhorn DRO model is upper-bounded as 
\[
\lambda\orho+\lambda\Reg\log\left(
\bE_{x\sim\hP}
\bE_{z\sim \bQ_{x,\Reg}}\left[ 
e^{f(z)/(\lambda\Reg)}
\right]\right),
\]
which corresponds to the dual objective for the following KL-divergence DRO~\citep{Ben13} 
\[
\sup_{\bP}~\left\{
\mathbb{E}_{z\sim\bP}[f(z)]:~
D_{\text{KL}}(\bP\|\bP^0)\le \orho/\Reg
\right\}.
\]
Here $\bP^0$ satisfies $\diff\bP^0(z) = \bE_{x\sim \hP}\left[\diff\bQ_{x,\Reg}(z)\right]$, which can be viewed as a non-parametric kernel density estimation constructed from $\hP$.
Particularly, when $\hP=\frac{1}{n}\sum_{i=1}^n \delta_{\hx_i}, \cZ=\mathbb{R}^d$ and $c(x,y)=\|x-y\|_2^2$, $\bP^0$ is kernel density estimator with Gaussian kernel and bandwidth $\epsilon$:
\[
\frac{\diff \bP^0(z)}{\diff z}
 = \frac{1}{n}\sum_{i=1}^nK_{\Reg}\left( 
z-x_i
\right),\quad z\in\mathbb{R}^d,
\]
where $K_\Reg(x)\propto \exp(-\|x\|_2^2/\Reg)$ represents the Gaussian kernel.
By Lemma~\ref{Lemma:refor:Inf-W-E} and divergence inequality~\citep[Theorem~2.6.3]{Cover06}, 
we can see the Sinkhorn DRO with $\orho=0$ is reduced to the following SAA model based on the distribution $\bP^0$:
\begin{equation}
V_{\Primal} = \mathbb{E}_{z\sim\bP^0}[f(z)] = \mathbb{E}_{x\sim\hP}\bE_{z\sim\bQ_{x,\Reg}}[f(z)].
\label{Eq:KDE:SAA}
\end{equation}
In non-parameteric statistics, the optimal bandwidth to minimize the mean-squared-error between the estimated distribution $\bP_0$ and the underlying true one is at rate $\Reg=O(n^{-1/(d+4)})$~\citep[Theorem~4.2.1]{hardle1990applied}.
However, such an optimal choice for the kernel density estimator may not be the optimal choice for optimizing the out-of-sample performance of the Sinkhorn DRO.
In our numerical experiments in Section \ref{Sec:numerical}, we select $\Reg$ based on cross-validation unless otherwise stated.

\Jie{
We also note that data-driven KL-divergence DRO is typically more efficient to solve than Sinkhorn DRO. This is because KL-divergence DRO yields a worst-case distribution supported on the same set as the empirical observations, leading to a finite-dimensional convex reformulation that can be efficiently solved using off-the-shelf solvers (see the numerical study of CPU times in Appendix~\ref{Sec:baseline:run}).
}
\QEG

\end{remark}

\begin{remark}[Connection with Bayesian DRO]
Bayesian DRO~\citep{shapiro2021bayesian}  proposed to solve
\[
\mathbb{E}_{x\sim \hP}~\left[
\sup_{\bP}~\Big\{ 
\mathbb{E}_{z\sim \bP}[f(z)]:~\bP\in \mathcal{P}_x
\Big\}
\right],
\]
where $\hP$ is a special posterior distribution constructed from collected observations, and the ambiguity set $\mathcal{P}_x$ is typically constructed as a KL-divergence ball, i.e., $\mathcal{P}_x:=\{\bP:~D_{\text{KL}}(\bP\|\bQ_x)\le \eta\}$, with $\bQ_x$ being the parametric distribution conditioned on $x$.
According to \citep[Section~2.1.3]{shapiro2021bayesian}, a relaxation of the Bayesian DRO dual formulation is given by
\[
\inf_{\lambda\ge0}~
\left\{ 
\lambda\eta + \lambda\bE_{x\sim \hP}\Big[ 
\log
\bE_{z\sim \bQ_x}\left[ 
e^{f(z) / \lambda}
\right]
\Big]
\right\}.
\]
When specifying the parametric distribution $\bQ_x$ as the kernel probability distribution in \eqref{Eq:bQ:x:Gibbs} and applying the change-of-variable technique such that $\lambda$ is replaced with $\lambda\Reg$, this relaxed formulation becomes
\[
\inf_{\lambda\ge0}~
\left\{ 
\lambda(\eta\Reg) + \lambda\Reg
\bE_{x\sim \hP}\Big[ 
\log
\bE_{z\sim\bQ_{x}}\left[ 
e^{f(z) / (\lambda\Reg)}
\right]
\Big]
\right\}.
\]
In comparison with \eqref{Eq:dual:SDRO:general:unified}, we find the Sinkhorn DRO model can be viewed as a special relaxation formulation of the Bayesian DRO model.
\QEG
\end{remark}

\subsection{Examples}
\label{Sec:example:reformulation}
In the following, we provide several cases in which our strong dual reformulation~\eqref{Eq:dual:SDRO:general:unified} can be simplified into more tractable formulations.

\begin{example}[Linear loss]\label{Example:DRO}
Suppose that the loss function $f(z)=a\trans z$, support $\cZ=\mathbb{R}^d$, $\nu$ is the corresponding Lebesgue measure, and the transport cost is the Mahalanobis distance, i.e., $c(x,y)=\frac{1}{2}(x-y)\trans\Omega(x-y)$, where $\Omega$ is a positive definite matrix.
In this case, the kernel probability distribution
$
\bQ_{x,\Reg}=\mathcal{N}(x, \Reg\Omega^{-1}).
$, 
and the dual problem can be written as 
\[
V_{\Dual} = \inf_{\lambda>0}~\left\{\lambda\orho+\lambda\Reg\bE_{x\sim\hP}\Big[\Lambda_x(\lambda)\Big]\right\},
\]
where 
\[
\Lambda_x(\lambda) = \log
\mathbb{E}_{z\sim \mathcal{N}(x,\Reg\Omega^{-1})}\left[ 
e^{a^\top z/(\lambda\Reg)}
\right]
= \frac{a\trans x}{\lambda\Reg} + \frac{a\trans\Omega^{-1}a}{2\lambda^2\Reg}.
\]
Therefore
\[
V_{\Dual}=a\trans\mathbb{E}_{x\sim\hP}[x] + \sqrt{2\orho}\sqrt{a\trans\Omega^{-1}a} := \mathbb{E}_{x\sim\hP}[a\trans x] + \sqrt{2\orho}\cdot \|a\|_{\Omega^{-1}}.
\]
This indicates that the Sinkhorn DRO is equivalent to an empirical risk minimization with norm regularization, and can be solved efficiently using algorithms for the second-order cone program.
\QEG
\end{example}

\begin{example}[Quadratic loss]\label{Example:DRO:LR}
Consider the example of linear regression with quadratic loss $f_{\theta}(z)=(a\trans\theta - b)^2$, where $z:=(a,b)$ denotes the predictor-response pair, $\theta\in\bR^d$ denotes the fixed parameter choice, and $\cZ=\mathbb{R}^{d+1}$.
Taking $\nu$ as the Lebesgue measure and the transport cost as $c((a,b), (a', b')) = \frac{1}{2}\|a-a'\|_2^2 + \infty|b-b'|$.
In this case, 
the dual problem becomes
\[
V_{\Dual} = \bE_{z\sim\hP}[(a\trans\theta-b)^2] + \inf_{
\lambda>2\|\theta\|_2^2
}~\left\{\lambda\orho+ \frac{\bE_{z\sim\hP}[(a\trans\theta-b)^2]}{\frac{1}{2}\lambda\|\theta\|_2^{-2}-1} - \frac{\lambda\Reg}{2}\log\det\left( 
I - \frac{\theta\theta\trans}{\frac{1}{2}\lambda}
\right)\right\}.
\]
In comparison with the corresponding Wasserstein DRO formulation with radius $\rho$~(see, e.g., \citep[Example~4]{blanchet2021statistical})
\[
V_{\Dual}^{\text{WDRO}} = \bE_{z\sim\hP}[(a\trans\theta-b)^2] + \inf_{
\lambda>2\|\theta\|_2^2
}~\left\{\lambda\rho+ \frac{\bE_{z\sim\hP}[(a\trans\theta-b)^2]}{\frac{1}{2}\lambda\|\theta\|_2^{-2}-1}\right\},
\]
one can check in this case the Sinkhorn DRO formulation is equivalent to the Wasserstein DRO with log-determinant regularization.
\QEG
\end{example}

When the support $\cZ$ is finite, the following result presents a conic programming reformulation.
\begin{corollary}[Conic Reformulation for Finite Support]\label{Corollary:conic:reformulation}
Suppose that the support contains $L_{\max}$ elements, i.e., $\cZ = \{z_{\ell}\}_{\ell=1}^{L_{\max}}$, and the nominal distribution $\hP = \frac{1}{n}\sum_{i=1}^n\delta_{\hx_i}$.
If Condition~\ref{Assumption:light:tailed:f} holds and $\orho\ge0$, the dual problem~\eqref{Eq:dual:SDRO:general:unified} can be formulated as the following conic optimization:
\begin{equation}\label{Eq:Sinkhorn:Conic:reformulation}
\begin{aligned}
V_{\Dual}=\min_{\substack{\lambda\ge0, s\in\mathbb{R}^n, \\ a\in\mathbb{R}^{n\times L}}}
&\quad 
\lambda\orho + \frac{1}{n}\sum_{i=1}^ns_i\\
\mbox{s.t.}&\quad \lambda\Reg\ge \sum_{\ell=1}^{L_{\max}}q_{i,\ell}a_{i,\ell}, i\in[n],\\
&\quad (\lambda\Reg, a_{i,\ell}, f(z_{\ell})-s_i)\in \mathcal{K}_{\exp}, i\in[n],\ell\in[L].
\end{aligned}
\end{equation}
where $q_{i,\ell}:=\Pr_{z\sim\bQ_{\hx_i,\Reg}}\{z = z_{\ell}\}$, with the distribution $\bQ_{\hx_i,\Reg}$ defined in \eqref{Eq:bQ:x:Gibbs}, and $\mathcal{K}_{\exp}$ denotes the exponential cone
$
\mathcal{K}_{\exp}:= 
\left\{
(\nu,\lambda,\delta)\in\mathbb{R}_+\times\mathbb{R}_+\times\mathbb{R}:~
\exp(\delta/\nu)\le \lambda/\nu
\right\}$.
\end{corollary}
Problem~\eqref{Eq:Sinkhorn:Conic:reformulation} is a convex program that minimizes a linear function with respect to linear and conic constraints, which can be solved using interior point algorithms~\citep{nesterov1994interior, vandenberghe1995semidefinite}.
We will develop an efficient first-order optimization algorithm in Section~\ref{Sec:first:order} that is able to solve a more general problem (without a finite support).

\subsection{Proof of Theorem~\ref{Theorem:strong:duality}}
\label{Proof:Theorem:strong:duality:main}
In this subsection, we present a sketch of the proof for Theorem~\ref{Theorem:strong:duality}.
We begin with the weak duality result in Lemma~\ref{Lemma:DRO:weak:duality}, which can be shown by  applying the Lagrangian weak duality.
\begin{lemma}[Weak Duality]
\label{Lemma:DRO:weak:duality}
Under Assumption~\ref{Assumption:distance:measure:function}, it holds that
\[
V_{\Primal}\le \inf_{\lambda\ge0}\;\bigg\{\lambda\rho + \lambda\Reg\; \mathbb{E}_{x\sim\hP}\Big[ \log\mathbb{E}_{z\sim\nu}\big[
   e^{
  (f(z) - \lambda c(x,z))/(\lambda\Reg)}
  \big] \Big]\bigg\}=V_{\Dual}.
\]
\end{lemma}

\proof{Proof of Lemma~\ref{Lemma:DRO:weak:duality}.}
Based on Definition~\ref{EntW:Wasserstein} of Sinkhorn distance, we reformulate $V_{\Primal}$ as 
\begin{align*}
V_{\Primal}&= 
\sup_{\gamma\in\cP(\cZ\times\cZ): \proj_{1\#}\gamma=\hP} 
\left\{ 
\mathbb{E}_{z\sim\bP}[f(z)]: \mathbb{E}_{(x,z)\sim\gamma}\left[ 
c(x,z) + \Reg\log\left(
\frac{\diff\gamma(x,z)}{\diff\hP(x)\diff\nu(z)}
\right)
\right]\le \rho
\right\}.
\end{align*}
By Assumption~\ref{Assumption:distance:measure:function},  the constraint is equivalent to
\[
\bE_{x\sim\hP}
\mathbb{E}_{z\sim \gamma_x}\left[ 
c(x, z) + \Reg\log\left(
\frac{\diff\gamma_x(z)}{\diff\nu(z)}
\right)
\right]
\le \rho,%
\]
and the primal problem is equivalent to
\begin{equation}
V_{\Primal}=\sup_{\{\gamma_x\}_{x\in\mathrm{supp}\,\hP}\subset \cP(\cZ)} \left\{
\bE_{x\sim\hP}\mathbb{E}_{z\sim\gamma_x}[f(z)]:\
\bE_{x\sim\hP}
\mathbb{E}_{z\sim \gamma_x}\left[ 
c(x, z) + \Reg\log\left(
\frac{\diff\gamma_x(z)}{\diff\nu(z)}
\right)
\right]
\le \rho
\right\}.
\label{Eq:V:primal:revision}
\end{equation}
Introducing the Lagrange multiplier $\lambda$ associated to the constraint, we reformulate $V_{\Primal}$ as 
\begin{align*}
V_{\Primal}=\sup_{\{\gamma_x\}_{x\in\mathrm{supp}\,\hP}\subset \cP(\cZ)}~\left\{\inf_{\lambda\ge0}~
\left\{\lambda\rho + 
\bE_{x\sim\hP}
\mathbb{E}_{z\sim\gamma_x}
\left[ 
f(z) - \lambda c(x,z) - \lambda\Reg\log\left(\frac{\diff \gamma_x(z)}{\diff\nu(z)}\right)
\right]\right\}
\right\}.
\end{align*}
Interchanging the order of the supremum and infimum operators, we have that
\begin{equation}\label{Eq:v:primal:upper}
V_{\Primal}\le \inf_{\lambda\ge0}~\left\{
\lambda\rho + \sup_{\{\gamma_x\}_{x\in\mathrm{supp}\,\hP}\subset \cP(\cZ)}
\left\{
\bE_{x\sim\hP}
\mathbb{E}_{z\sim\gamma_x}
\left[ 
f(z) - \lambda c(x,z) - \lambda\Reg\log\left(\frac{\diff \gamma_x(z)}{\diff\nu(z)}\right)
\right]
\right\}
\right\}.
\end{equation}
\Jie{
For $x\in\mathrm{supp}\,\hP$ and $\lambda\ge0$, define
\begin{equation}
v_x(\lambda) := \sup_{\gamma_x\in\cP(\cZ)}~
\left\{\mathbb{E}_{z\sim\gamma_x}\left[ 
f(z) - \lambda c(x,z) - \lambda\Reg\log\left(\frac{\diff \gamma_x(z)}{\diff\nu(z)}\right)
\right]\right\}.\label{Eq:v:x:lambda}
\end{equation}
Note that this function is measurable for any choice of $\lambda$ (we omit its proof in Lemma~\ref{Lemma:measure:vx}).
One can swap the supremum and the expectation operator on the right-hand-side of \eqref{Eq:v:primal:upper} to further upper bound it as
\[
V_{\Primal}\le \inf_{\lambda\ge0}~\left\{\lambda\rho + \bE_{x\sim\hP}\Big[v_{x}(\lambda)\Big]\right\}.
\]
By Lemma~\ref{Lemma:expression:vxlambda}, when there exists $\lambda>0$ such that Condition~\ref{Assumption:light:tailed:f} is satisfied,
it holds that
\[
v_x(\lambda) = 
\lambda\Reg\log\mathbb{E}_{z\sim\nu}\big[
   e^{
  (f(z) - \lambda c(x,z))/(\lambda\Reg)}
  \big]<\infty,
\]
and the desired result holds.
Otherwise, for any $\lambda>0$,
\[
\hP\left\{
x:~\mathbb{E}_{z\sim\bQ_{x,\Reg}}\left[ 
e^{f(z)/(\lambda\Reg)}
\right]=\infty
\right\}
=
\hP\left\{
x:~
\mathbb{E}_{z\sim\nu}\big[
   e^{
  (f(z) - \lambda c(x,z))/(\lambda\Reg)}
  \big]=\infty
\right\}
>0,
\]
then intermediately we obtain 
\[
V_{\Primal}\le \inf_{\lambda\ge0}\;\bigg\{\lambda\rho + \lambda\Reg\; \mathbb{E}_{x\sim\hP}\Big[ \log\mathbb{E}_{z\sim\nu}\big[
   e^{
  (f(z) - \lambda c(x,z))/(\lambda\Reg)}
\big] \Big]\bigg\}=\infty,
\]
and the weak duality still holds.}
\QED\endproof

Next, we show the feasibility result in Theorem~\ref{Theorem:strong:duality}\ref{Theorem:strong:duality:I}.
The key observation is that the primal problem \eqref{Inf-W-E} can be reformulated as a generalized KL-divergence DRO problem.
Consequently, Theorem~\ref{Theorem:strong:duality}\ref{Theorem:strong:duality:I} holds because of the non-negativity of KL-divergence.
\begin{lemma}[Reformulation of \eqref{Inf-W-E}]
\label{Lemma:refor:Inf-W-E}
Under Assumption~\ref{Assumption:distance:measure:function}, it holds that
\[
V_{\Primal}=\sup_{\{\gamma_x\}_{x\in\mathrm{supp}\,\hP}\subset \cP(\cZ)} \left\{\mathbb{E}_{x\sim\hP}\bE_{z\sim\gamma_x}[f(z)]:\
 \Reg
 \mathbb{E}_{x\sim\hP}\bE_{z\sim\gamma_x}\left[ 
\log\left(
\frac{\diff\gamma_x(z)}{\diff\bQ_{x,\Reg}(z)}
\right)
\right]
\le \orho
\right\}.
\]
\end{lemma}
\Jie{
\begin{remark}[Comparison with Infinite-dimensional Convex Analysis]
If assuming that $\hP=\frac{1}{n}\sum_{i=1}^n\delta_{x_i}$ is finitely supported, by Lemma~\ref{Lemma:refor:Inf-W-E}, \eqref{Inf-W-E} can be reformulated as a conic linear program
\[
V = \sup_{\{\gamma_i\}_{i\in[n]}\subset \cP(\cZ)} \left\{\frac{1}{n}\sum_{i\in[n]}\bE_{z\sim\gamma_i}[f(z)]:\
\frac{1}{n}\sum_{i\in[n]}D_{\text{KL}}(\gamma_i\|\bQ_{x_i,\Reg})
\le \frac{\orho}{\Reg}
\right\}.
\]
Thus, strong duality from infinite-dimensional convex analysis~(e.g., \citep{shapiro2001duality}) can be applied to show Theorem~\ref{Theorem:strong:duality}.
However, our strong duality proof, as described below, differs from this one in several key aspects.
First, our approach imposes less restrictive assumptions, holding for any measurable sample space $\cZ$, measurable loss function $f$, and nominal distribution $\hP$.
The strong duality result in \citep{shapiro2001duality} requires $\cZ$ to be convex, $f$ to be upper semicontinuous, and $\hP$ to be finitely supported.
Second, our approach is constructive: we explicitly characterize the worst-case distribution for Sinkhorn DRO, whereas a nonconstructive method was employed in \citep{shapiro2001duality}.
Third, our approach provides a byproduct — an explicit necessary and sufficient condition for when the Sinkhorn ambiguity constraint is binding~(Theorem~\ref{Theorem:strong:duality}\ref{Theorem:strong:duality:IV}). This insight offers practical guidance on choosing the ambiguity set size to avoid over-conservativeness.
\QEG
\end{remark}
}
Finally, we develop the strong duality.
\Jie{The general proof idea involves deriving the optimality condition of the dual minimizer, which then guide the construction of the worst-case distribution of \eqref{Inf-W-E}.
}
In the following, we provide the proof of the first part of Theorem~\ref{Theorem:strong:duality}\ref{Theorem:strong:duality:III} for the most representative case where $\orho>0$, the dual minimizer $\lambda^*$ exists with $\lambda^*>0$, and Condition~\ref{Assumption:light:tailed:f} holds. 
Proofs of other cases are moved in Appendix~\ref{Sec:appendix:proof}.
We first develop the optimality condition when the dual minimizer $\lambda^*>0$, by setting the derivative of the dual objective function to zero.

\begin{lemma}[First-order Optimality Condition when $\lambda^*>0$]\label{Lemma:lambda:positive}
Suppose $\orho>0$ and Condition~\ref{Assumption:light:tailed:f} is satisfied, and assume further that \Jie{there exists a dual minimizer $\lambda^*>0$, then the dual minimizer is unique and }$\lambda^*$  satisfies
\begin{equation}\label{Eq:optimality:lambda:positive:lemma:4}
\begin{aligned}
\Jie{%
\frac{1}{\lambda^*}\bE_{x\sim\hP}\left[
\frac{\mathbb{E}_{z\sim\nu}\left[ 
e^{(f(z)-\lambda^*c(x,z))/(\lambda^*\Reg)}f(z)
\right]}{\mathbb{E}_{z\sim\nu}\left[ 
e^{(f(z)-\lambda^*c(x,z))/(\lambda^*\Reg)}
\right]}
\right]
-\Reg\bE_{x\sim\hP}\left[ 
\log\bE_{z\sim\nu}\left[e^{(f(z) - \lambda^{\ast}c(x,z))/(\lambda^{\ast}\Reg)}\right]
\right]=\rho.
}
\end{aligned}
\end{equation}
\end{lemma}

\proof{Proof of Theorem~\ref{Theorem:strong:duality}\ref{Theorem:strong:duality:III} 
for the case where Condition~\ref{Assumption:light:tailed:f} holds and $\orho>0, \lambda^*>0$.}
We take the transport mapping $\gamma_*$ such that
\[
\frac{\diff\gamma_*(x,z)}{\diff\hP(x)\diff\nu(z)} = \alpha_x\cdot
  \exp\Big(\big( 
  f(z) - \lambda^\ast c(x,z) \big)/(\lambda^\ast\Reg)\Big),
\]
and $\alpha_x:=\left(
\bE_{z\sim\nu}\left[e^{(f(z) - \lambda^{\ast}c(x,z))/(\lambda^{\ast}\Reg)}\right]
  \right)^{-1}$ is a normalizing constant such that $\proj_{1\#}\gamma_* = \hP$.
Also define the primal (approximate) optimal distribution $\bP_\ast :=\proj_{2\#}\gamma_*.$
Recall the expression of the Sinkhorn distance in Definition~\ref{EntW:Wasserstein}, one can verify that
\[
\begin{aligned}
&\hphantom{=}\mathcal{W}_{\Reg}(\hP,\bP_\ast)
=\inf_{\gamma\in\Gamma(\hP,\bP_\ast)}~\left\{ 
\mathbb{E}_{(x,z)\sim\gamma}\left[ 
c(x,z) + \Reg\log\left( 
\frac{\diff\gamma(x,z)}{\diff\hP(x)\diff\nu(z)}
\right)
\right]
\right\}
\\
&=\inf_{\gamma\in\Gamma(\hP,\bP_\ast)}~\left\{ 
\mathbb{E}_{(x,z)\sim \gamma}\left[ 
\Reg\log\left(
\frac{e^{c(x,z)/\Reg}\diff\gamma(x,z)}{\diff\hP(x)\diff\nu(z)}
\right)
\right]
\right\}\\
&\le \mathbb{E}_{(x,z)\sim\gamma_*}\left[ 
\Reg\log\left(
\frac{e^{c(x,z)/\Reg}\diff\gamma_*(x,z)}{\diff\hP(x)\diff\nu(z)}
\right)
\right]
=\mathbb{E}_{(x,z)\sim\gamma_*}\left\{
\frac{1}{\lambda^*}f(z) + \Reg\log(\alpha_x)
\right\}\\
&=
\frac{1}{\lambda^*}\bE_{x\sim\hP}\left[
\frac{\mathbb{E}_{z\sim\nu}\left[ 
e^{(f(z)-\lambda^*c(x,z))/(\lambda^*\Reg)}f(z)
\right]}{\mathbb{E}_{z\sim\nu}\left[ 
e^{(f(z)-\lambda^*c(x,z))/(\lambda^*\Reg)}
\right]}
\right]
-\Reg\bE_{x\sim\hP}\left[ 
\log\bE_{z\sim\nu}\left[e^{(f(z) - \lambda^{\ast}c(x,z))/(\lambda^{\ast}\Reg)}\right]
\right]
\end{aligned}
\]
where the inequality relation is because $\gamma_*$ is a feasible solution in $\Gamma(\hP,\bP_\ast)$, and the last two relations are by substituting the expression of $\gamma_*$.
Since $\orho>0$ and the dual minimizer $\lambda^*>0$, the optimality condition in \eqref{Eq:optimality:lambda:positive:lemma:4} holds, which implies that $\mathcal{W}_{\Reg}(\hP,\bP_\ast)\le\rho$, i.e., the distribution $\bP_\ast$ is primal feasible for the problem~\eqref{Inf-W-E}.
Moreover, we can see that the primal optimal value is lower bounded by the dual optimal value:
\begin{align*}
V_{\Primal}~\ge&~\mathbb{E}_{\bP_\ast}[f(z)]=
\bE_{(x,z)\sim\gamma_{\ast}}[f(z)]\\
=&
\bE_{x\sim\hP}\bE_{z\sim \nu}\left[ 
f(z)\left(\frac{\diff\gamma_*(x,z)}{\diff\hP(x)\diff\nu(z)}\right)
\right]
=\bE_{x\sim\hP}\left[
\frac{\mathbb{E}_{z\sim\nu}\left[ 
e^{(f(z)-\lambda^*c(x,z))/(\lambda^*\Reg)}f(z)
\right]}{\mathbb{E}_{z\sim\nu}\left[ 
e^{(f(z)-\lambda^*c(x,z))/(\lambda^*\Reg)}
\right]}
\right]\\
=&\lambda^*\left(
\orho + \Reg\bE_{x\sim\hP}\left[\log
\mathbb{E}_{z\sim\bQ_{x,\Reg}}\left[ 
e^{f(z)/(\lambda^*\Reg)}
\right]\right]
\right)
=V_{\Dual},
\end{align*}
where the third equality is by substituting the expression of $\gamma_{\ast}$, and the last equality is based on the optimality condition in \eqref{Eq:optimality:lambda:positive:lemma:4}.
This, together with the weak duality, completes the proof.
\QED\endproof

\section{Efficient First-order Algorithm for Sinkhorn Robust Optimization}
\label{Sec:first:order}
Consider the Sinkhorn robust optimization problem
\begin{equation}\label{Eq:statistical:learning}
\begin{aligned}
\inf_{\theta\in\Theta}\sup_{\mathbb{P}\in\mathbb{B}_{\rho,\Reg}(\hP)}~
\mathbb{E}_{z\sim \mathbb{P}}[f_{\theta}(z)].
\end{aligned}
\end{equation}
Here the feasible set $\Theta\subseteq\mathbb{R}^{d_{\theta}}$ is \emph{closed and convex} containing all possible candidates of decision vector $\theta$, and the Sinkhorn uncertainty set is centered around a given nominal distribution $\hP$.
Based on our strong dual expression \eqref{Eq:dual:SDRO:general:unified}, we reformulate \eqref{Eq:statistical:learning} as
\begin{equation}\label{Eq:reformulate:statistical:learning}
\inf_{\lambda\ge0}~\left\{\lambda\orho+\inf_{\theta\in\Theta}~
\mathbb{E}_{x\sim\hP}\bigg[ 
\lambda\Reg\log
\mathbb{E}_{z\sim\bQ_{x,\Reg}}\left[ 
e^{f_{\theta}(z)/(\lambda\Reg)}
\right]
\bigg]
\right\},
\tag{\textsf{D}}
\end{equation}
where the constant $\orho$ and the distribution $\bQ_{x,\Reg}$ are defined in \eqref{Eq:def:orho:general:unified} and \eqref{Eq:bQ:x:Gibbs}, respectively.
In Examples~\ref{Example:DRO} and \ref{Example:DRO:LR}, we have seen special instances of \eqref{Eq:reformulate:statistical:learning} where we can get closed-form expressions for the above integration.
In this section, we develop an efficient algorithm for solving \eqref{Eq:reformulate:statistical:learning} for general loss functions where a closed-form expression is not available.

A typical approach for solving a stochastic optimization is the stochastic (sub)gradient method such as stochastic mirror descent~(SMD)~\cite{Blair85}.
Unlike many other stochastic optimization problems, one salient feature of \eqref{Eq:reformulate:statistical:learning} is that its inner objective involves a nonlinear transformation of the expectation.
Consequently, based on a batch of simulated samples from $\bQ_{x,\Reg}$, an unbiased subgradient estimate could be challenging to obtain.
In Section \ref{sec:alg}, we will combine SMD with biased subgradient estimators and bisection search to solve \eqref{Eq:reformulate:statistical:learning}.
We will analyze its computational complexity in Section~\ref{sec:convergence}.

\subsection{\Jie{Algorithm Framework}}\label{sec:alg}
\Jie{
We define 
\begin{equation}
  F(\theta; \lambda) := \mathbb{E}_{x\sim\hP}\bigg[ 
  \lambda\Reg\log
  \mathbb{E}_{z\sim\bQ_{x,\Reg}}\left[ 
  e^{f_{\theta}(z)/(\lambda\Reg)}
  \right]
  \bigg],
  \label{Expression:F:theta}
\end{equation}
and define the objective value of the outer minimization in \eqref{Eq:reformulate:statistical:learning} as
\begin{equation}
\Jie{\Psi(\lambda)}
:=
\lambda\orho+\inf_{\theta\in\Theta}%
~F(\theta;\lambda).
\label{Eq:F:lambda}
\end{equation}
}
\Jie{
Solving \eqref{Eq:reformulate:statistical:learning} involves determining the optimal Lagrange multiplier $\lambda$ of the function $\Psi$, where evaluating $\Psi$ requires solving the minimization problem in \eqref{Eq:F:lambda} to obtain the optimal decision $\theta$. We first introduce a biased SMD~(BSMD) algorithm for finding the optimal decision $\theta$ given a fixed Lagrange multiplier $\lambda$ in Section~\ref{Sec:configure:BSMD}. Then, in Section~\ref{sec:alg:overall}, we present a bisection search algorithm to find the optimal Lagrange multiplier.
\Jie{
Throughout this section, we assume the loss function $f_{\theta}(z)$ is convex in $\theta$ but it can be a potentially nonsmooth function.
For any function $r(\theta)$ that is subdifferentiable in $\theta$, we use the notation $\nabla_{\theta}r(\theta)$ to denote an arbitrary subgradient from its subdifferential, unless otherwise specified.
}
}

\subsubsection{BSMD.} \label{Sec:configure:BSMD}

\Jie{
In this part, we omit the dependence of $\lambda$ when defining objective or subgradient terms, e.g., we write $F(\theta)$ for $F(\theta;\lambda)$.}
We first introduce several notations that are standard in the mirror descent algorithm.
Let $\omega:~\Theta\to\mathbb{R}$ be a distance generating function that is continuously differentiable and $\kappa$-strongly convex on $\Theta$ with respect to norm $\|\cdot\|$, which induces the Bregman divergence $D_{\omega}(\theta, \theta'):~\Theta\times\Theta\to\mathbb{R}_+$:
$D_{\omega}(\theta, \theta') = \omega(\theta') - 
\omega(\theta) -\inp{\nabla\omega(\theta)}{\theta' - \theta}.$
Define the \emph{prox-mapping} $\prox:~\mathbb{R}^{d_{\theta}}\to\Theta$ as
\[
\prox_{\theta}(y) = \argmin_{\theta'\in\Theta}\big\{ 
\inp{y}{\theta' - \theta} + D_{\omega}(\theta, \theta')
\big\}.
\]
With these notations in hand, we present our algorithm in Algorithm \ref{alg:BSMD:sampling}\Jie{, which} iteratively obtains a biased stochastic (sub)gradient estimator and performs a proximal update.
\begin{algorithm}[!ht]
\caption{
BSMD for finding the optimal solution of \eqref{Eq:F:lambda} while fixing $\lambda$}
\label{alg:BSMD:sampling} 
\begin{algorithmic}[1] %
\REQUIRE
{
Maximal iteration $\Jie{T}$, constant step size $\Jie{h}$, initial guess $\theta_0$, fixed multiplier $\lambda$.
}
\FOR{$t=0,1,\ldots,\Jie{T}-1$}
\STATE{Construct a (biased) subgradient estimator $v(\theta_t)$ of $F(\theta_t)$ using \eqref{Eq:SGD:scheme} or \eqref{Eq:RTMLMC:scheme}.
}
\STATE{
Update $\theta_{t+1} = \prox_{\theta_t}\big(\Jie{h} v(\theta_t)\big)
$.
}
\ENDFOR\\
\textbf{Output} the estimate of optimal solution $\Jie{\widehat{\theta}} = \frac{1}{\Jie{T}}\sum_{t=1}^{\Jie{T}}\theta_t$.
\end{algorithmic}
\end{algorithm}

\Jie{
At the core of Algorithm~\ref{alg:BSMD:sampling} lies the crucial task of efficiently simulating the subgradient estimator in Step~2.
It is noteworthy that the minimization in \eqref{Eq:F:lambda} is a special conditional stochastic optimization~(CSO), as studied in \citep{hu2021biasvar, Yifan20, hu2020sample}.
CSO typically has the formulation
\begin{equation}
\min_{\theta\in\Theta}~\bE_{x\sim\hP}[H^1(\bE_{z\sim\bQ_{x,\Reg}}[H^2(\theta; z)])],\label{Eq:CSO:speific}
\end{equation}
and we specify $H^1(\boldsymbol{\cdot})=\lambda\Reg\log(\boldsymbol{\cdot})$ and $H^2(\boldsymbol{\cdot};z)=\exp(f_{\boldsymbol{\cdot}}(z) / (\lambda\Reg))$ to convert \eqref{Eq:F:lambda} into \eqref{Eq:CSO:speific}.
This structure suggests that ideas from CSO-related literature, particularly multi-level Monte Carlo (MLMC) estimators, can be applied to generate biased subgradient estimators with controlled bias, enhancing computational efficiency in Step~2.
Our framework and analysis differs from the aforementioned references in several aspects, and see the discussion in Remark~\ref{Remark:compare:hu}.}

To generate the subgradient estimator, we first construct a function $F^{\ell}(\theta), \ell\in\mathbb{N}$ that approximate the original objective function $F(\theta)$ with $\cO(2^{-\ell})$-gap:
\begin{equation}\label{Eq:stat:robust:formula:approximation}
F^{\ell}(\theta) = \bE_{x^{\ell}\sim\hP}~\bE_{\{z^{\ell}_j\}_{j\in[2^{\ell}]}\sim \bQ_{x^{\ell}, \Reg}}\bigg[ \lambda\Reg\log\Big(
\frac{1}{2^{\ell}}\sum_{j\in[2^{\ell}]}e^{ 
f_{\theta}(z^{\ell}_j)/(\lambda\Reg)
}
\Big)
\bigg],
\end{equation}
where the random variable $x^{\ell}$ follows distribution $\hP$, and given a realization of $x^{\ell}$, $\{z^{\ell}_j\}_{j\in[2^{\ell}]}$ are independent
and identically distributed~(i.i.d.) samples from $\bQ_{x^{\ell}, \Reg}$.
Unlike the original objective $F(\theta)$, unbiased subgradient estimators of its approximation $F^{\ell}(\theta)$ can be easily obtained. 
Denote by $\zeta^{\ell}=(x^{\ell}, \{z^{\ell}_j\}_{j\in[2^{\ell}]})$ the collection of random sampling parameters, and 
\Jie{
\[
\begin{aligned}
U_{n_1:n_2}(\theta,\zeta^{\ell}) &:= \lambda\Reg\log\Big(
\frac{1}{n_2-n_1+1}\sum_{j\in[n_1:n_2]}e^{ 
f_{\theta}(z^{\ell}_j)/(\lambda\Reg)
}\Big),\\
A^{\ell}(\theta,\zeta^{\ell}) &:=U_{1:2^{\ell}}(\theta,\zeta^{\ell}) - \frac{1}{2}U_{1:2^{\ell-1}}(\theta,\zeta^{\ell})
- \frac{1}{2}U_{2^{\ell-1}+1:2^{\ell}}(\theta,\zeta^{\ell}).
\end{aligned}
\]}
\Jie{We take $U_{n_1:n_2}(\theta,\zeta^{\ell})=0$ if $[n_1:n_2]=\emptyset$.}
For fixed $\theta$ and $\ell\in\mathbb{N}$, we define 
\Jie{
\begin{align*}
   g^{\ell}(\theta,\zeta^{\ell})&:=\nabla_{\theta} U_{1:2^{\ell}}(\theta,\zeta^{\ell}), \qquad    G^{\ell}(\theta,\zeta^{\ell}):=
   \nabla_{\theta}A^{\ell}(\theta,\zeta^{\ell}).
\end{align*}}
The random vector $g^{\ell}(\theta,\zeta^{\ell})$ is an unbiased estimator of $\nabla_{\theta} F^{\ell}(\theta)$, while the random vector $G^{\ell}(\theta, \zeta^{\ell})$ is an unbiased estimator of $\nabla_{\theta} F^{\ell}(\theta) - \nabla_{\theta} F^{\ell-1}(\theta)$.
\Jie{
We note that computing $G^{\ell}(\theta,\zeta^{\ell})$ involves computing subgradient vectors $\nabla_{\theta}f_{\theta}(z_j^{\ell}), j\in[2^{\ell}]$, and we use the same subgradient computation across $U_{1:2^{\ell}}, U_{1:2^{\ell-1}}, U_{2^{\ell-1}+1:2^{\ell}}$ to facilitate the reduction of the second-order moment of $G^{\ell}(\theta,\zeta^{\ell})$.
Using these components, we define two types of subgradient estimators below.
}

\Jie{
\noindent{-- \textit{Stochastic (sub)Gradient~(SG) Estimator:}}
Fix the maximum level $L\in\mathbb{N}_+$.
We first generate the sample set $\zeta^L$ and next construct the SG estimator
\begin{equation}\label{Eq:SGD:scheme}
v^{\text{SG}}(\theta) =g^{L}(\theta,\zeta^{L}). %
\end{equation}
\noindent{-- \textit{Randomized Truncation MLMC~(RT-MLMC) Estimator~\citep{blanchet2015unbiased}: }}
Fix the maximum level $L\in\mathbb{N}_+$.
We first sample a random level $\widehat{\ell}$ following a truncated geometric distribution 
\begin{equation}
p_{\ell}:=\Pr(\widehat{\ell}=\ell)=\frac{2^{-\ell}}{2 - 2^{-L}}, \ell=0,1,\ldots,L.\label{Eq:p:ell:dist}
\end{equation}
Next, we construct the RT-MLMC estimator
\begin{equation}\label{Eq:RTMLMC:scheme}
v^{\text{RT-MLMC}}(\theta)=
p_{\widehat{\ell}}^{-1}
\cdot 
G^{\widehat{\ell}}(\theta,\zeta^{\widehat{\ell}}).
\end{equation}
}

\begin{remark}[Sampling from $\bQ_{x,\Reg}$]\label{Remark:sample:Q}
\Jie{
Sampling from $\bQ_{x,\Reg}$ is crucial for the construction of subgradient estimators.
In many cases, it is an easy task:
}
When the transport cost $c(\cdot,\cdot)=\frac{1}{2}\|\cdot-\cdot\|_2^2$ and $\cZ=\mathbb{R}^d$, the distribution $\bQ_{x,\Reg}$ becomes a Gaussian distribution $\mathcal{N}(x,\Reg I_d)$.
When the transport cost $c(\cdot,\cdot)$ is decomposable in each coordinate, we can apply the acceptance-rejection method~\citep{asmussen2007stochastic} to generate samples in each coordinate independently, the complexity of which only increases linearly in the data dimension.
\QEG
\end{remark}

\subsubsection{Bisection Search}
\label{sec:alg:overall}

\Jie{In this part, we introduce a bisection search algorithm to solve the one-dimensional convex minimization problem in \eqref{Eq:reformulate:statistical:learning}. 
The algorithm relies on an efficient oracle to estimate the objective value of \eqref{Eq:reformulate:statistical:learning}. We first define this oracle in Algorithm~\ref{Alg:inexact:obj}: Given a fixed multiplier $\lambda$, it solves problem~\eqref{Eq:F:lambda} using Algorithm~\ref{alg:BSMD:sampling} and then estimates the corresponding objective value. 
}
It has $m$ independent repetitions, whose value will be determined later in Section \ref{sec:convergence} to achieve the optimal complexity.

\begin{algorithm}[!ht]
\caption{
Evaluating the objective value of \eqref{Eq:reformulate:statistical:learning}
}
\label{Alg:inexact:obj}
\begin{algorithmic}[1] %
\REQUIRE{
Fixed multiplier $\lambda$, error tolerance $\delta$, batch size $m$.
}
\FOR{$j=1,2,\ldots,m$}
\STATE{
Obtain a $\delta$-optimal solution $\widehat{\theta}_j$ of problem~\eqref{Eq:F:lambda} using Algorithm \ref{alg:BSMD:sampling}.
}
\STATE{
Estimate the objective in \eqref{Expression:F:theta} with $\theta\equiv \widehat{\theta}_j$ using RT-MLMC estimator~\eqref{Eq:RTMLMC:est}, denoted as $\hF(\widehat{\theta}_j;\lambda)$
.
}
\ENDFOR\\
{\bf Output} 
\Jie{$\widehat{\Psi}(\lambda):=\lambda\orho + \underset{\theta\in\{
\widehat{\theta}_1,\ldots,\widehat{\theta}_m
\}}{\min}\hF(\theta;\lambda)$.}
\end{algorithmic}
\end{algorithm}

\Jie{
To implement Step~3 of Algorithm \ref{Alg:inexact:obj}, we again leverage RT-MLMC to efficiently estimate the objective in~\eqref{Expression:F:theta}.
For given $(\theta,\lambda)$, let $m'$ denote the mini-batch size.
For $i\in[m']$, we sample $\widehat{\ell}_i$ following the distribution defined in \eqref{Eq:p:ell:dist} and sample an i.i.d. copy of $\zeta^{\widehat{\ell}_i}$ that is denoted as $\zeta^{\widehat{\ell}_i}_i$.
Next, we construct the objective estimator at the $i$-th trial as
\[
\hF_i(\theta;\lambda)=
p_{\widehat{\ell}_i}^{-1}\cdot A^{\widehat{\ell}_i}(\theta,\zeta^{\widehat{\ell}_i}_i;\lambda).
\]
To reduce the variance of objective estimator, the final estimator of $F(\theta;\lambda)$ is constructed by averaging the outcomes over all trials, denoted as 
\begin{equation}\label{Eq:RTMLMC:est}
\hF(\theta;\lambda) = \frac{1}{m'}\sum_{i=1}^{m'}~\hF_i(\theta;\lambda).
\end{equation}
}
\Jie{
Given an inexact objective oracle of \eqref{Eq:reformulate:statistical:learning} (e.g., by querying Algorithm~\ref{Alg:inexact:obj}), we use bisection search to find a near-optimal multiplier in \eqref{Eq:reformulate:statistical:learning}; see Algorithm~\ref{Alg:outer:bisection} for details.
Unlike conventional bisection that relies on gradient information, this algorithm leverages an inexact objective oracle $\widehat{\Psi}$ to iteratively shrink the search interval.
It begins by dividing the interval into five evenly spaced points and selecting the minimum among the three central points.
In each iteration, it updates the left, middle, and right points based on the the current minimum (from among the three middle points) and its two nearest neighbors. Then, it adjusts the middle-left and middle-right points to maintain evenly spacing. 
The oracle is queried twice per iteration, reusing one value from the previous iteration to identify the new minimum efficiently.
This algorithm is adopted from~\citep[Algorithm~8]{cohen2016geometric}, but is more efficient, as the original algorithm only shrinks the interval by $1/3$ with each iteration, whereas one can improve it to factor $1/2$.
See its performance guarantee in Theorem~\ref{Proposition:complexity:bisection}.
}

\begin{algorithm}[!ht]
\caption{
Bisection search for finding the optimal multiplier of \eqref{Eq:reformulate:statistical:learning}}
\label{Alg:outer:bisection}
\Jie{
\begin{algorithmic}[1] %
\REQUIRE
{
Interval $[\lambda_{l}, \lambda_u]$ such that $\lambda_l<\lambda^*<\lambda_u$, maximum iterations $\Jie{T'}$
}
\REQUIRE
{
Inexact objective oracle $\widehat{\Psi}(\cdot):~\bR_+\to\bR$
}
\STATE{
Update 
$\beta_i^{(0)} = \frac{5-i}{4}\lambda_l + \frac{i-1}{4}\lambda_u$ for $i=1,\ldots,5$
}
\hfill{\COMMENT{\textit{Divide interval using 5 grid points}}}
\STATE{
Query oracle to obtain $\widehat{\Psi}(\beta_j^{(0)}), j=2,3,4$
}
\STATE{
Specify $i^{(1)}=\underset{j=2,3,4}{\argmin}~\widehat{\Psi}(\beta_j^{(0)})$
}
\STATE{\textbf{for $t=1,\ldots,\Jie{T'}$ do}}
\STATE{
~~~~~~Update $(\beta_1^{(t)}, \beta_3^{(t)}, \beta_5^{(t)})
=(\beta^{(t-1)}_{i^{(t)}-1}, \beta^{(t-1)}_{i^{(t)}}, \beta^{(t-1)}_{i^{(t)}+1})
$
}
\hfill{\COMMENT{\textit{Move left, middle, right points}}}
\STATE{
~~~~~~Update $(\beta_2^{(t)}, \beta_4^{(t)})
=\left(
\frac{\beta_1^{(t)} + \beta_3^{(t)}}{2}, 
\frac{\beta_3^{(t)} + \beta_5^{(t)}}{2}
\right)
$
}
\hfill{\COMMENT{\textit{Move middle-left, middle-right points}}}
\STATE{
~~~~~~Query oracle to obtain $\widehat{\Psi}(\beta_j^{(t)}), j=2,4$, and update $\widehat{\Psi}(\beta_3^{(t)})=\widehat{\Psi}(\beta^{(t-1)}_{i^{(t)}})$.
}
\STATE{
~~~~~~Specify $i^{(t+1)}=\underset{j=2,3,4}{\argmin}~\widehat{\Psi}(\beta_j^{(t)})$
}
\STATE{\textbf{end for}}\\
\textbf{Output} approximate optimal multiplier $\beta_{i^{(T'+1)}}^{(\Jie{T'})}$.
\end{algorithmic}
}
\end{algorithm}

\Jie{An alternative approach to solving \eqref{Eq:reformulate:statistical:learning} is to jointly optimize $(\lambda,\theta)$ using BSMD.
Although this approach is theoretically sound, updating $\lambda$ could lead to oscillations or divergence if the associated stepsize is not carefully tuned, due to the high variance of gradient estimators when $\lambda$ is small.
In our algorithm, we have developed a bisection method to update $\lambda$, which requires only specifying the maximum iterations, initial interval, and an inexact objective oracle but does not require tuning the stepsize for updating $\lambda$.
}
We also remark that, as a practical alternative, one can solve \eqref{Eq:F:lambda} using Algorithm \ref{alg:BSMD:sampling} alone and tune the hyperparameter $\lambda$, as tuning the radius $\orho$ is equivalent to tuning the Lagrangian multiplier $\lambda$ in \eqref{Eq:F:lambda}.
This corresponds to the Sinkhorn robust learning problem with a soft Sinkhorn constraint.

\subsection{Convergence Analysis}\label{sec:convergence}
In this subsection, we analyze the convergence properties of the proposed algorithms.
We begin with the following assumptions on the loss function $f_{\theta}$:
\begin{assumption}\label{Assumption:throughout:loss}
\begin{enumerate}
    \item\label{Assumption:throughout:loss:cvx}(Convexity): The loss function $f_{\theta}(z)$ is convex in $\theta$.
    \item\label{Assumption:throughout:loss:lip}(Lipschitz Continuity): For any fixed $z$ and $\theta_1,\theta_2$, it holds that $|f_{\theta_1}(z) - f_{\theta_2}(z)|\le L_f\|\theta_1 - \theta_2\|_2$.
    \item\label{Assumption:throughout:loss:bound}(Boundedness): The loss function $f_{\theta}(z)$ satisfies $0\le f_{\theta}(z)\le B$ for any $\theta\in\Theta$ and $z\in\cZ$.
\end{enumerate}
\end{assumption}

\Jie{
Assumption~\ref{Assumption:throughout:loss}\ref{Assumption:throughout:loss:cvx} ensures the convexity of the objective in Sinkhorn robust optimization, and enables us to develop globally convergent optimization algorithms.
Assumption~\ref{Assumption:throughout:loss}\ref{Assumption:throughout:loss:lip} is crucial for establishing the bounded subgradient norm condition required by the BSMD algorithm and deriving its global convergence rate.
Assumption~\ref{Assumption:throughout:loss}\ref{Assumption:throughout:loss:bound} guarantees the Lipschitz continuity of the nonlinear operator $H^1(\boldsymbol{\cdot})=\lambda\Reg\log(\boldsymbol{\cdot})$ in \eqref{Expression:F:theta}, ensuring that the objective in \eqref{Eq:stat:robust:formula:approximation} approximates \eqref{Expression:F:theta} with $\mathcal{O}(2^{-\ell})$ gap.
This assumption is restrictive in practice, whereas one can replace it by the following conditions (see the argument in \citep[footnote~2]{levy2020large}): (III-1): $\Theta$ is bounded; (III-2): Assumption~\ref{Assumption:throughout:loss}\ref{Assumption:throughout:loss:lip}; and (III-3): $\inf\limits_{\theta\in\Theta}~f_{\theta}(z_1) - \inf\limits_{\theta\in\Theta}~f_{\theta}(z_2)\le B_0, \forall z_1,z_2\in\cZ$.
}

\Jie
{
We note that our algorithm is assumed to have access to two sampling oracles:
(i) Oracle $\mathbf{O}(\hP)$ that generates a sample from $\hP$;
(ii) Oracle $\mathbf{O}(\bQ_{x,\Reg})$ that, based on the input $x\in\mathrm{supp}\, \hP$, generates a sample from $\bQ_{x,\Reg}$.
In practical implementations, the cost of generating samples from these two distributions can differ. In data-driven applications, sampling from $\hP$ often reduces to randomly selecting observed data points, whereas sampling from $\bQ_{x,\Reg}$ (described in Remark~\ref{Remark:sample:Q}) usually involves stochastic noises such as Gaussian.
In the subsequent analysis, we report the sample complexities from $\hP$ or $\bQ_{x,\Reg}$ individually.
Based on our algorithm design~(see, e.g., objective and subgradient estimators in \eqref{Eq:SGD:scheme}, \eqref{Eq:RTMLMC:scheme}, \eqref{Eq:RTMLMC:est}), it is also easy to check that the total computational time is roughly proportional to the sum of these two sample complexities.
}

\subsubsection{Complexity of BSMD.
}\label{Sec:inner:step:1}
In this part, we discuss the complexity of Algorithm~\ref{alg:BSMD:sampling}.
We say $\theta$ is a $\delta$-optimal solution if \Jie{$\bE[F(\theta;\lambda)] - F(\theta^*;\lambda)\le \delta$}, where $\theta^*$ is the optimal solution of \eqref{Eq:F:lambda}.
By properly tuning hyper-parameters to balance the trade-off between bias and second-order moment of the subgradient estimate, we establish its performance guarantees in Theorem~\ref{Theorem:complexity:BSMD}.
The explicit constants and proof can be found in Appendix~\ref{Appendix:convergence:analysis}.
\Jie{
Let us define the constant $K_{\lambda,\Reg,B}=\frac{B}{\lambda\Reg}$ that  depends on $\lambda,\Reg,B$.
}

\begin{theorem}
\label{Theorem:complexity:BSMD}
Under \Jie{Assumption~\ref{Assumption:throughout:loss},}
when using BSMD~(Algorithm~\ref{alg:BSMD:sampling}) to find a $\delta$-optimal solution of \eqref{Eq:F:lambda}, the following results hold:
\Jie{
\begin{enumerate}
\item\label{Theorem:complexity:BSMD:nonsmooth}
If using SG subgradient estimator,
the sample complexity from $\hP$ is $\cO(\delta^{-2})$, and 
that from $\bQ_{x,\Reg}$ is $\cO(\lambda\Reg\exp(2K_{\lambda,\Reg,B})\cdot \delta^{-3})$, with $\cO(\cdot)$ hiding constants depending on $L_f, \theta_0, \kappa$.
    \item \label{Theorem:complexity:BSMD:smooth}
If using RT-MLMC subgradient estimator,
the sample complexity from $\hP$ is $\tO(K_{\lambda,\Reg,B}\exp(4K_{\lambda,\Reg,B})\cdot\delta^{-2})$, and 
that from $\bQ_{x,\Reg}$ is $\tO(K_{\lambda,\Reg,B}^2\exp(4K_{\lambda,\Reg,B})\cdot\delta^{-2})$, with $\tO(\cdot)$ hiding constants depending on $L_f, \theta_0, \kappa$ and linearly depending on $(\log\frac{\lambda\Reg}{\delta})^2$.
\end{enumerate}   
}
\end{theorem}

\Jie{
Theorem~\ref{Theorem:complexity:BSMD} shows that the sample complexity from $\hP$ of BSMD, whether using the SG or RT-MLMC subgradient estimator, is of the same order with respect to the error tolerance $\delta$. 
This rate matches the known lower bound for general convex stochastic programming problems~\citep{Blair85}. However, this complexity associated with the RT-MLMC estimator has a worse constant dependence on the parameters $\lambda, \Reg, $ and $B$.
Despite this, the RT-MLMC estimator has a lower-order sample complexity from $\bQ_{x,\Reg}$ compared to the SG estimator. Our numerical experiments in Appendices~\ref{Appendix:compare:opt:alg} and~\ref{Appendix:compare:opt:alg:port} further demonstrate that the RT-MLMC estimator exhibits a significantly faster empirical convergence rate than the SG estimator.
}

\begin{remark}[Comparison with Biased Sample Average Approximation]\label{Remark:BSAA}
An alternative way to solve \eqref{Eq:F:lambda} is to approximate the objective using finite samples for both expectations. This leads to a biased sample estimate, called Biased Sample Average Approximation~(BSAA).
\Jie{
Under Assumption~\ref{Assumption:throughout:loss} and apply~\citep[Corollary~4.2]{hu2020sample}, it can be shown that for BSAA, the sample complexity from $\hP$ is $n_1=\tO\left(
d_{\theta}
B^2\exp(2K_{\lambda,\Reg,B})
\cdot \delta^{-2}
\right)$ and that from $\bQ_{x,\Reg}$ is $\cO(\lambda\Reg \exp(2K_{\lambda,\Reg,B})\cdot n_1\cdot \delta^{-1})$.
Our proposed BSMD with the RT-MLMC-based subgradient estimator has smaller order of the sample complexity from $\bQ_{x,\Reg}$ (in terms of error tolerance $\delta$).
}
Also, the BSAA method still requires computing the optimal solution of the approximated optimization problem as the output.
Hence, it typically takes considerably less time and memory to run the BSMD step rather than solving the BSAA formulation.
\QEG
\end{remark}

\subsubsection{
Complexity of Bisection Search
}
\label{Sec:inner:iteration:step:2}

We first provide the complexity analysis for Algorithm~\ref{Alg:inexact:obj}, which produces an estimator of the objective value of the outer minimization in \eqref{Eq:reformulate:statistical:learning}.
\Jie{Define the constant $H_{\lambda,\Reg,B} = \max(\exp(2K_{\lambda,\Reg,B}), \lambda^2\Reg^2)$.}

\begin{proposition}
\label{Proposition:com:estimate:optval}
Let $\eta\in(0,1)$ and set the batch size $m=\lceil \log_2\frac{2}{\eta}\rceil$.
Assume Assumption~\ref{Assumption:throughout:loss} holds, %
and we choose hyper-parameters in Step~3 of Algorithm~\ref{Alg:inexact:obj} as
\[
\Jie{
L=\left\lceil\log_2\frac{2\lambda\Reg\exp(2K_{\lambda,\Reg,B})}{\delta}\right\rceil,\quad 
m' = \cO(1)\frac{\lambda^2\Reg^2\exp(2K_{\lambda,\Reg,B})(L+1)}{\delta^2}\cdot\log\frac{m}{\eta}.
}
\]
With probability at least $1-\eta$, the output in Algorithm~\ref{Alg:inexact:obj} %
satisfies \Jie{$|\widehat{\Psi}(\lambda) - \Psi(\lambda)|\le \delta$.}
\Jie{
When using RT-MLMC subgradient estimator~\eqref{Eq:RTMLMC:scheme} in the BSMD step and RT-MLMC objective estimator~\eqref{Eq:RTMLMC:est}, the sample complexity from $\hP$ is $\tO\left(H_{\lambda,\Reg,B}K_{\lambda,\Reg,B}\exp(2K_{\lambda,\Reg,B})\cdot \delta^{-2}\right)$ and that from $\bQ_{x,\Reg}$ is $\tO\left(H_{\lambda,\Reg,B}K_{\lambda,\Reg,B}^2\exp(2K_{\lambda,\Reg,B})\cdot \delta^{-2}\right)$,
 with $\tO(\cdot)$ hiding constants depending on $L_f, \theta_0, \kappa$ and linearly depending on $(\log\frac{\lambda\Reg}{\delta})^2, (\log\frac{1}{\eta})^2$.
}
\end{proposition}

Next, we provide the convergence analysis for Algorithm~\ref{Alg:outer:bisection}.
\begin{theorem}
\label{Proposition:complexity:bisection}
Let $\eta\in(0,1)$.
Assume \Jie{Assumption~\ref{Assumption:throughout:loss} holds} and $0<\lambda_{l}\le \lambda^*\le \lambda_u<\infty$. Specify hyper-parameters in Algorithm~\ref{Alg:outer:bisection} as
\[
\Jie{
T' = \left\lceil \log_2\big(\frac{4L_{\Psi}(\lambda_u - \lambda_l)}{\delta}\big)\right\rceil,\quad 
\eta' = \frac{\eta}{3 + 2T'},\quad 
L_{\Psi}=\orho + \frac{B}{\lambda_{l}}\left[1 + \exp(K_{\lambda_l,\Reg,B})\right].
}
\]
Suppose there exists an oracle $\Jie{\widehat{\Psi}}$ such that for any $\lambda>0$, it estimates $\Psi$ defined in \eqref{Eq:F:lambda} with accuracy level $\delta/4$ with probability at least $1-\eta'$,
then with probability at least $1-\eta$, Algorithm~\ref{Alg:outer:bisection} finds the optimal multiplier with accuracy level $\delta$ (i.e., it finds $\lambda$ such that $\Jie{\Psi(\lambda) - \min\limits_{\lambda_{l}\le \lambda\le \lambda_u}~\Psi(\lambda)\le \delta}$) %
by calling the inexact oracle $\Jie{\widehat{\Psi}}$ for $\tO(K_{\lambda_l,\Reg,B})$ times, \Jie{where $\tO(\cdot)$ hides constants depending on $\orho$ and linearly depending on $\log\frac{\lambda_u-\lambda_l}{\delta}$ and $\log\frac{B}{\lambda_l}$.}
\end{theorem}

\Jie{
\begin{remark}[Selection of $\lambda_u$ and $\lambda_l$]
Algorithm~\ref{Alg:outer:bisection} requires the upper and lower bounds $\lambda_u,\lambda_l$ on the optimal Lagrange multiplier $\lambda^*$ as inputs.
Under Assumption~\ref{Assumption:throughout:loss}, we have a theoretical upper bound $\lambda_u:=\orho^{-1}B$~(see the proof in Lemma~\ref{Lemma:last} in Appendix~\ref{proof:Sec:inner:iteration:step:2}).
For the lower bound $\lambda_l$, it can be shown that as long as the condition in Theorem~\ref{Theorem:strong:duality}\ref{Theorem:strong:duality:IV} does not hold for $f_{\theta}(\cdot)$ for any $\theta$, a valid lower bound $\lambda_l>0$ exists.
Unfortunately, deriving the closed-form expression of $\lambda_l$ is infeasible.
In practice, choosing an excessively small multiplier values may lead to solutions that are too conservative (where the Sinkhorn distance constraint becomes nearly unbinding). To mitigate this, we recommend empirical tuning of $\lambda_l$ to ensure it remains sufficiently bounded away from zero, thereby avoiding degenerate cases.
For examples in Section~\ref{Sec:newsvendor} and \ref{Sec:portfolio}, we set \Jie{$\lambda_l=0.01$ and $\lambda_u=500$}.
\QEG
\end{remark}
}
Combining Proposition~\ref{Proposition:com:estimate:optval} and Theorem~\ref{Proposition:complexity:bisection}, the sample complexity from $\hP$ for obtaining a $\delta$-optimal solution of \eqref{Eq:reformulate:statistical:learning} with high probability is \Jie{$
\tO\left(H_{\lambda,\Reg,B}K_{\lambda,\Reg,B}^2\exp(2K_{\lambda,\Reg,B})\cdot \delta^{-2}\right)
$, and sample complexity from $\bQ_{x,\Reg}$ is $
\tO\left(H_{\lambda,\Reg,B}K_{\lambda,\Reg,B}^3\exp(2K_{\lambda,\Reg,B})\cdot \delta^{-2}\right)
$.}

\begin{remark}[Comparison with Empirical Risk Minimization]
The minimax lower bound of sample complexity from $\hP$ for obtaining a $\delta$-optimal solution from the empirical risk minimization~(ERM) $\inf_{\theta\in\Theta}~\bE_{x\sim\hP}[f_{\theta}(x)]$ with a convex loss function $f_{\theta}(z)$ (regardless of the smoothness assumption) is $\cO(\delta^{-2})$~\citep{Blair85}.
\Jie{
The sample complexity from $\hP$ of solving the Sinkhorn DRO model matches with its ERM counterpart, differing only by a (near-)constant factor in terms of error tolerance $\delta$.
However, we highlight that the constant factor has non-negligible dependence on related parameters $\lambda_l, \Reg, B$.
}
\QEG
\end{remark}

\begin{remark}[Comparison with Wasserestein DRO]
Recall from Table~\ref{Tab:compare:DRO} that
Wasserstein DRO is \Jie{computationally efficient to solve} for a restricted family of loss functions (Table~\ref{Tab:compare:DRO}). 
Specially, Wasserstein DRO with $\hP = \frac{1}{n}\sum_{i=1}^n\delta_{\hx_i}$ can be formulated as a minimax problem
\begin{equation*}
\min_{\theta\in\Theta, \lambda\ge0}~\max_{z_i\in\mathbb{R}^d, i\in[n]}~
\lambda\rho + \frac{1}{n}\sum_{i=1}^n\big[ 
f_{\theta}(z_i) - \lambda c(\hx_i, z_i)
\big].
\end{equation*}
When $f_{\theta}(z)$ is not piecewise concave in $z$, the above problem generally reduces to the convex-non-concave saddle point problem, whose global optimality is difficult to obtain.
\Jie{
In comparison, we provided complexity guarantees for solving Sinkhorn DRO model for a broader class of loss functions.
}

\Jie{
We also remark that the regularization parameter $\Reg$ is treated as a fixed intentional design choice when solving Sinkhorn DRO formulation.
Our goal is not to use it as a computational approximation of Wasserstein DRO, even though it converges to Wasserstein DRO as $\Reg\to0$.
Otherwise the constant part in our complexity bounds will explode to infinity.
}
\QEG
\end{remark}

Finally, we compare our algorithm design and analysis with existing references on CSO below.
\begin{remark}[Comparison with {{\citep{hu2020sample, hu2021biasvar, Yifan20}}}]\label{Remark:compare:hu}
\Jie{
Our algorithm and analysis differ from existing references on CSO \citep{hu2020sample, hu2021biasvar, Yifan20} in several aspects.
}

\Jie{
Recall that we apply BSMD with the RT-MLMC subgradient estimator to solve \eqref{Eq:F:lambda}. While the BSAA approach proposed in \citep{hu2020sample} is applicable, it is less efficient (see Remark~\ref{Remark:BSAA}). The BSMD method using the SG subgradient estimator from \citep{Yifan20} can also be applied, but it results in worse sample complexity from $\bQ_{x,\Reg}$. 
}

\Jie{
Although one can apply \citet{hu2021biasvar} to consider BSMD with RT-MLMC or vanilla MLMC (V-MLMC) gradient estimators, their analysis focuses on unconstrained optimization with $S_f$-Lipschitz smooth loss functions, i.e., functions that are continuously differentiable with Lipschitz continuous gradients: Their analysis requires carefully tuning the parameters of the gradient estimators to balance the trade-off between bias and gradient variance, with stepsize that depends on the smoothness constant $S_f$~(see \citep[Theorem~C.1]{hu2021biasvar}). 
This greatly limits the applicability of their theoretical guarantees, 
such as the newsvendor problem and portfolio optimization examples in our numerical study.
In contrast, we demonstrate that the RT-MLMC estimator results in the same complexity bounds even for nonsmooth loss functions by appropriately balancing the trade-off between bias and subgradient second-order moment~(Lemma~\ref{Lemma:SGD:scheme:convex:nonsmooth}).
Notably, the V-MLMC estimator in \citep{hu2021biasvar} no longer has performance guarantees in nonsmooth case.
}

\Jie{
Another component of our algorithm requires estimating the objective in \eqref{Eq:F:lambda}. 
Reference~\citep{Yifan20} used the SG objective estimator that leads to worse sample complexity from $\bQ_{x,\Reg}$;
References~\citep{hu2021biasvar, Yifan20} did not investigate this problem.
Our work provided the RT-MLMC estimator to estimate the Sinkhorn robust optimization objective, demonstrating its superior efficiency~(Lemma~\ref{Proposition:MLMC:obj}).
}
\QEG
\end{remark}

\section{Applications}\label{Sec:numerical}

In this section, we apply our methodology to three applications: the newsvendor \Jie{problem}, mean-risk portfolio optimization, and adversarial classification.
We compare our model with three benchmarks: (i) the classical sample average approximation~(SAA) model; (ii) the Wasserstein DRO model; and (iii) the KL-divergence DRO model.
We choose the transport cost $c(\cdot, \cdot) = \|\cdot - \cdot\|_1$ for $1$-Wasserstein or $1$-Sinkhorn DRO model, and $c(\cdot,\cdot)=\frac{1}{2}\|\cdot-\cdot\|^2_2$ for $2$-Wasserstein or $2$-Sinkhorn DRO model.
Throughout this section, we take the reference measure $\nu$ in the Sinkhorn distance to be the Lebesgue measure.
The hyper-parameters are selected using \Jie{the holdout method following from \citep{Mohajerin18}.} 
\Jie{All experiments were conducted on a Mac mini computer with 24GB of memory and M4 Pro GPU with 20 cores running Python~3.9}. 
Further implementation details and experiments are included in Appendices~\ref{Sec:experiment:setup} and \ref{Appendix:add:exp}, respectively.

\subsection{Newsvendor \Jie{Problem}}\label{Sec:newsvendor}

Consider the following distributionally robust newsvendor \Jie{problem}:
\begin{equation*}
    \min_{\theta}\max_{\bP\in\mathbb{B}_{\rho,\Reg}(\hP)}~\mathbb{E}_{z\sim\bP}\big[k\theta - u\min(\theta,z)\big],
\end{equation*}
where the random variable $z$ stands for the random demand, whose empirical distribution $\hP$ consists of $n$ independent samples from the underlying data distribution; the decision variable $\theta$ represents the inventory level; and $k=5,u=7$ are constants corresponding to overage and underage costs, respectively.

In this experiment, we examine the performance of DRO models for various sample sizes $n\in\{10,30,100\}$ and under three different types of data distribution: (i) the exponential distribution with rate parameter $1$, 
(ii) the gamma distribution with shape parameter $2$ and scale parameter $1.5$, (iii) the equiprobable mixture of two truncated normal distributions $\mathcal{N}(\mu=1,\sigma=1,a=0,b=10)$ and $\mathcal{N}(\mu=6,\sigma=1,a=0,b=10)$.
We do not report the performance for $1$-Wasserstein DRO model in this example, because it is identical to the SAA approach~\citep[Remark~6.7]{Mohajerin18}.
As $2$-Wasserstein DRO is computationally intractable for this example, we solve the corresponding formulation by discretizing the support of the distributions.

We measure the out-of-sample performance of a solution $\theta$ based on training dataset $\mathcal{D}$ using the percentage of improvement (a.k.a., coefficient of prescriptiveness) in~\citep{bertsimas2020predictive}:
\begin{equation}
\Jie{\text{Prescriptiveness}(\theta) = \max\left(1 - \frac{J(\theta) - J^*}{J(\theta^{\text{SAA}}_{\mathcal{D}}) - J^*}, -1\right)\times100\%,}
\label{Eq:coef:pre}
\end{equation}
where $J^*$ denotes the true optimal value when the true distribution is known, $\theta^{\text{SAA}}_{\mathcal{D}}$ denotes the decision from the SAA approach with dataset $\mathcal{D}$, and $J(\theta)$ denotes the expected loss of the solution $\theta$ under the true distribution, estimated through an SAA objective value with $10^5$ testing samples.
This coefficient is always bounded between $-100\%$ and $100\%$, and the higher this coefficient is, the better the solution's out-of-sample performance.

\begin{figure}[!ht]
    \centering\Jie{
     \subfigure[Exponential distribution]{\centering
     \includegraphics[width=0.3\textwidth]{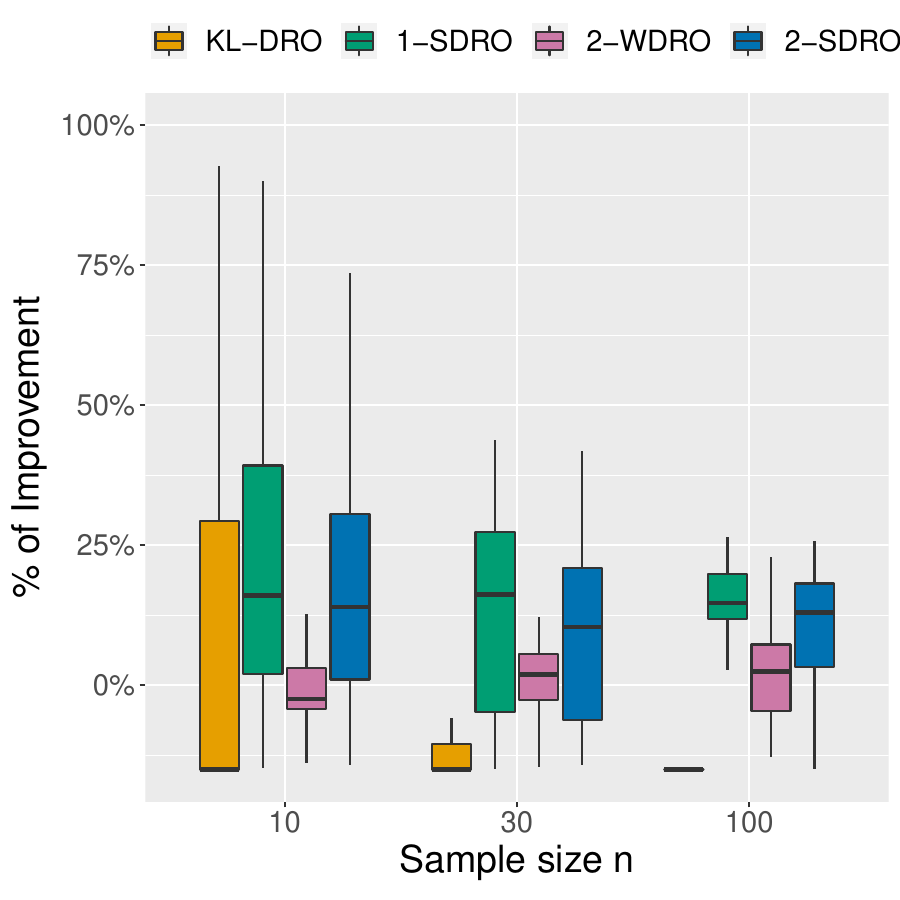}
        \label{fig:exp:newsvendor}
    }
    \hfill
    \subfigure[Gamma distribution]{\centering
    \includegraphics[width=0.3\textwidth]{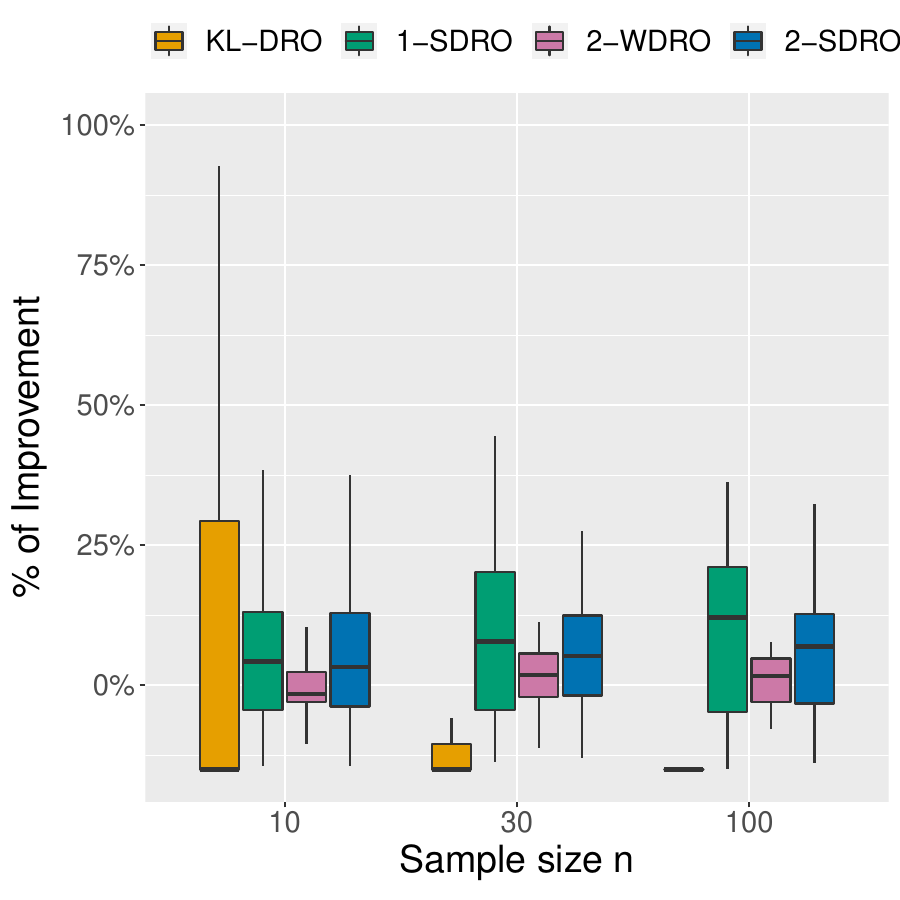}
    }
    \hfill
    \subfigure[Mixture of truncated normal distributions]{\centering
    \includegraphics[width=0.3\textwidth]{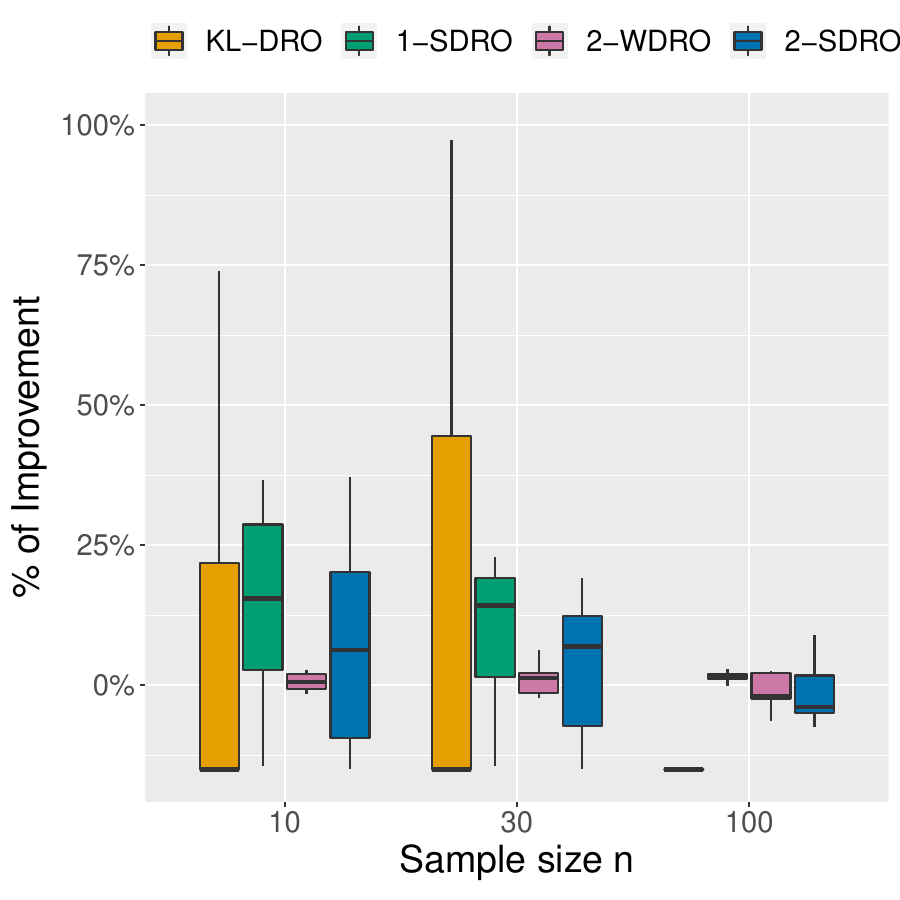}
    }
    \hfill
    \caption{
    Experiment results of the newsvendor problem for different sample sizes and different data distributions in box plots.
    }
    \label{fig:newsvendor:plot}
    }
\end{figure}

\begin{figure}[!ht]
    \centering\Jie{
     \subfigure[Exponential distribution]{\centering
     \includegraphics[width=0.3\textwidth]{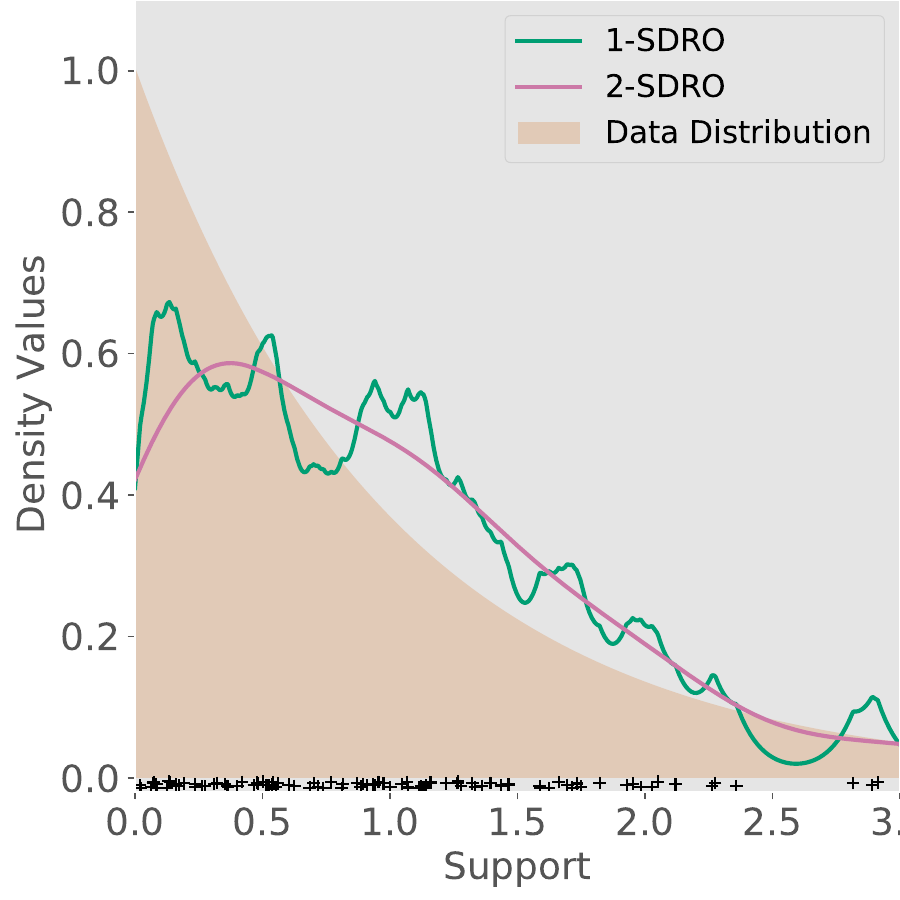}
    }
    \hfill
    \subfigure[Gamma distribution]{\centering
    \includegraphics[width=0.3\textwidth]{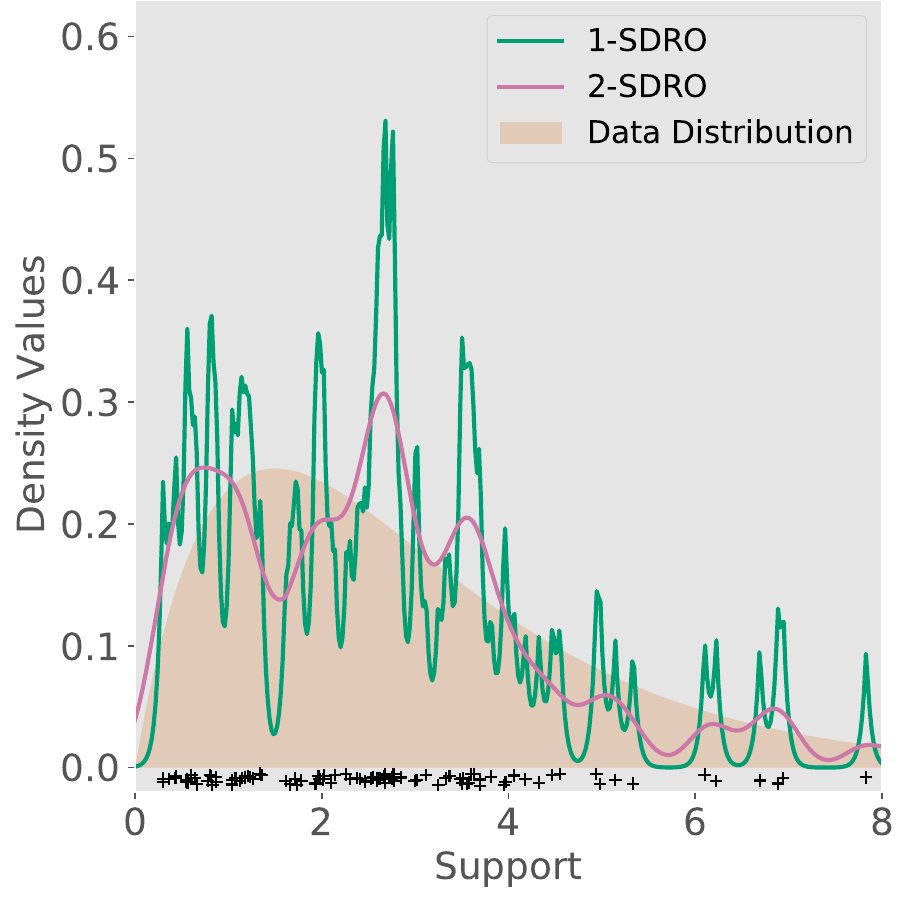}
    }
    \hfill
    \subfigure[Mixture of truncated normal distributions]{\centering
    \includegraphics[width=0.3\textwidth]{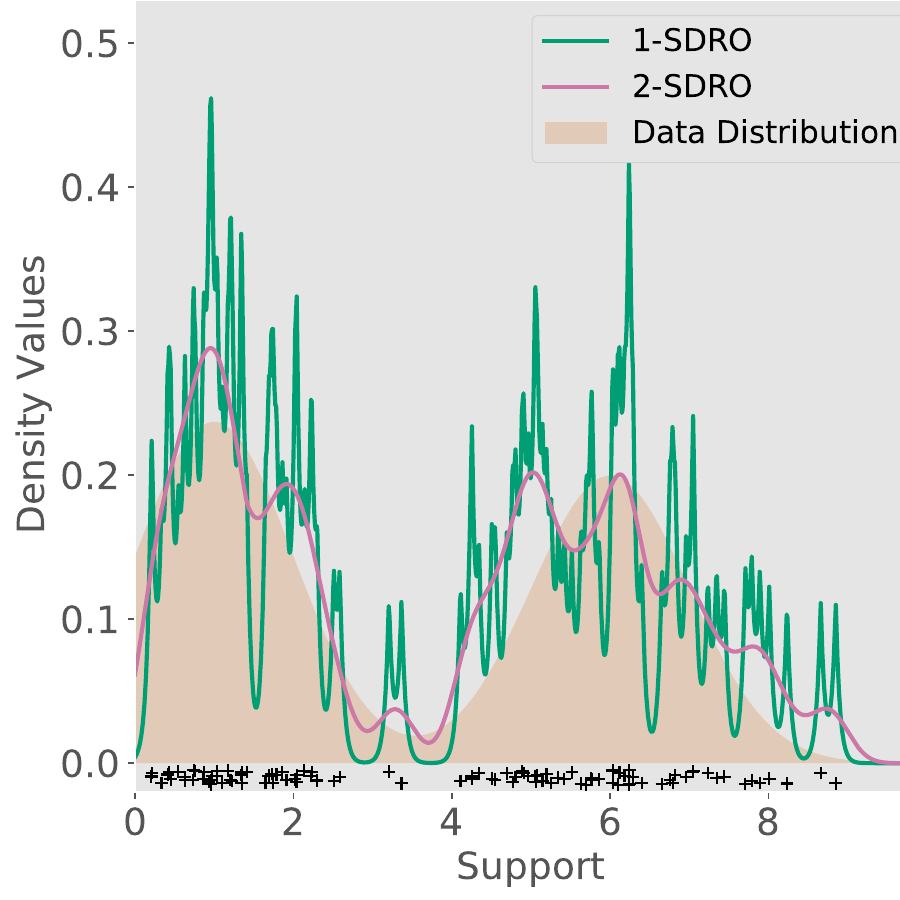}
    }
    \hfill
    \caption{
    Plots for the density of worst-case distributions generated by the $1$-SDRO or $2$-SDRO model for newsvendor problem with different data distributions.
    }
    \label{fig:newsvendor:plot:LFD}
    }
\end{figure}

\begin{figure}[!ht]   \centering
    \includegraphics[width=\linewidth]{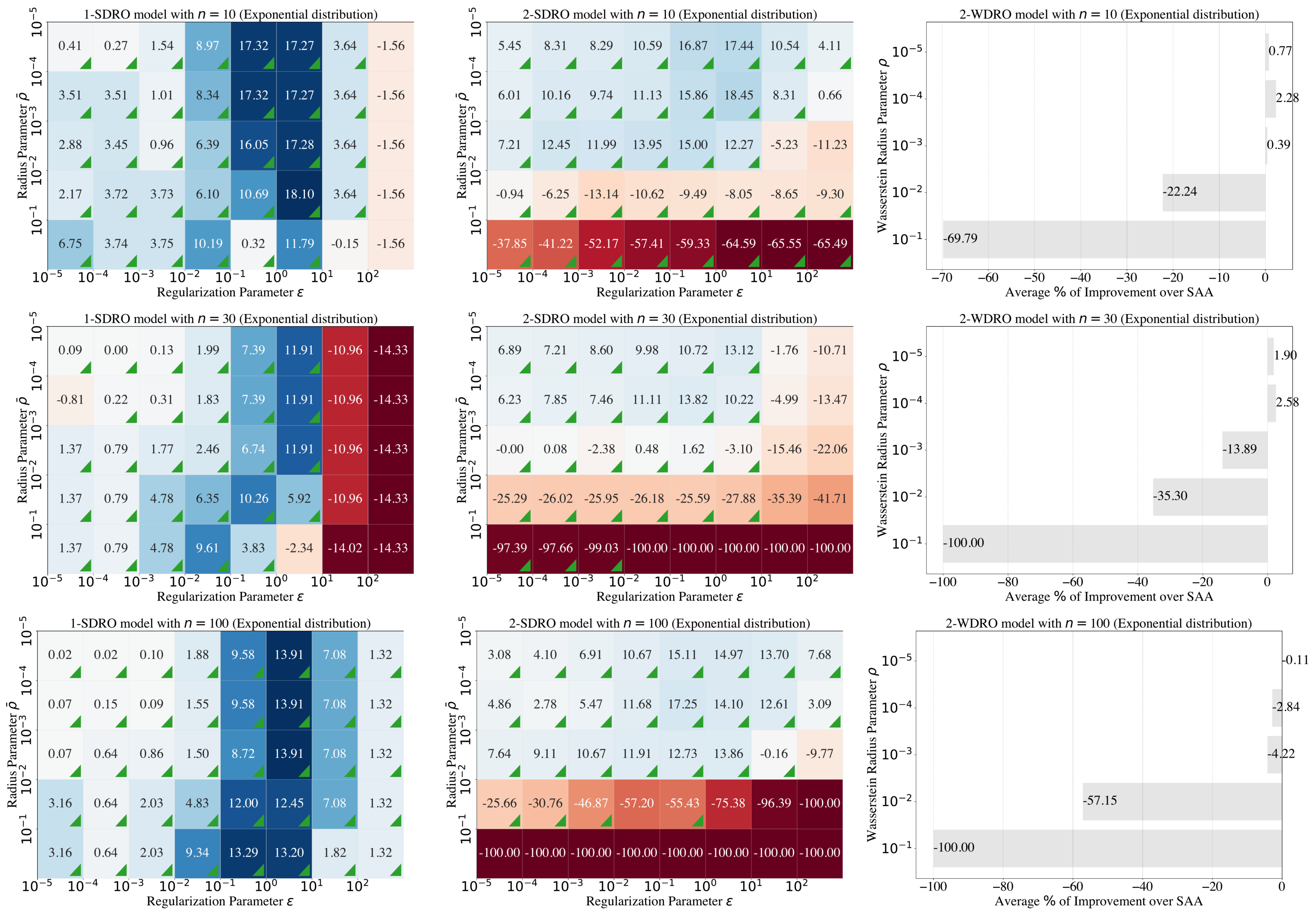}
    \caption{\Jie{
    Experiment results of the newsvendor problem for exponential data distribution.
    Subplots from different rows correspond to different training sample sizes $n\in\{10,30,100\}$.
    Subplots from the first and second columns correspond to the heatmap plot of the coefficient of prescriptiveness for $1$-SDRO and $2$-SDRO models with different radius and regularization parameters, and the subplots from the last column correspond to the histogram plot of the coefficient of prescriptiveness for $2$-WDRO model with different radius parameters.
    Each instance is taken the average of the simulation results over $50$ independent trials.
    For SDRO models, we add a green triangle for each radius-regularization combination that outperforms the corresponding WDRO models with the same radius choice.
   }
    }
    \label{fig:newsvendor:exponential:merged}
\end{figure}

We report the box-plots of the coefficients of prescriptiveness in Fig.~\ref{fig:newsvendor:plot} using \Jie{$500$} independent trials.
We find that either $1$-SDRO or $2$-SDRO model achieve the best out-of-sample performance over all sample sizes and data distributions listed, as it consistently scores higher than other benchmarks in the box plots.
In contrast, the KL-DRO model does not achieve satisfactory performance, and sometimes even underperforms the SAA model.
While the $2$-WDRO model demonstrates some improvement over the SAA model, the $2$-SDRO model shows more clear improvement.
We plot the density of worst-case distributions for $1$-SDRO or $2$-SDRO model in Fig.~\ref{fig:newsvendor:plot:LFD}.
When specifying the data distribution as exponential, gamma, or Gaussian mixture, the corresponding worst-case distributions capture the shape of the ground truth distribution reasonably well, which partly explains why the Sinkhorn DRO model achieves superior performance when the data distribution is absolutely continuous.
\Jie{We report the coefficient of prescriptiveness of SDRO and WDRO models with different parameters when the data distribution is exponential in Fig.~\ref{fig:newsvendor:exponential:merged}~(plots for other distributions are presented in Appendix~\ref{Appendix:coeff:p:combination}).
We observe that there exists a large range of parameter choices of SDRO models that lead to the superior performance over the WDRO models, since almost each row of the heatmap includes many green triangles.
This fact justifies the benefits of adding entropic regularization.
}

\subsection{Mean-risk Portfolio Optimization}\label{Sec:portfolio}

\begin{figure}[!ht]
 \centering\Jie{
     \subfigure[fixing data dimension $d=30$ and varying sample size $n\in\{30, 50,100,150,200,400\}$]{\centering
     \includegraphics[width=0.75\textwidth]{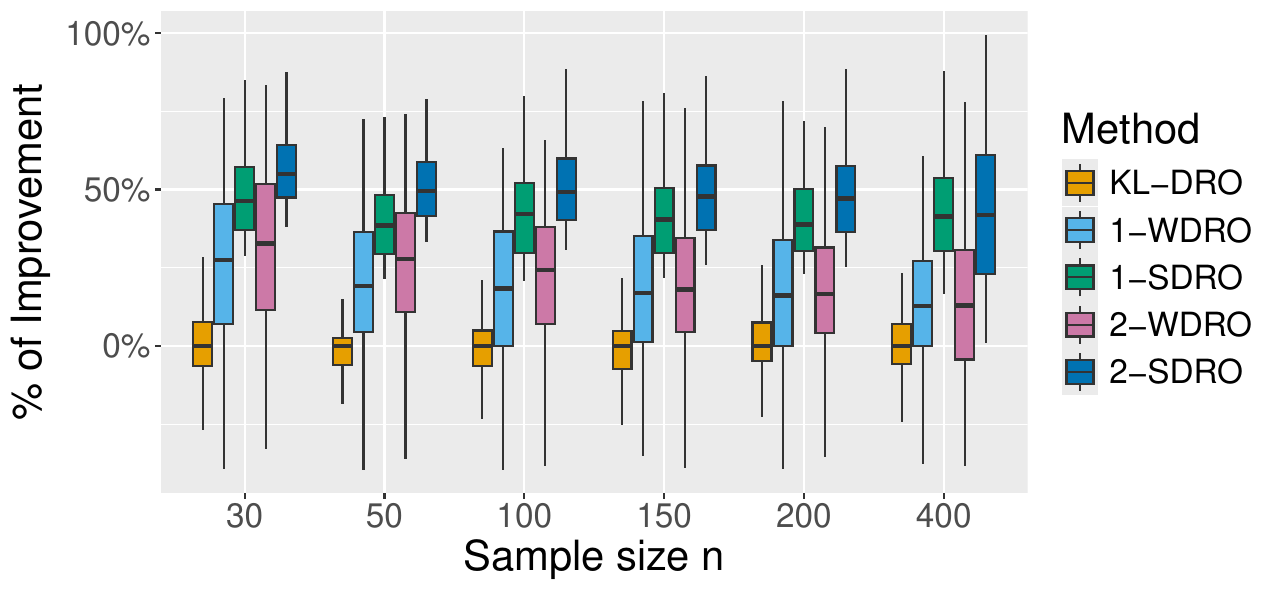}
    }
    \hfill
    \subfigure[fixing sample size $n=100$ and varying data dimension $d\in\{5, 10, 20, 40, 80, 100\}$]{\centering
    \includegraphics[width=0.75\textwidth]{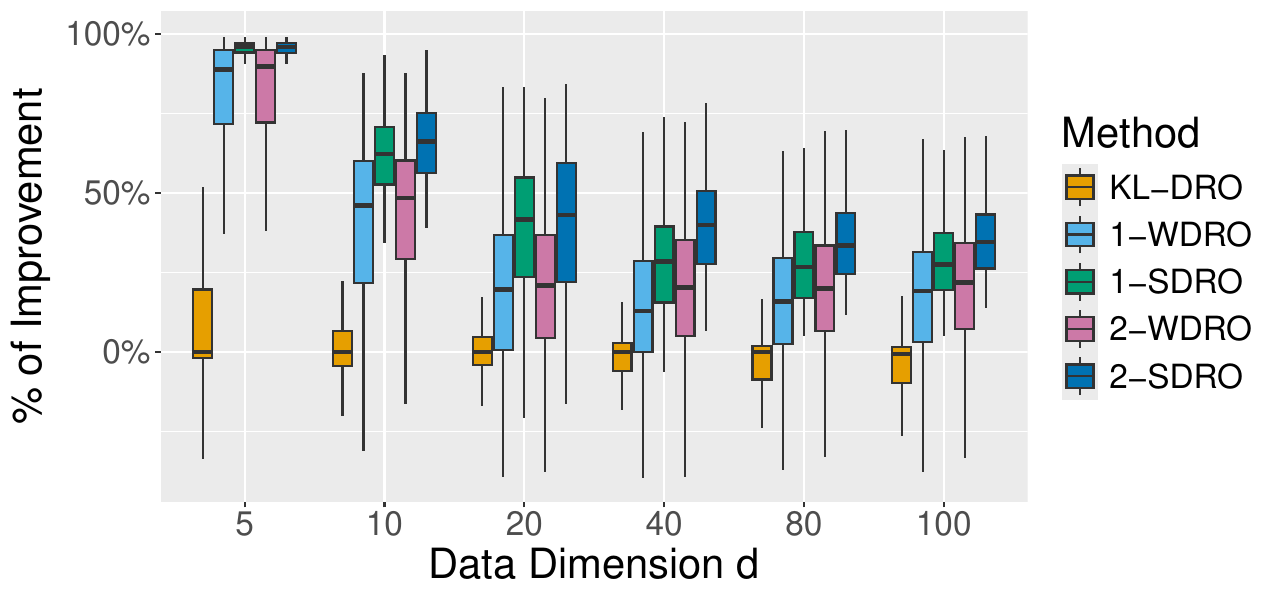}
    }
    \caption{
    Experiment results of the portfolio optimization problem for different sample sizes and dimensions in box plots.
    }}
    \label{fig:portfolio:plot}
\end{figure}

Consider the following distributionally robust mean-risk portfolio optimization problem
\begin{align*}
\min_{\theta}\max_{\bP\in\mathbb{B}_{\rho,\Reg}(\hP)}&\quad \mathbb{E}_{z\sim\bP}[-\theta\trans z]+\varrho\cdot\text{$\bP$-CVaR}_{\alpha}(-\theta\trans z)\\
\mbox{s.t.}&\quad \theta\in\Theta=\{\theta\in\mathbb{R}^d_+:\,\theta\trans 1=1\},
\end{align*}
where the random vector $z\in\mathbb{R}^d$ stands for the returns of assets; the decision variable $\theta\in\Theta$ represents the portfolio strategy that invests a certain percentage $\theta_i$ of the available capital in the $i$-th asset; and the term $\text{$\bP$-CVaR}_{\alpha}(-\theta\trans z)$ quantifies conditional value-at-risk \citep{rockafellar1999optimization}, i.e., the average of the $\alpha\times100\%$ worst portfolio losses under the distribution $\bP$.
We follow a similar setup as in \citet{Mohajerin18}.
Specifically, we set $\alpha=0.2,\varrho=10$.
The underlying true random return can be decomposed into a systematic risk factor $\psi\in\mathbb{R}$ and idiosyncratic risk factors $\epsilon\in\mathbb{R}^d$:
\[
z_i = \psi + \epsilon_i,\quad i=1,2,\ldots,d,
\]
where $\psi\sim\mathcal{N}(0,0.02)$ and $\epsilon_i\sim \mathcal{N}(i\times 0.03, i\times 0.025)$.
When solving the Sinkhorn DRO formulation, we take the Bregman divergence $D_{\omega}$ as the KL-divergence when performing BSMD algorithm in Algorithm~\ref{alg:BSMD:sampling}, allowing for efficient implementation~\citep{Blair85}.

We quantify the performance of a given solution using the same criterion defined in Section~\ref{Sec:newsvendor} and generate box plots using \Jie{$500$} independent trials.
Fig.~\ref{fig:portfolio:plot}a) reports the scenario where the data dimension $d=30$ is fixed and sample size $n\in\{30, 50, 100, 150, 200, 400\}$,
and Fig.~\ref{fig:portfolio:plot}b) reports the scenario where the sample size $n=100$ is fixed and the number of assets $d\in\{5, 10, 20, 40, 80, 100\}$.
We find that the KL-DRO model does not have competitive performance compared to other DRO models, especially as the data dimension $d$ increases. This is because the ambiguity set of KL-DRO model only takes into account those distributions sharing the same support as the nominal distribution, which seems to be restrictive, especially for high-dimensional settings.
While $1$-WDRO or $2$-WDRO model has better out-of-sample performance than the SAA model, the corresponding $1$-SDRO or $2$-SDRO model has clearer improvements, as it consistently scores higher in the box plots.
\Jie{Finally, we show the coefficient of prescriptiveness of WDRO/SDRO models with different parameters for the instance $(n,d)=(30,30)$ in Fig.~\ref{fig:portfolio:heatmap:init}~(other instances can be found in Appendix~\ref{Appendix:coeff:p:combination}).
Similar as in the newsvendor problem, we observe that there exists a large number of parameter choices of SDRO models that outperform WDRO models, as indicated by many green triangles in heatmap plots.}

\begin{figure}[!ht]
    \centering\Jie{
     \centering
     \includegraphics[width=0.8\textwidth]{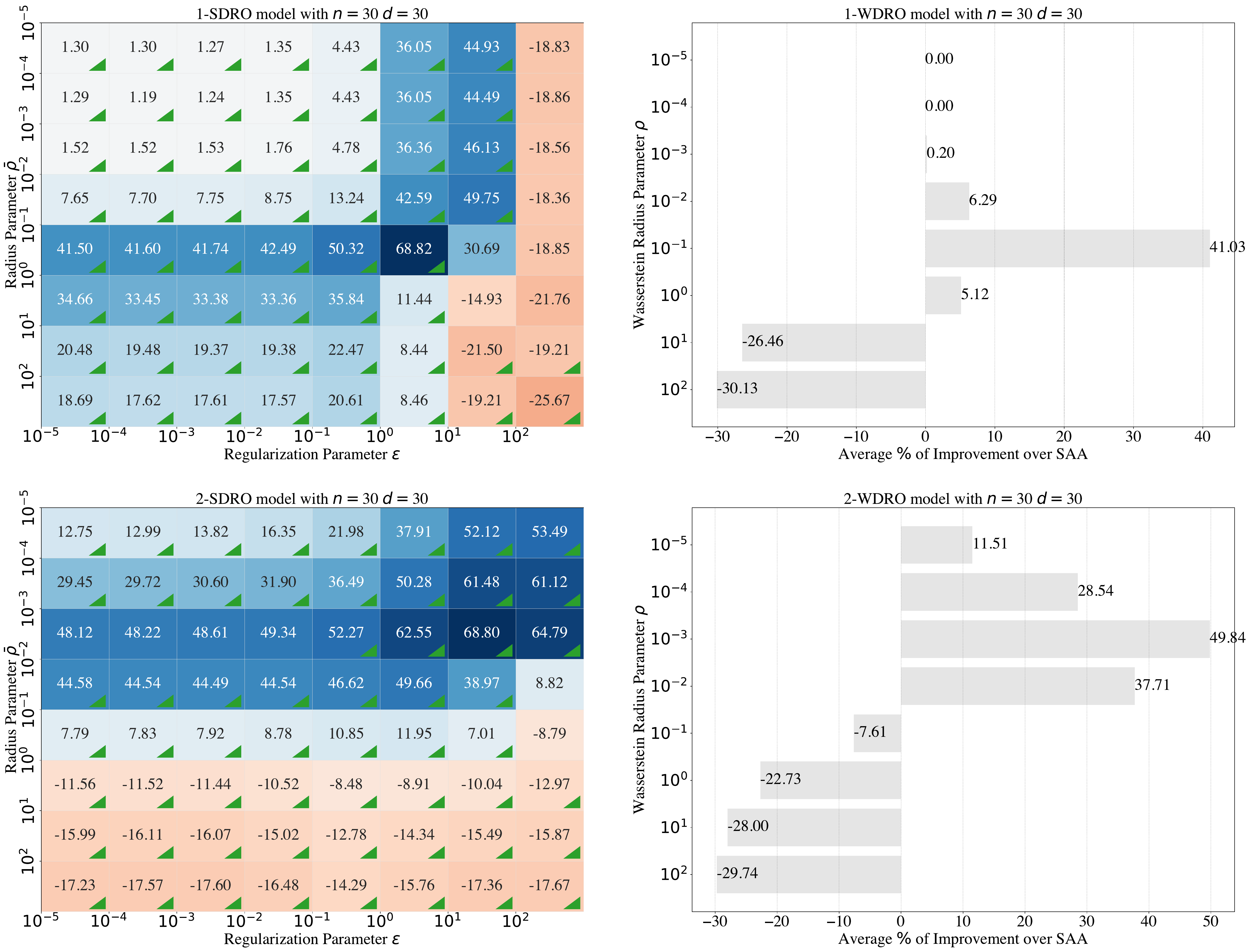}
    \caption{
    Experiment results of the portfolio optimization model 
    with $(n,d)=(30,30)$ in heatmaps.
    The four subplots correspond to the heatmap plot of the coefficient of prescripiveness for 1-SDRO, 1-WDRO, 2-SDRO, and 2-WDRO models with varying parameters.
    For SDRO models, we add a green triangle for each radius-regularization combination that outperforms the corresponding WDRO models with the same radius choice.
    }
    \label{fig:portfolio:heatmap:init}
    }
\end{figure}

\subsection{Adversarial Multi-class Logistic Regression}
\label{Sec:classification}
\begin{figure}[!ht]
    \centering
    \includegraphics[width=\textwidth]{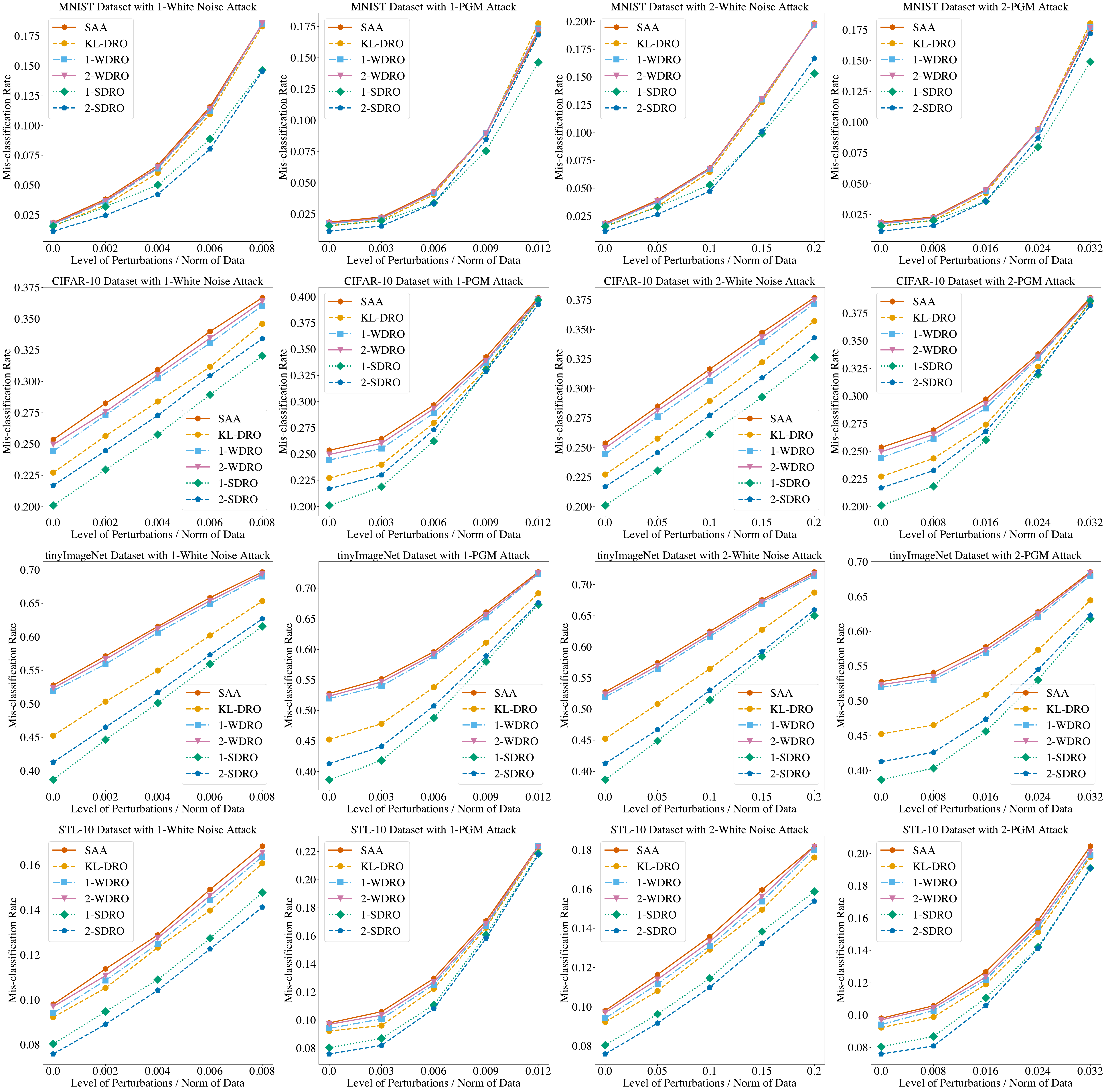}
    \caption{Results of adversarial training on various image datasets with different types of adversarial attack.
    From left to right, the figures correspond to 
    \Jie{(a) white Laplacian noise attack; (b) $\ell_1$-norm PGD attack;
    (c) white Gaussian noise attack; and (d) $\ell_2$-norm PGD attack.}
    From top to bottom, the figures correspond to \Jie{(a) MNIST dataset; (b) CIFAR-10 dataset; (c) tinyImageNet dataset; and (d) STL-10 dataset}.}
    \label{fig:adv:cifar10}
\end{figure}

\Jie{
Adversarial machine learning is an emerging topic in artificial intelligence, aiming to develop models that are robust against (potentially adversarial) data perturbations.
It has been observed that small perturbations to the data can cause well-trained machine learning models to produce unexpectedly inaccurate predictions~\citep {goodfellow2014explaining}.
In real applications involving high-stake environments, such as self-driving and automated tumor detection, ensuring model robustness is essential for reliability and safety.
Among existing approaches that produce robust machine learning models~\citep{papernot2017practical, papernot2016limitations, papernot2016distillation, rozsa2016towards, sinha2018certifiable, he2017adversarial, madry2017towards, tramer2017ensemble}, 
stands out as a particularly effective method and provides certifiable robustness~\citep{sinha2018certifiable}.
}
In this subsection, \Jie{we examine the performance of various DRO approaches for multi-class logistic regression with data perturbations.}
Given a feature vector $x\in\mathbb{R}^d$ and its label $y\in[C]$, we denote ${\bm y}\in\{0,1\}^C$ as the corresponding one-hot label vector, and define the negative likelihood loss
\[
h_B(x, {\bm y}) = -{\bm y}\trans B\trans x + \log\big(1\trans e^{B\trans x}\big),
\]
where $B:=[w_1,\ldots,w_K]$ \Jie{denotes the parameters of the linear classifier}.
Let $\hP$ be the empirical distribution from training samples.
Since the testing samples may have slightly different data distributions than the training samples, the DRO model aims to solve the following optimization problem to mitigate the impact of \Jie{data perturbations}:
\[
\begin{aligned}
\min_{B}~\max_{\bP\in\mathbb{B}_{\rho,\Reg}(\hP)}&\quad \mathbb{E}_{(x, {\bm y})\sim \bP}\big[h_B(x, {\bm y})\big].
\end{aligned}
\]
\Jie{It is assumed that the data perturbation only happens for the feature vector $x$ but not the label ${\bm y}$.}

We conduct experiments on four large-scale datasets: 
\Jie{MNIST~\citep{lecun1998gradient},
CIFAR-10~\citep{krizhevsky2009learning}, 
tinyImageNet~\citep{tinyimagenet},
and STL-10~\citep{coates2011analysis}.}
We pre-process these datasets using the ResNet-18 network~\citep{he2016deep} pre-trained on the ImageNet dataset to extract linear features.
Since this network has learned a rich set of hierarchical features from the large and diverse ImageNet dataset, it typically extracts useful features for other image datasets.
We then add different types of perturbations to the testing datasets, such as \emph{white Laplace noise}, \emph{white Gaussian noise}, and \emph{$\ell_p$-norm adversarial projected gradient descent~(PGD) attacks~\citep{madry2017towards} with $p\in\{1,2\}$.
}
The level of perturbation is normalized by the averaged $\ell_2$ norm of the feature vectors from testing dataset.
See the detailed procedure for generating data perturbations and statistics on pre-processed datasets in Appendix~\ref{Appendix:adv:re}.
We use the mis-classification rate on testing dataset to measure the performance of the obtained classifers.

\Jie{
For baseline DRO models, we solve their Lagrangian relaxation, which adds the penalty of the statistical distance into the objective to ensure efficient implementation.
To make fair comparisons, we tune the penalty parameter for each method such that the $2$-Wasserstein distance between the nominal distribution and its perturbed one is controlled within $\varrho:=0.05\cdot\mathbb{E}[\|\bm x\|_2^2]$~(expectation taken with respect to the training dataset).
For SDRO methods, we fix $\Reg=0.1$ unless otherwise stated.
It is noteworthy that solving the Lagrangian relaxation of $2$-WDRO has global convergence guarantees only when the penalty parameter is sufficiently large~\citep{sinha2018certifiable}, which is not the case for this example. 
In general, solving $1$- or $2$-WDRO model reduces to solving a convex-non-concave minimax game, and we try gradient descent ascent~\citep[Algorithm~1]{sinha2018certifiable} as a heuristic.
}

\begin{figure}[!ht]
    \centering
    \includegraphics[width=\linewidth]{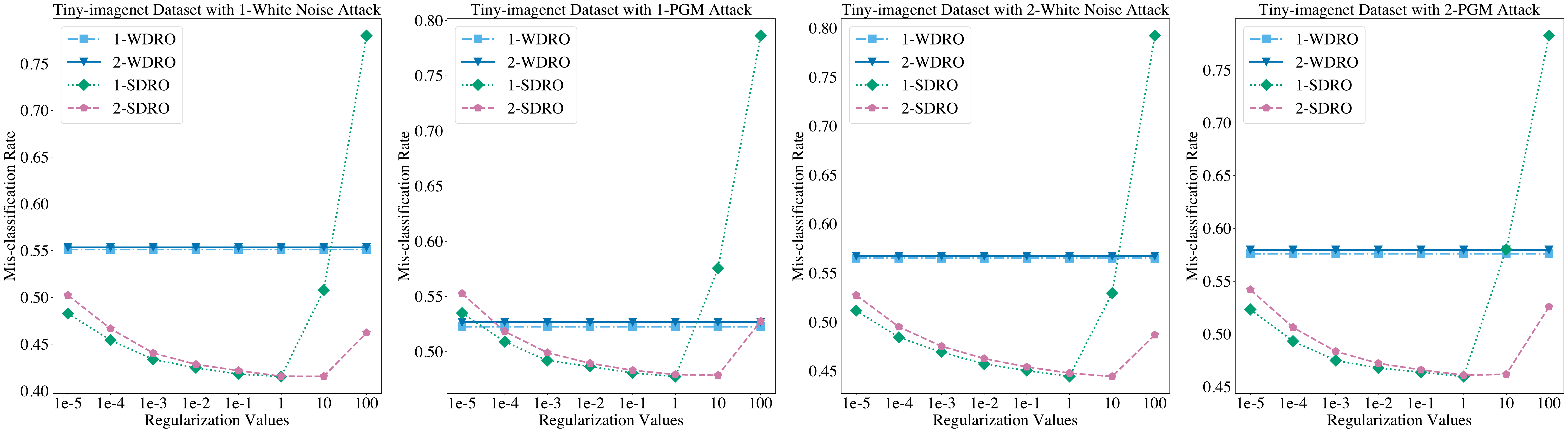}
    \caption{\Jie{Experiment results of the adversarial classification problem with tinyImagenet dataset.
    The subplots from left to right correspond to the misclassification errors of SDRO and WDRO models with different types of adversarial attack.
    For SDRO models we vary the regularization parameter $\Reg$.
    }}
    \label{fig:adversarial:ablation}
\end{figure}
Fig.~\ref{fig:adv:cifar10} presents the classification results for different types of adversarial attacks with varying levels of perturbations on the datasets. We observe that as the level of perturbations on the testing samples increases, all methods tend to perform worse. However, both the $1$-SDRO and $2$-SDRO models show a slower trend of increasing error rates than other benchmarks across all types of adversarial attacks and all datasets. This suggests that SDRO models can be a suitable choice for adversarial robust training.
\Jie{
Fig.\ref{fig:adversarial:ablation} presents the misclassification rates of SDRO models under different types of adversarial attacks on the TinyImageNet dataset (additional results are provided in Appendix~\ref{Appendix:coeff:p:combination}). For baseline comparison, we also report the performance of WDRO models. In all subplots, the perturbation levels (normalized by the data norm) are fixed at $0.004$, $0.006$, $0.01$, and $0.016$ for the four adversarial attacks, while the penalty parameter is set to $\lambda=10$ for both WDRO and SDRO models. The figure demonstrates that SDRO outperforms WDRO across a wide range of regularization parameters.
}

\section{Concluding Remarks}\label{Sec:conclusion}
In this paper, we investigated a new distributionally robust optimization framework based on the Sinkhorn distance.
By developing a strong dual reformulation and a biased stochastic mirror descent algorithm, we have shown that the resulting problem is efficient to solve under mild assumptions.
Analysis of the worst-case distribution indicates that Sinkhorn DRO hedges a more reasonable set of adverse scenarios and is thus less conservative than Wasserstein DRO.
Extensive numerical experiments  demonstrated that Sinkhorn DRO is a promising candidate for modeling distributional ambiguities in decision-making under uncertainty.

In the meantime, several topics worth investigating are left for future work. 
\Jie{
For example, it is desirable to study the statistical performance guarantees under suitable choices of hyper-parameters.
It is also of research interest to develop and analyze optimization algorithms with less restrictive assumptions and sharper complexity bounds.
Exploring and discovering the benefits of Sinkhorn DRO in other applications may also be of future interest.
}

\theendnotes
\bibliographystyle{informs2014} 
\bibliography{shortbib}

\begin{thebibliography}{125}
\providecommand{\natexlab}[1]{#1}
\providecommand{\url}[1]{\texttt{#1}}
\providecommand{\urlprefix}{URL }

\bibitem[{Agrawal et~al.(2012)Agrawal, Ding, Saberi, \protect\BIBand{}
  Ye}]{agrawal2012price}
Agrawal S, Ding Y, Saberi A, Ye Y (2012) Price of correlations in stochastic
  optimization. \emph{Operations Research} 60(1):150--162.

\bibitem[{Altschuler et~al.(2017)Altschuler, Weed, \protect\BIBand{}
  Rigollet}]{Altschuler17}
Altschuler J, Weed J, Rigollet P (2017) Near-linear time approximation
  algorithms for optimal transport via sinkhorn iteration. \emph{Advances in
  Neural Information Processing Systems}, 1961–1971.

\bibitem[{ApS(2021)}]{mosek}
ApS M (2021) Mosek modeling cookbook 3.2.3.
  \url{https://docs.mosek.com/modeling-cookbook/index.html#}.

\bibitem[{Asmussen \protect\BIBand{} Glynn(2007)}]{asmussen2007stochastic}
Asmussen S, Glynn PW (2007) \emph{Stochastic simulation: algorithms and
  analysis}, volume~57 (Springer Science \& Business Media).

\bibitem[{Azizian et~al.(2023)Azizian, Iutzeler, \protect\BIBand{}
  Malick}]{azizian2022regularization}
Azizian W, Iutzeler F, Malick J (2023) Regularization for wasserstein
  distributionally robust optimization. \emph{ESAIM: Control, Optimisation and
  Calculus of Variations} 29:33.

\bibitem[{Bacharach(1965)}]{bacharach1965estimating}
Bacharach M (1965) Estimating nonnegative matrices from marginal data.
  \emph{International Economic Review} 6(3):294--310.

\bibitem[{Bai et~al.(2020)Bai, Wu, \protect\BIBand{}
  Ozgur}]{bai2020information}
Bai Y, Wu X, Ozgur A (2020) Information constrained optimal transport: From
  talagrand, to marton, to cover. \emph{2020 IEEE International Symposium on
  Information Theory (ISIT)}, 2210--2215.

\bibitem[{Bayraksan \protect\BIBand{} Love(2015)}]{bayraksan2015data}
Bayraksan G, Love DK (2015) Data-driven stochastic programming using
  phi-divergences. \emph{The Operations Research Revolution}, 1--19 (INFORMS).

\bibitem[{Ben-Tal et~al.(2013)Ben-Tal, den Hertog, De~Waegenaere, Melenberg,
  \protect\BIBand{} Rennen}]{Ben13}
Ben-Tal A, den Hertog D, De~Waegenaere A, Melenberg B, Rennen G (2013) Robust
  solutions of optimization problems affected by uncertain probabilities.
  \emph{Management Science} 59(2):341--357.

\bibitem[{Bertsimas \protect\BIBand{} Kallus(2020)}]{bertsimas2020predictive}
Bertsimas D, Kallus N (2020) From predictive to prescriptive analytics.
  \emph{Management Science} 66(3):1025--1044.

\bibitem[{Bertsimas et~al.(2006)Bertsimas, Natarajan, \protect\BIBand{}
  Teo}]{bertsimas2006persistence}
Bertsimas D, Natarajan K, Teo CP (2006) Persistence in discrete optimization
  under data uncertainty. \emph{Mathematical programming} 108(2):251--274.

\bibitem[{Bertsimas et~al.(2019)Bertsimas, Sim, \protect\BIBand{}
  Zhang}]{Bertsimas19}
Bertsimas D, Sim M, Zhang M (2019) Adaptive distributionally robust
  optimization. \emph{Management Science} 65(2):604--618.

\bibitem[{Blackwell \protect\BIBand{} Ryll-Nardzewski(1963)}]{blackwell1963non}
Blackwell D, Ryll-Nardzewski C (1963) Non-existence of everywhere proper
  conditional distributions. \emph{The Annals of Mathematical Statistics}
  34(1):223--225.

\bibitem[{Blanchet et~al.(2022{\natexlab{a}})Blanchet, Chen, \protect\BIBand{}
  Zhou}]{blanchet2018distributionallyportfolio}
Blanchet J, Chen L, Zhou XY (2022{\natexlab{a}}) Distributionally robust
  mean-variance portfolio selection with wasserstein distances.
  \emph{Management Science} 68(9):6382--6410.

\bibitem[{Blanchet et~al.(2019{\natexlab{a}})Blanchet, Glynn, Yan,
  \protect\BIBand{} Zhou}]{blanchet2019multivariate}
Blanchet J, Glynn PW, Yan J, Zhou Z (2019{\natexlab{a}}) Multivariate
  distributionally robust convex regression under absolute error loss.
  \emph{Advances in Neural Information Processing Systems}, volume~32,
  11817--11826.

\bibitem[{Blanchet \protect\BIBand{} Kang(2020)}]{blanchet2020semi}
Blanchet J, Kang Y (2020) Semi-supervised learning based on distributionally
  robust optimization. \emph{Data Analysis and Applications 3: Computational,
  Classification, Financial, Statistical and Stochastic Methods} 5:1--33.

\bibitem[{Blanchet et~al.(2019{\natexlab{b}})Blanchet, Kang, \protect\BIBand{}
  Murthy}]{blanchet2019robust}
Blanchet J, Kang Y, Murthy K (2019{\natexlab{b}}) Robust wasserstein profile
  inference and applications to machine learning. \emph{Journal of Applied
  Probability} 56(3):830--857.

\bibitem[{Blanchet \protect\BIBand{} Murthy(2019)}]{blanchet2019quantifying}
Blanchet J, Murthy K (2019) Quantifying distributional model risk via optimal
  transport. \emph{Mathematics of Operations Research} 44(2):565--600.

\bibitem[{Blanchet et~al.(2021)Blanchet, Murthy, \protect\BIBand{}
  Nguyen}]{blanchet2021statistical}
Blanchet J, Murthy K, Nguyen VA (2021) Statistical analysis of wasserstein
  distributionally robust estimators. \emph{Tutorials in Operations Research:
  Emerging Optimization Methods and Modeling Techniques with Applications},
  227--254 (INFORMS).

\bibitem[{Blanchet et~al.(2022{\natexlab{b}})Blanchet, Murthy,
  \protect\BIBand{} Si}]{blanchet2019confidence}
Blanchet J, Murthy K, Si N (2022{\natexlab{b}}) Confidence regions in
  wasserstein distributionally robust estimation. \emph{Biometrika}
  109(2):295--315.

\bibitem[{Blanchet et~al.(2022{\natexlab{c}})Blanchet, Murthy,
  \protect\BIBand{} Zhang}]{blanchet2021optimal}
Blanchet J, Murthy K, Zhang F (2022{\natexlab{c}}) Optimal transport-based
  distributionally robust optimization: Structural properties and iterative
  schemes. \emph{Mathematics of Operations Research} 47(2):1500--1529.

\bibitem[{Blanchet \protect\BIBand{} Glynn(2015)}]{blanchet2015unbiased}
Blanchet JH, Glynn PW (2015) Unbiased monte carlo for optimization and
  functions of expectations via multi-level randomization. \emph{2015 Winter
  Simulation Conference (WSC)}, 3656--3667.

\bibitem[{Chang \protect\BIBand{} Lin(2011)}]{libsvmregression}
Chang CC, Lin CJ (2011) Libsvm: a library for support vector machines.
  \emph{ACM transactions on intelligent systems and technology (TIST)}
  2(3):1--27.

\bibitem[{Chen \protect\BIBand{} Paschalidis(2019)}]{chen2019selecting}
Chen R, Paschalidis IC (2019) Selecting optimal decisions via distributionally
  robust nearest-neighbor regression. \emph{Advances in Neural Information
  Processing Systems}.

\bibitem[{Chen et~al.(2020)Chen, Sun, \protect\BIBand{}
  Xu}]{chen2021decomposition}
Chen Y, Sun H, Xu H (2020) Decomposition and discrete approximation methods for
  solving two-stage distributionally robust optimization problems.
  \emph{Computational Optimization and Applications} 78(1):205--238.

\bibitem[{Chen et~al.(2022)Chen, Kuhn, \protect\BIBand{}
  Wiesemann}]{chen2018data}
Chen Z, Kuhn D, Wiesemann W (2022) Data-driven chance constrained programs over
  wasserstein balls. \emph{Operations Research} .

\bibitem[{Chen et~al.(2019)Chen, Sim, \protect\BIBand{} Xu}]{Chen19}
Chen Z, Sim M, Xu H (2019) Distributionally robust optimization with infinitely
  constrained ambiguity sets. \emph{Operations Research} 67(5):1328--1344.

\bibitem[{Cherukuri \protect\BIBand{}
  Cort{\'e}s(2019)}]{cherukuri2019cooperative}
Cherukuri A, Cort{\'e}s J (2019) Cooperative data-driven distributionally
  robust optimization. \emph{IEEE Transactions on Automatic Control}
  65(10):4400--4407.

\bibitem[{Coates \protect\BIBand{} Ng(2011)}]{coates2011analysis}
Coates A, Ng AY (2011) Analysis of large-scale visual recognition.
  \emph{Advances in neural information processing systems} 24:873--881.

\bibitem[{Cohen et~al.(2016)Cohen, Lee, Miller, Pachocki, \protect\BIBand{}
  Sidford}]{cohen2016geometric}
Cohen MB, Lee YT, Miller G, Pachocki J, Sidford A (2016) Geometric median in
  nearly linear time. \emph{Proceedings of the forty-eighth annual ACM
  symposium on Theory of Computing}, 9--21.

\bibitem[{Courty et~al.(2017)Courty, Flamary, Habrard, \protect\BIBand{}
  Rakotomamonjy}]{courty2017joint}
Courty N, Flamary R, Habrard A, Rakotomamonjy A (2017) Joint distribution
  optimal transportation for domain adaptation. \emph{Advances in Neural
  Information Processing Systems}.

\bibitem[{Courty et~al.(2014)Courty, Flamary, \protect\BIBand{}
  Tuia}]{courty2014domain2}
Courty N, Flamary R, Tuia D (2014) Domain adaptation with regularized optimal
  transport. \emph{Joint European Conference on Machine Learning and Knowledge
  Discovery in Databases}, 274--289.

\bibitem[{Courty et~al.(2016)Courty, Flamary, Tuia, \protect\BIBand{}
  Rakotomamonjy}]{courty2016optimal}
Courty N, Flamary R, Tuia D, Rakotomamonjy A (2016) Optimal transport for
  domain adaptation. \emph{IEEE Transactions on Pattern Analysis and Machine
  Intelligence} 39(9):1853--1865.

\bibitem[{Cover \protect\BIBand{} Thomas(2006)}]{Cover06}
Cover TM, Thomas JA (2006) \emph{Elements of Information Theory}
  (Wiley-Interscience).

\bibitem[{Cuturi(2013)}]{cuturi2013sinkhorn}
Cuturi M (2013) Sinkhorn distances: Lightspeed computation of optimal
  transport. \emph{Advances in neural information processing systems},
  volume~26, 2292--2300.

\bibitem[{Delage \protect\BIBand{} Ye(2010)}]{Delage10}
Delage E, Ye Y (2010) Distributionally robust optimization under moment
  uncertainty with application to data-driven problems. \emph{Operations
  Research} 58(3):595--612.

\bibitem[{Deming \protect\BIBand{} Stephan(1940)}]{deming1940least}
Deming WE, Stephan FF (1940) On a least squares adjustment of a sampled
  frequency table when the expected marginal totals are known. \emph{The Annals
  of Mathematical Statistics} 11(4):427--444.

\bibitem[{Doan \protect\BIBand{} Natarajan(2012)}]{doan2012complexity}
Doan XV, Natarajan K (2012) On the complexity of nonoverlapping multivariate
  marginal bounds for probabilistic combinatorial optimization problems.
  \emph{Operations research} 60(1):138--149.

\bibitem[{Duchi et~al.(2021)Duchi, Glynn, \protect\BIBand{} Namkoong}]{Duchi21}
Duchi JC, Glynn PW, Namkoong H (2021) Statistics of robust optimization: A
  generalized empirical likelihood approach. \emph{Mathematics of Operations
  Research} 0(0).

\bibitem[{Eckstein et~al.(2020)Eckstein, Kupper, \protect\BIBand{}
  Pohl}]{eckstein2020robust}
Eckstein S, Kupper M, Pohl M (2020) Robust risk aggregation with neural
  networks. \emph{Mathematical Finance} 30(4):1229--1272.

\bibitem[{Esfahani \protect\BIBand{} Kuhn(2018)}]{esfahani2018data}
Esfahani PM, Kuhn D (2018) Data-driven distributionally robust optimization
  using the wasserstein metric: Performance guarantees and tractable
  reformulations. \emph{Mathematical Programming} 171(1):115--166.

\bibitem[{Feng \protect\BIBand{} Schl{\"o}gl(2018)}]{feng2018model}
Feng Y, Schl{\"o}gl E (2018) Model risk measurement under wasserstein distance.
  \emph{arXiv preprint arXiv:1809.03641} .

\bibitem[{Fr{\'e}chet(1960)}]{frechet1960tableaux}
Fr{\'e}chet M (1960) Sur les tableaux dont les marges et des bornes sont
  donn{\'e}es. \emph{Revue de l'Institut international de statistique} 10--32.

\bibitem[{Gao(2022)}]{gao2020finitesample}
Gao R (2022) Finite-sample guarantees for wasserstein distributionally robust
  optimization: Breaking the curse of dimensionality. \emph{Operations
  Research} .

\bibitem[{Gao et~al.(2022)Gao, Chen, \protect\BIBand{}
  Kleywegt}]{gao2020wasserstein}
Gao R, Chen X, Kleywegt AJ (2022) Wasserstein distributionally robust
  optimization and variation regularization. \emph{Operations Research} .

\bibitem[{Gao \protect\BIBand{} Kleywegt(2022)}]{gao2016distributionally}
Gao R, Kleywegt A (2022) Distributionally robust stochastic optimization with
  wasserstein distance. \emph{Mathematics of Operations Research} .

\bibitem[{Genevay et~al.(2016)Genevay, Cuturi, Peyr\'{e}, \protect\BIBand{}
  Bach}]{aude2016stochastic}
Genevay A, Cuturi M, Peyr\'{e} G, Bach F (2016) Stochastic optimization for
  large-scale optimal transport. \emph{Advances in Neural Information
  Processing Systems}, volume~29.

\bibitem[{Genevay et~al.(2018)Genevay, Peyre, \protect\BIBand{}
  Cuturi}]{Aude18a}
Genevay A, Peyre G, Cuturi M (2018) Learning generative models with sinkhorn
  divergences. \emph{Proceedings of the Twenty-First International Conference
  on Artificial Intelligence and Statistics}, volume~84 of \emph{Proceedings of
  Machine Learning Research}, 1608--1617 (PMLR).

\bibitem[{Goh \protect\BIBand{} Sim(2010)}]{Goh10}
Goh J, Sim M (2010) Distributionally robust optimization and its tractable
  approximations. \emph{Operations Research} 58(4-part-1):902--917.

\bibitem[{Goodfellow et~al.(2014)Goodfellow, Shlens, \protect\BIBand{}
  Szegedy}]{goodfellow2014explaining}
Goodfellow IJ, Shlens J, Szegedy C (2014) Explaining and harnessing adversarial
  examples. \emph{arXiv preprint arXiv:1412.6572} .

\bibitem[{H{\"a}rdle(1990)}]{hardle1990applied}
H{\"a}rdle W (1990) \emph{Applied nonparametric regression} (Cambridge
  university press).

\bibitem[{He et~al.(2016)He, Zhang, Ren, \protect\BIBand{} Sun}]{he2016deep}
He K, Zhang X, Ren S, Sun J (2016) Deep residual learning for image
  recognition. \emph{Proceedings of the IEEE Conference on Computer Vision and
  Pattern Recognition}, 770--778.

\bibitem[{He et~al.(2017)He, Wei, Chen, Carlini, \protect\BIBand{}
  Song}]{he2017adversarial}
He W, Wei J, Chen X, Carlini N, Song D (2017) Adversarial example defense:
  Ensembles of weak defenses are not strong. \emph{WOOT}, 15--15.

\bibitem[{Hu et~al.(2020{\natexlab{a}})Hu, Chen, \protect\BIBand{}
  He}]{hu2020sample}
Hu Y, Chen X, He N (2020{\natexlab{a}}) Sample complexity of sample average
  approximation for conditional stochastic optimization. \emph{SIAM Journal on
  Optimization} 30(3):2103--2133.

\bibitem[{Hu et~al.(2021)Hu, Chen, \protect\BIBand{} He}]{hu2021biasvar}
Hu Y, Chen X, He N (2021) On the bias-variance-cost tradeoff of stochastic
  optimization. \emph{Advances in Neural Information Processing Systems}.

\bibitem[{Hu et~al.(2020{\natexlab{b}})Hu, Zhang, Chen, \protect\BIBand{}
  He}]{Yifan20}
Hu Y, Zhang S, Chen X, He N (2020{\natexlab{b}}) Biased stochastic first-order
  methods for conditional stochastic optimization and applications in meta
  learning. \emph{Advances in Neural Information Processing Systems},
  volume~33, 2759--2770.

\bibitem[{Hu \protect\BIBand{} Hong(2012)}]{hu2013kullback}
Hu Z, Hong LJ (2012) Kullback-leibler divergence constrained distributionally
  robust optimization. \emph{Optimization Online preprint Optimization
  Online:2012/11/3677} .

\bibitem[{Huang et~al.(2021)Huang, Ma, \protect\BIBand{}
  Lai}]{Riemannianhuang21}
Huang M, Ma S, Lai L (2021) A riemannian block coordinate descent method for
  computing the projection robust wasserstein distance. \emph{Proceedings of
  the 38th International Conference on Machine Learning}, 4446--4455.

\bibitem[{Kallenberg(1997)}]{kallenberg1997foundations}
Kallenberg O (1997) \emph{Foundations of modern probability}, volume~2
  (Springer).

\bibitem[{Kleywegt et~al.(2002)Kleywegt, Shapiro, \protect\BIBand{} Homem-de
  Mello}]{kleywegt2002sample}
Kleywegt AJ, Shapiro A, Homem-de Mello T (2002) The sample average
  approximation method for stochastic discrete optimization. \emph{SIAM Journal
  on optimization} 12(2):479--502.

\bibitem[{Krizhevsky \protect\BIBand{} Hinton(2009)}]{krizhevsky2009learning}
Krizhevsky A, Hinton G (2009) Learning multiple layers of features from tiny
  images. Technical report, Citeseer.

\bibitem[{Kruithof(1937)}]{kruithof1937telefoonverkeersrekening}
Kruithof J (1937) Telefoonverkeersrekening. \emph{De Ingenieur} 52:15--25.

\bibitem[{Kuhn et~al.(2024)Kuhn, Shafieezadeh-Abadeh, \protect\BIBand{}
  Wiesemann}]{daneilsurvey}
Kuhn D, Shafieezadeh-Abadeh S, Wiesemann W (2024) Distributionally robust
  optimization. \emph{arXiv preprint arXiv:2411.02549} .

\bibitem[{LeCun et~al.(1998)LeCun, Bottou, Bengio, \protect\BIBand{}
  Haffner}]{lecun1998gradient}
LeCun Y, Bottou L, Bengio Y, Haffner P (1998) Gradient-based learning applied
  to document recognition. \emph{Proceedings of the IEEE} 86(11):2278--2324.

\bibitem[{Levy et~al.(2020)Levy, Carmon, Duchi, \protect\BIBand{}
  Sidford}]{levy2020large}
Levy D, Carmon Y, Duchi JC, Sidford A (2020) Large-scale methods for
  distributionally robust optimization. \emph{Advances in Neural Information
  Processing Systems} 33:8847--8860.

\bibitem[{Li et~al.(2022)Li, Lin, Blanchet, \protect\BIBand{}
  Nguyen}]{li2022tikhonov}
Li J, Lin S, Blanchet J, Nguyen VA (2022) Tikhonov regularization is optimal
  transport robust under martingale constraints. \emph{Advances in Neural
  Information Processing Systems} 35:17677--17689.

\bibitem[{Lin et~al.(2020)Lin, Fan, Ho, Cuturi, \protect\BIBand{}
  Jordan}]{lin2020projection2}
Lin T, Fan C, Ho N, Cuturi M, Jordan M (2020) Projection robust wasserstein
  distance and riemannian optimization. \emph{Advances in Neural Information
  Processing Systems}, volume~33, 9383--9397.

\bibitem[{Lin et~al.(2022)Lin, Ho, \protect\BIBand{}
  Jordan}]{lin2022efficiency}
Lin T, Ho N, Jordan MI (2022) On the efficiency of entropic regularized
  algorithms for optimal transport. \emph{Journal of Machine Learning Research}
  23(137):1--42.

\bibitem[{Liu et~al.(2021)Liu, Yuan, \protect\BIBand{} Zhang}]{liu2021discrete}
Liu Y, Yuan X, Zhang J (2021) Discrete approximation scheme in distributionally
  robust optimization. \emph{Numer Math Theory Methods Appl} 14(2):285--320.

\bibitem[{Luise et~al.(2018)Luise, Rudi, Pontil, \protect\BIBand{}
  Ciliberto}]{luise2018differential}
Luise G, Rudi A, Pontil M, Ciliberto C (2018) Differential properties of
  sinkhorn approximation for learning with wasserstein distance. \emph{Advances
  in Neural Information Processing Systems}.

\bibitem[{Luo \protect\BIBand{} Mehrotra(2019)}]{Luo19}
Luo F, Mehrotra S (2019) Decomposition algorithm for distributionally robust
  optimization using wasserstein metric with an application to a class of
  regression models. \emph{European Journal of Operational Research}
  278(1):20--35.

\bibitem[{Madry et~al.(2018)Madry, Makelov, Schmidt, Tsipras, \protect\BIBand{}
  Vladu}]{madry2017towards}
Madry A, Makelov A, Schmidt L, Tsipras D, Vladu A (2018) Towards deep learning
  models resistant to adversarial attacks. \emph{International Conference on
  Learning Representations}.

\bibitem[{Mensch \protect\BIBand{} Peyr{\'e}(2020)}]{mensch2020online}
Mensch A, Peyr{\'e} G (2020) Online sinkhorn: Optimal transport distances from
  sample streams. \emph{Advances in Neural Information Processing Systems}
  33:1657--1667.

\bibitem[{Mohajerin~Esfahani \protect\BIBand{} Kuhn(2017)}]{Mohajerin18}
Mohajerin~Esfahani P, Kuhn D (2017) Data-driven distributionally robust
  optimization using the wasserstein metric: performance guarantees and
  tractable reformulations. \emph{Mathematical Programming} 171(1):115--166.

\bibitem[{Namkoong \protect\BIBand{} Duchi(2016)}]{namkoong2016stochastic}
Namkoong H, Duchi JC (2016) Stochastic gradient methods for distributionally
  robust optimization with f-divergences. \emph{Advances in Neural Information
  Processing Systems}, volume~29, 2208--2216.

\bibitem[{Natarajan et~al.(2009)Natarajan, Song, \protect\BIBand{}
  Teo}]{natarajan2009persistency}
Natarajan K, Song M, Teo CP (2009) Persistency model and its applications in
  choice modeling. \emph{Management Science} 55(3):453--469.

\bibitem[{Nemirovski et~al.(2009)Nemirovski, Juditsky, Lan, \protect\BIBand{}
  Shapiro}]{nemirovski2009robust}
Nemirovski A, Juditsky A, Lan G, Shapiro A (2009) Robust stochastic
  approximation approach to stochastic programming. \emph{SIAM Journal on
  optimization} 19(4):1574--1609.

\bibitem[{Nemirovsky \protect\BIBand{} Yudin(1983)}]{Blair85}
Nemirovsky A, Yudin D (1983) Problem complexity and method efficiency in
  optimization. \emph{John Wiley \& Sons} .

\bibitem[{Nesterov \protect\BIBand{} Nemirovskii(1994)}]{nesterov1994interior}
Nesterov Y, Nemirovskii A (1994) \emph{Interior-point polynomial algorithms in
  convex programming} (SIAM).

\bibitem[{Nguyen et~al.(2020)Nguyen, Si, \protect\BIBand{}
  Blanchet}]{nguyen2020robust}
Nguyen VA, Si N, Blanchet J (2020) Robust bayesian classification using an
  optimistic score ratio. \emph{International Conference on Machine Learning},
  7327--7337.

\bibitem[{Nguyen et~al.(2024)Nguyen, Zhang, Wang, Blanchet, Delage,
  \protect\BIBand{} Ye}]{nguyen2021robustifying}
Nguyen VA, Zhang F, Wang S, Blanchet J, Delage E, Ye Y (2024) Robustifying
  conditional portfolio decisions via optimal transport. \emph{Operations
  Research} .

\bibitem[{Papernot et~al.(2017)Papernot, McDaniel, Goodfellow, Jha, Celik,
  \protect\BIBand{} Swami}]{papernot2017practical}
Papernot N, McDaniel P, Goodfellow I, Jha S, Celik ZB, Swami A (2017) Practical
  black-box attacks against machine learning. \emph{Proceedings of the 2017 ACM
  on Asia conference on computer and communications security}, 506--519.

\bibitem[{Papernot et~al.(2016{\natexlab{a}})Papernot, McDaniel, Jha,
  Fredrikson, Celik, \protect\BIBand{} Swami}]{papernot2016limitations}
Papernot N, McDaniel P, Jha S, Fredrikson M, Celik ZB, Swami A
  (2016{\natexlab{a}}) The limitations of deep learning in adversarial
  settings. \emph{2016 IEEE European symposium on security and privacy
  (EuroS\&P)}, 372--387 (IEEE).

\bibitem[{Papernot et~al.(2016{\natexlab{b}})Papernot, McDaniel, Wu, Jha,
  \protect\BIBand{} Swami}]{papernot2016distillation}
Papernot N, McDaniel P, Wu X, Jha S, Swami A (2016{\natexlab{b}}) Distillation
  as a defense to adversarial perturbations against deep neural networks.
  \emph{2016 IEEE symposium on security and privacy (SP)}, 582--597 (IEEE).

\bibitem[{Patrini et~al.(2020)Patrini, van~den Berg, Forre, Carioni, Bhargav,
  Welling, Genewein, \protect\BIBand{} Nielsen}]{patrini2020sinkhorn}
Patrini G, van~den Berg R, Forre P, Carioni M, Bhargav S, Welling M, Genewein
  T, Nielsen F (2020) Sinkhorn autoencoders. \emph{Uncertainty in Artificial
  Intelligence}, 733--743.

\bibitem[{Petzka et~al.(2018)Petzka, Fischer, \protect\BIBand{}
  Lukovnikov}]{petzka2018on}
Petzka H, Fischer A, Lukovnikov D (2018) On the regularization of wasserstein
  {GAN}s. \emph{International Conference on Learning Representations}.

\bibitem[{Peyre \protect\BIBand{} Cuturi(2019)}]{Computational19}
Peyre G, Cuturi M (2019) Computational optimal transport: With applications to
  data science. \emph{Foundations and Trends in Machine Learning}
  11(5-6):355--607.

\bibitem[{Pflug \protect\BIBand{} Wozabal(2007)}]{pflug2007ambiguity}
Pflug G, Wozabal D (2007) Ambiguity in portfolio selection. \emph{Quantitative
  Finance} 7(4):435--442.

\bibitem[{Pichler \protect\BIBand{} Shapiro(2021)}]{pichler2021mathematical}
Pichler A, Shapiro A (2021) Mathematical foundations of distributionally robust
  multistage optimization. \emph{SIAM Journal on Optimization}
  31(4):3044--3067.

\bibitem[{Popescu(2005)}]{popescu2005semidefinite}
Popescu I (2005) A semidefinite programming approach to optimal-moment bounds
  for convex classes of distributions. \emph{Mathematics of Operations
  Research} 30(3):632--657.

\bibitem[{Qi et~al.(2023)Qi, Lyu, Chan, Bai, \protect\BIBand{}
  Yang}]{qi2022stochastic}
Qi Q, Lyu J, Chan KS, Bai EW, Yang T (2023) Stochastic constrained {DRO} with a
  complexity independent of sample size. \emph{Transactions on Machine Learning
  Research} ISSN 2835-8856.

\bibitem[{Rockafellar et~al.(1999)Rockafellar, Uryasev
  et~al.}]{rockafellar1999optimization}
Rockafellar RT, Uryasev S, et~al. (1999) Optimization of conditional
  value-at-risk. \emph{Journal of risk} 2:21--42.

\bibitem[{Rozsa et~al.(2018)Rozsa, Gunther, \protect\BIBand{}
  Boult}]{rozsa2016towards}
Rozsa A, Gunther M, Boult TE (2018) Towards robust deep neural networks with
  bang. \emph{2018 IEEE Winter Conference on Applications of Computer Vision
  (WACV)}, 803--811.

\bibitem[{Scarf(1957)}]{scarf1958min}
Scarf H (1957) A min-max solution of an inventory problem. \emph{Studies in the
  mathematical theory of inventory and production} .

\bibitem[{Selvi et~al.(2022)Selvi, Belbasi, Haugh, \protect\BIBand{}
  Wiesemann}]{selvi2022wasserstein}
Selvi A, Belbasi MR, Haugh MB, Wiesemann W (2022) Wasserstein logistic
  regression with mixed features. \emph{Advances in Neural Information
  Processing Systems}.

\bibitem[{Shafieezadeh-Abadeh et~al.(2023)Shafieezadeh-Abadeh, Aolaritei,
  D{\"o}rfler, \protect\BIBand{} Kuhn}]{shafieezadeh2023new}
Shafieezadeh-Abadeh S, Aolaritei L, D{\"o}rfler F, Kuhn D (2023) New
  perspectives on regularization and computation in optimal transport-based
  distributionally robust optimization. \emph{arXiv preprint arXiv:2303.03900}
  .

\bibitem[{Shafieezadeh-Abadeh et~al.(2019)Shafieezadeh-Abadeh, Kuhn,
  \protect\BIBand{} Esfahani}]{shafieezadeh2019regularization}
Shafieezadeh-Abadeh S, Kuhn D, Esfahani PM (2019) Regularization via mass
  transportation. \emph{Journal of Machine Learning Research} 20(103):1--68.

\bibitem[{Shafieezadeh~Abadeh et~al.(2015)Shafieezadeh~Abadeh,
  Mohajerin~Esfahani, \protect\BIBand{} Kuhn}]{Shafieezadeh15}
Shafieezadeh~Abadeh S, Mohajerin~Esfahani PM, Kuhn D (2015) Distributionally
  robust logistic regression. \emph{Advances in Neural Information Processing
  Systems}, volume~28.

\bibitem[{Shapiro(2001)}]{shapiro2001duality}
Shapiro A (2001) On duality theory of conic linear problems.
  \emph{Semi-infinite programming}, 135--165 (Springer).

\bibitem[{Shapiro(2017)}]{shapiro2017distributionally}
Shapiro A (2017) Distributionally robust stochastic programming. \emph{SIAM
  Journal on Optimization} 27(4):2258--2275.

\bibitem[{Shapiro et~al.(2023)Shapiro, Zhou, \protect\BIBand{}
  Lin}]{shapiro2021bayesian}
Shapiro A, Zhou E, Lin Y (2023) Bayesian distributionally robust optimization.
  \emph{SIAM Journal on Optimization} 33(2):1279--1304.

\bibitem[{Singh \protect\BIBand{} Zhang(2021)}]{SINGH2021121}
Singh D, Zhang S (2021) Distributionally robust profit opportunities.
  \emph{Operations Research Letters} 49(1):121--128.

\bibitem[{Singh \protect\BIBand{} Zhang(2022)}]{singh2020tight}
Singh D, Zhang S (2022) Tight bounds for a class of data-driven
  distributionally robust risk measures. \emph{Applied Mathematics \&
  Optimization} 85(1):1--41.

\bibitem[{Sinha et~al.(2018)Sinha, Namkoong, \protect\BIBand{}
  Duchi}]{sinha2018certifiable}
Sinha A, Namkoong H, Duchi J (2018) Certifiable distributional robustness with
  principled adversarial training. \emph{International Conference on Learning
  Representations}.

\bibitem[{Sinkhorn(1964)}]{sinkhorn1964relationship}
Sinkhorn R (1964) A relationship between arbitrary positive matrices and doubly
  stochastic matrices. \emph{The annals of mathematical statistics}
  35(2):876--879.

\bibitem[{Song et~al.(2023)Song, He, Ding, \protect\BIBand{}
  Zhao}]{song2022efficient}
Song J, He N, Ding L, Zhao C (2023) Provably convergent policy optimization via
  metric-aware trust region methods. \emph{Transactions on Machine Learning
  Research} ISSN 2835-8856.

\bibitem[{Staib \protect\BIBand{} Jegelka(2019)}]{staib2019distributionally}
Staib M, Jegelka S (2019) Distributionally robust optimization and
  generalization in kernel methods. \emph{Advances in Neural Information
  Processing Systems} 32:9134--9144.

\bibitem[{TinyImageNet(2014)}]{tinyimagenet}
TinyImageNet (2014) {TinyImageNet Visual Recognition Challenge}.

\bibitem[{Tram{\`e}r et~al.(2018)Tram{\`e}r, Kurakin, Papernot, Goodfellow,
  Boneh, \protect\BIBand{} McDaniel}]{tramer2017ensemble}
Tram{\`e}r F, Kurakin A, Papernot N, Goodfellow I, Boneh D, McDaniel P (2018)
  Ensemble adversarial training: Attacks and defenses. \emph{International
  Conference on Learning Representations}.

\bibitem[{Van~Parys et~al.(2015)Van~Parys, Goulart, \protect\BIBand{}
  Kuhn}]{van2015generalized}
Van~Parys BP, Goulart PJ, Kuhn D (2015) Generalized gauss inequalities via
  semidefinite programming. \emph{Mathematical Programming} 156(1-2):271--302.

\bibitem[{Vandenberghe \protect\BIBand{}
  Boyd(1995)}]{vandenberghe1995semidefinite}
Vandenberghe L, Boyd S (1995) Semidefinite programming. \emph{SIAM review}
  38(1):49--95.

\bibitem[{Wang et~al.(2018)Wang, Gao, Qiu, Wang, \protect\BIBand{}
  Xin}]{Wang18}
Wang C, Gao R, Qiu F, Wang J, Xin L (2018) Risk-based distributionally robust
  optimal power flow with dynamic line rating. \emph{IEEE Transactions on Power
  Systems} 33(6):6074--6086.

\bibitem[{Wang et~al.(2022{\natexlab{a}})Wang, Gao, \protect\BIBand{}
  Xie}]{wang2020kerneltwosample}
Wang J, Gao R, Xie Y (2022{\natexlab{a}}) Two-sample test with kernel projected
  wasserstein distance. \emph{Proceedings of The 25th International Conference
  on Artificial Intelligence and Statistics}, volume 151, 8022--8055 (PMLR).

\bibitem[{Wang et~al.(2022{\natexlab{b}})Wang, Gao, \protect\BIBand{}
  Zha}]{wang2021reliable}
Wang J, Gao R, Zha H (2022{\natexlab{b}}) Reliable off-policy evaluation for
  reinforcement learning. \emph{Operations Research} .

\bibitem[{Wang et~al.(2015)Wang, Glynn, \protect\BIBand{}
  Ye}]{wang2016likelihood}
Wang Z, Glynn PW, Ye Y (2015) Likelihood robust optimization for data-driven
  problems. \emph{Computational Management Science} 13(2):241--261.

\bibitem[{Wiesemann et~al.(2014)Wiesemann, Kuhn, \protect\BIBand{}
  Sim}]{wiesemann2014distributionally}
Wiesemann W, Kuhn D, Sim M (2014) Distributionally robust convex optimization.
  \emph{Operations Research} 62(6):1358--1376.

\bibitem[{Wozabal(2012)}]{wozabal2012framework}
Wozabal D (2012) A framework for optimization under ambiguity. \emph{Annals of
  Operations Research} 193(1):21--47.

\bibitem[{Xie(2019)}]{xie2019distributionally}
Xie W (2019) On distributionally robust chance constrained programs with
  wasserstein distance. \emph{Mathematical Programming} 186(1):115--155.

\bibitem[{Yang(2017)}]{yang2017convex}
Yang I (2017) A convex optimization approach to distributionally robust markov
  decision processes with wasserstein distance. \emph{IEEE control systems
  letters} 1(1):164--169.

\bibitem[{Yang(2020)}]{YangWDRO20}
Yang I (2020) Wasserstein distributionally robust stochastic control: A
  data-driven approach. \emph{IEEE Transactions on Automatic Control}
  66(8):3863--3870.

\bibitem[{Yu et~al.(2022)Yu, Lin, Mazumdar, \protect\BIBand{}
  Jordan}]{yu2022fast}
Yu Y, Lin T, Mazumdar EV, Jordan M (2022) Fast distributionally robust learning
  with variance-reduced min-max optimization. \emph{International Conference on
  Artificial Intelligence and Statistics}, 1219--1250.

\bibitem[{Yule(1912)}]{yule1912methods}
Yule GU (1912) On the methods of measuring association between two attributes.
  \emph{Journal of the Royal Statistical Society} 75(6):579--652.

\bibitem[{Zhao \protect\BIBand{} Guan(2018)}]{zhao2018data}
Zhao C, Guan Y (2018) Data-driven risk-averse stochastic optimization with
  wasserstein metric. \emph{Operations Research Letters} 46(2):262--267.

\bibitem[{Zhu et~al.(2021)Zhu, Jitkrittum, Diehl, \protect\BIBand{}
  Sch{\"o}lkopf}]{Kernelzhu}
Zhu J, Jitkrittum W, Diehl M, Sch{\"o}lkopf B (2021) Kernel distributionally
  robust optimization: Generalized duality theorem and stochastic
  approximation. \emph{Proceedings of The 24th International Conference on
  Artificial Intelligence and Statistics}, 280--288.

\bibitem[{Zymler et~al.(2013)Zymler, Kuhn, \protect\BIBand{} Rustem}]{Zymler13}
Zymler S, Kuhn D, Rustem B (2013) Distributionally robust joint chance
  constraints with second-order moment information. \emph{Mathematical
  Programming} 137(1):167--198.

\end{thebibliography}

\ECSwitch

\ECHead{Supplementary for \emph{``Sinkhorn Distributionally Robust Optimization''}}

\section{Detailed Experiment Setup}\label{Sec:experiment:setup}
Unless stated otherwise, we solved the SAA, Wasserstein DRO, and KL-divergence DRO baseline models exactly using the off-the-shelf solver Mosek~\citep{mosek}. Optimization hyperparameters, such as step size, maximum iterations, and number of levels, were tuned to minimize training error after 10 outer iterations. 
We use RT-MLMC subgradient estimator to solve the Sinkhorn DRO model.
We employed the \emph{warm starting} strategy during the iterative procedure: 
we set the initial guess of parameter $\theta$ at the beginning of outer iteration as the one obtained from the SAA approach.
At other outer iterations, the initial guess of parameter $\theta$ is set to be the final obtained solution $\theta$ at the last outer iteration.
The following subsections outline some special reformulations, optimization algorithms used to solve the baseline models. %

\subsection{Setup for Newsvendor Problem and Running Time}
To solve the $2$-Wasserstein DRO model with radius $\rho$, we approximate the support of worst-case distribution using discrete grid points.
Denote by $\mathcal{D}_n=\{x_1,\ldots,x_n\}$ the set of observed $n$ samples and $\mathcal{G}_{200-n}$ the set of $200-n$ points evenly supported on the interval $[0,10]$.
Then the support of worst-case distribution is restricted to $\mathcal{D}_n\cup\mathcal{G}_{200-n}:=\{\hat{z}_1,\ldots,\hat{z}_{200}\}$.
The corresponding $2$-Wasserstein DRO problem has the following linear programming reformulation:
\[
\begin{aligned}
\min_{\theta, \lambda, s}&\quad \lambda\rho + \frac{1}{n}\sum_{i=1}^ns_i\\
\mbox{s.t.}&\quad k\theta - u\min(\theta, \hat{z}_j) - \lambda(x_i - \hat{z}_j)^2\le s_i,\quad \forall i\in[n], \forall j\in[200].
\end{aligned}
\]

\subsection{Setup for Mean-risk Portfolio Optimization}

From \citep[Eq.~(27)]{Mohajerin18} we can see that the $1$-Wasserstein DRO formulation with radius $\rho$ for the portfolio optimization problem becomes 
\[
\begin{aligned}
\min_{\theta,\tau, \lambda, s}&\quad
\lambda\rho + \frac{1}{n}\sum_{i=1}^ns_i\\
\mbox{s.t.}&\quad \theta\in\Theta,\quad b_j\tau + a_j\inp{\theta}{\hat{z}_i}\le s_i, i\in[n], j\in[H],\\
           &\quad \|a_j\theta\|_2\le \lambda, j\in[H].
\end{aligned}
\]
Also, we argue that the $2$-Wasserstein DRO formulation with radius $\rho$ for the portfolio optimization problem has a finite convex reformulation:
\begin{align*}
&\inf_{\theta\in\Theta, \tau}~\sup_{\bP:~W_2(\bP, \hP_n)\le \rho}\mathbb{E}_{\bP}\big[ 
\max_{j\in[H]}a_j\inp{\theta}{z} + b_j\tau
\big]\\
=&\inf_{\theta\in\Theta, \tau, \lambda\ge0}~
\left\{ 
\lambda\rho^2 + \frac{1}{n}\sum_{i=1}^n\sup_{s_i}~\left\{
\max_{j\in[H]}a_j\inp{\theta}{s_i} + b_j\tau - \lambda\|s_i - \hat{z}_i\|_2^2
\right\}
\right\}.
\end{align*}
In particular, the inner subproblem has the following reformulation:
\begin{align*}
&\sup_{s_i}~\left\{
\max_{j\in[H]}a_j\inp{\theta}{s_i} + b_j\tau - \lambda\|s_i - \hat{z}_i\|_2^2
\right\}\\
=&\max_{j\in[H]}~b_j\tau + \sup_{s_i}\left\{
a_j\inp{\theta}{s_i} - \lambda\|s_i - \hat{z}_i\|_2^2
\right\}\\
=&\max_{j\in[H]}~b_j\tau + \frac{a_j^2}{4\lambda}\|\theta\|^2_2 + a_j\inp{\theta}{\hat{z}_i}.
\end{align*}
Hence, the $2$-Wasserstein DRO can be reformulated as
\[
\begin{aligned}
\min_{\theta,\tau, \lambda, s}&\quad
\lambda\rho^2 + \frac{1}{n}\sum_{i=1}^ns_i\\
\mbox{s.t.}&\quad \theta\in\Theta,\quad 
           b_j\tau +a_j\inp{\theta}{\hat{z}_i} +  \frac{a_j^2}{4\lambda}\|\theta\|_2^2\le s_i
           ,\quad i\in[n], j\in[H].
\end{aligned}
\]

\subsection{Setup for Adversarial Multi-class Logistic Regression}
\label{Appendix:adv:re}

The procedure for generating various adversarial perturbations is reported in the following:
\begin{enumerate}
    \item 
    For a given classifer $B$ and data sample $(x,{\bm y})$, the $\ell_p$-norm ($p\in\{1,2\}$) adversarial attack based on projected gradient method~\citep{madry2017towards} iterates as follows: $x_0\leftarrow x$ and 
\[
\left\{ 
\begin{aligned}
\Delta x^{k+1}&\leftarrow \argmax_{\|\eta\|_p\le \xi}~\Big\{ 
\nabla_x h_B(x^k, {\bm y})\trans \eta
\Big\},\\ 
x^{k+1}&\leftarrow \text{Proj}_{\{x':~\|x-x'\|_p\le \xi\}}~\Big\{ 
x^k + \frac{\alpha}{\sqrt{k+1}} \Delta x^{k+1}
\Big\}.
\end{aligned}
\right.
\]
We perform the gradient update above for $15$ steps with initial learning rate $\alpha=1$.
\Jie{When $p=1$, the radius of attack $\xi\in\{\texttt{0}, \texttt{3e-3}, \texttt{6e-3}, \texttt{9e-3}, \texttt{1.2e-2}\}\cdot \varrho$; and when $p=2$, the radius $\xi\in\{\texttt{0}, \texttt{8e-3}, \texttt{1.6e-2}, \texttt{2.4e-2}, \texttt{3.2e-2}\}\cdot \varrho$.
}
\item
For a given feature vector $x$, the perturbed feature using white Laplacian noise becomes $x + \xi\cdot \zeta$, where the random vector $\zeta$ follows the isotropic Laplace distribution with zero mean and unit variance.
\Jie{The ratio $\xi\in\{\texttt{0}, \texttt{2e-3}, \texttt{4e-3}, \texttt{6e-3}, \texttt{8e-3}\}\cdot \varrho.$}
Similarly, the perturbed feature using white Gaussian noise becomes $x + \xi\cdot \zeta$, with $\zeta$ being the isotropic Gaussian distribution with  zero mean and unit variance.
\Jie{In this case, the ratio $\xi\in\{\texttt{0}, \texttt{5e-2}, \texttt{1e-1}, \texttt{1.5e-1}, \texttt{2e-1}\}\cdot \varrho.$
}
\end{enumerate}
In this example, we use stochastic gradient methods to solve the SAA formulation and all penalized DRO formulations.
We terminate the training of SAA or DRO models when the number of epoches, i.e., the number of times for processes each training sample, exceeds $30$.
In is worth mentioning that the Wasserstein DRO model with a fixed Lagrangian multiplier $\lambda$ using samples $\{x_i, {\bm y}_i\}_{i=1}^n$ can be reformulated as 
\begin{equation}
\min_{B}~
\frac{1}{n}\sum_{i=1}^n~\left[ 
\max_{x\in \mathbb{R}^d}~
\Big\{h_B(x, {\bm y}_i)  - \lambda c(x_i, x)\Big\}
\right].
\label{Eq:WDRO:minimax}
\end{equation}

\section{Additional Validation Experiments}
\label{Appendix:add:exp}
\subsection{Comparison of Optimization Algorithms: Linear Regression}
\label{Appendix:compare:opt:alg}
To examine the performance of different (sub)gradient estimators, we study the problem of distributionally robust linear regression (see the setup in Example~\ref{Example:DRO:LR}).
We take the nominal distribution $\hP$ as the empirical one based on samples $\{(a_i,b_i)\}_{i=1}^n$.
As a consequence, the inner objective function in \eqref{Eq:F:lambda} has the closed form expression: %
\[
F(\theta) = \frac{1}{n}\sum_{i=1}^n(a_i\trans\theta-b_i)^2 + \frac{\frac{1}{n}\sum_{i=1}^n(a_i\trans\theta-b_i)^2}{\frac{1}{2}\lambda\|\theta\|_2^{-2}-1} - \frac{\lambda\Reg}{2}\log\det\left( 
I - \frac{\theta\theta\trans}{\frac{1}{2}\lambda}
\right),\quad \text{if }\|\theta\|_2^2<\frac{\lambda}{2},
\]
and otherwise $F(\theta)=\infty$.
We take the constraint set $\Theta=\{\theta:~\|\theta\|_2^2\le0.999\cdot\frac{\lambda}{2}\}$.
Similar to the setup in \citep[Section~5.1]{li2022tikhonov}, we examine the performance using three LIBSVM regression real world datasets~\citep{libsvmregression}: housing, mg, and mpg.

\begin{figure}[!ht]
\centering
    \includegraphics[height=0.25\textwidth]{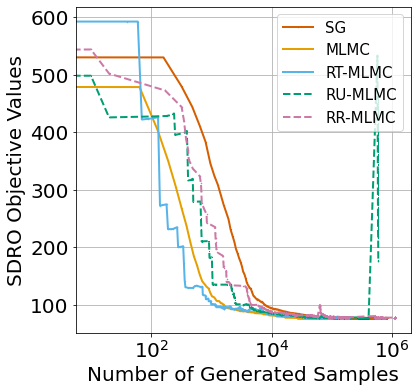}
    \includegraphics[height=0.25\textwidth]{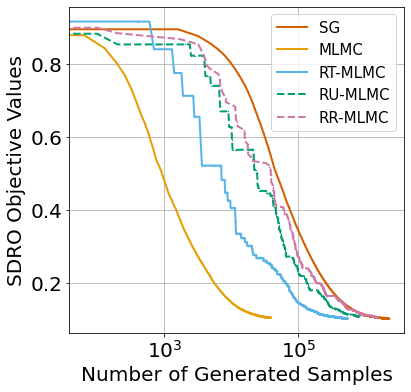}
    \includegraphics[height=0.25\textwidth]{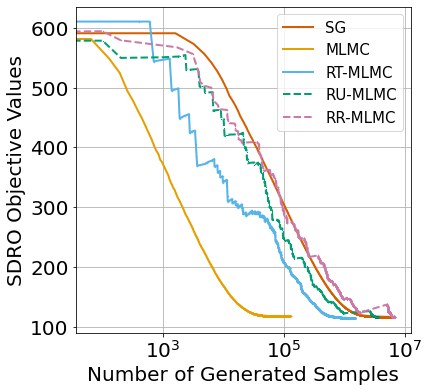}
    \includegraphics[height=0.25\textwidth]{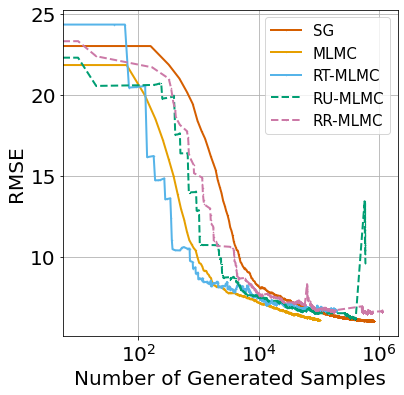}
    \includegraphics[height=0.25\textwidth]{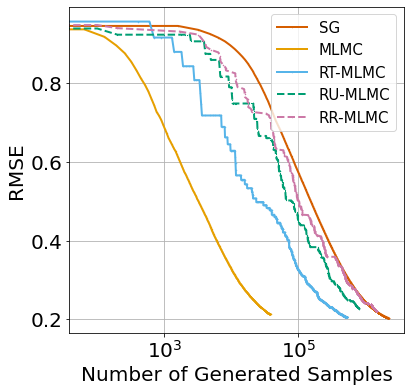}
    \includegraphics[height=0.25\textwidth]{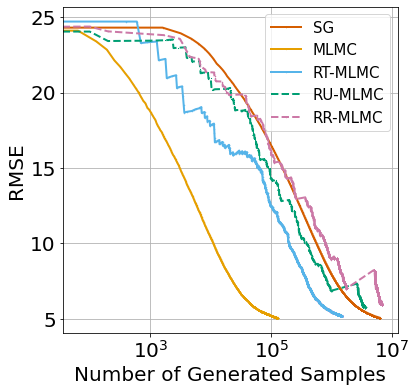}
\caption{Comparison results of SG, (V-)MLMC, RT-MLMC, RU-MLMC, and RR-MLMC on robust linear regression problem in terms of sample complexties from $\hP$ and $\bQ_{x,\Reg}$.
From left to right, the figures correspond to three different regression datasets: (a) housing; (b) mg; and (c) mpg.
From top to bottom, the figures correspond to plots of (a) Sinkhorn DRO objective values; and (b) RMSE of obtained solutions.}
\label{Fig:regression}
\end{figure}

The quality of proposed gradient estimators is examined in a single BSMD step with specified hyper-parameters $(\lambda,\Reg)=(10^3, 10^{-1})$.
For baseline comparison, we examine the SG, RT-MLMC estimators together with the (V-)MLMC, RU-MLMC, and RR-MLMC estimators that have been proposed in \citep{hu2021biasvar}.
We have validated in Theorem~\ref{Theorem:complexity:BSMD} that both SG and RT-MLMC estimators have convergence guarantees for smooth and nonsmooth loss functions, whereas SG estimator has slower convergence rate.
The (V-)MLMC estimator only have convergence guarantees for smooth loss functions, and RU-MLMC/RR-MLMC estimators do not have convergence guarantees as their (sub)gradient second-order moments are unbounded.

For a given solution $\theta$, we quantify its performance using the corresponding Sinkhorn DRO objective value.
Besides, we report its root-mean-square error~(RMSE) on training data. 
Thus, the smaller those two performance criteria are, the smaller the solution's optimization performance has.
Fig.~\ref{Fig:regression} shows the performance of various gradient estimators in terms of the number of generated samples from $\hP$ and $\bQ_{x,\Reg}, x\in\mathrm{supp}\hP$ based on these criteria. 
The results demonstrate that the SG scheme does not perform competitively, as expected from our theoretical analysis, which shows that SG has the worst complexity order. 
In contrast, using other four types of MLMC methods lead to faster convergence behavior.
While the RU-MLMC and RR-MLMC schemes exhibit competitive performance, the optimization procedure shows some oscillations. One possible explanation is that the variance values of those gradient estimators are unbounded, making these two approaches unstable.

\subsection{Comparison of Optimization Algorithms: Portfolio Optimization}
\label{Appendix:compare:opt:alg:port}
\begin{figure}[!ht]
\centering
    \includegraphics[height=0.25\textwidth]{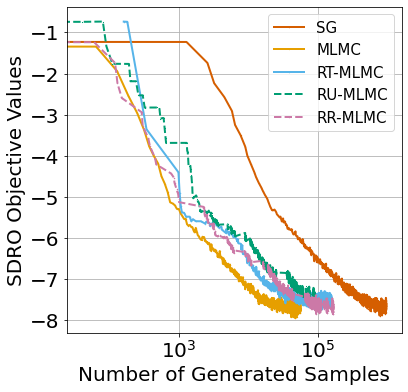}
    \includegraphics[height=0.25\textwidth]{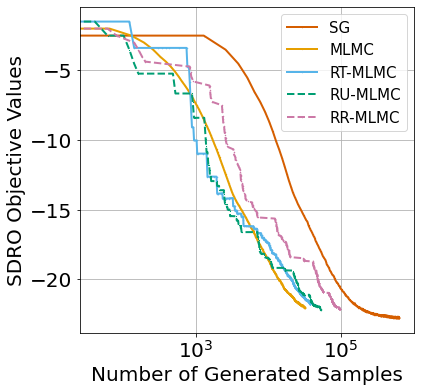}
    \includegraphics[height=0.25\textwidth]{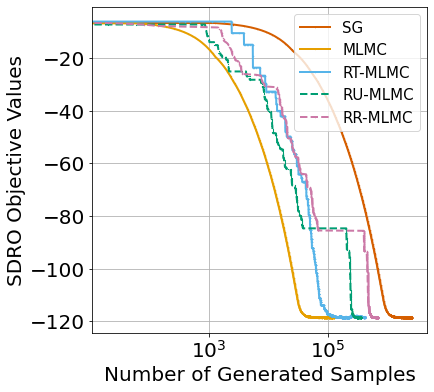}
\caption{Comparison results of SG, (V-)MLMC, RT-MLMC, RU-MLMC, and RR-MLMC on portfolio optimization problem.
From left to right, plots correspond to three different instances of $(n,d)\in\{(50, 50), (100, 100), (400, 400)\}$.
}
\label{Fig:portfolio}
\end{figure}

In this subsection, we validate the competitive performance of RT-MLMC gradient estimator on the case where the loss is convex and nonsmooth, and we try to solve the $2$-SDRO formulation.
We consider the portfolio optimization problem, and specify instances $(n,d) = (50, 50), (100, 100), (400, 400)$.
We quantify the performance of obtained solution using the Sinkhorn DRO objective value. Since in this problem setup no analytical expression of the objective value is available, we estimate the objective value using \eqref{Eq:RTMLMC:est} with hyper-parameters $L=8$ and $\nout_L=10^3$.
Fig.~\ref{Fig:portfolio} shows the performance in terms of the number of generated samples based on this criterion.
The results demonstrate that even for nonsmooth loss function, those listed MLMC-based gradient estimators have better performance than the SG estimator.
Besides, the proposed RT-MLMC and standard (V-)MLMC schemes have comparable performance, and in some cases (V-)MLMC estimator even has better performance.
It is an open question that whether the (V-)MLMC estimator will have the similar performance guarantees as the RT-MLMC estimator for convex nonsmooth optimization, which can be a topic for future study.

\subsection{\Jie{Comparison of Running Time for Different Baselines}}\label{Sec:baseline:run}

\Jie{
The computational time for the newsvendor problem in Section~\ref{Sec:newsvendor} is reported in Table~\ref{Table:time:newsvendor}.
We observe that the training time of $2$-Wasserstein DRO model increases quickly as the sample size increases, while the training time of other DRO models increases mildly in the training sample size.
}
\begin{table}[!ht]
  \centering
    \caption{Average computational time (in seconds) per problem instance for the newsvendor problem.}
\label{Table:time:newsvendor}
\Jie{
\begin{tabular}{l|ccccccccc}
\toprule
\multirow{2}{*}{Model} & \multicolumn{3}{c}{Exponential} & \multicolumn{3}{c}{Gamma} & \multicolumn{3}{c}{Gaussian Mixture} \\[6pt]
\cline{2-10}
                       & $n=10$     & $n=30$     & $n=100$     & $n=10$   & $n=30$   & $n=100$   & $n=10$       & $n=30$       & $n=100$      \\[3pt]
\midrule
\rowcolor{gray!0}
SAA                    & 4.11e-3         &  4.66e-3        &    4.67e-3       &   3.96e-3     &  4.57e-3      &   5.81e-3      &       3.82e-3     &     4.60e-3       &    4.79e-3      \\
KL-DRO                    &  6.92e-3         &  8.17e-3       &   1.15e-2      &  8.07e-3     &    8.24e-3  &         1.16e-2 &     7.77e-3       &   8.47e-3        &   1.12e-2      \\
\rowcolor{gray!0}
1-SDRO                 &   8.77e-2      &   8.88e-2     &    1.03e-1       &  2.76e-2     &   3.40e-2    &    4.72e-2    &   2.90e-2       &    3.13e-2      &     4.50e-2      \\
2-WDRO                 & 1.68e00       &   5.67e00       &    2.71e01       &  1.72e00     &  5.63e00    &    2.77e01    &     1.51e00       & 5.47e00          &   2.84e01       \\
\rowcolor{gray!0}
2-SDRO                 &   3.16e-2       &   3.77e-2     &   5.92e-2       &  2.64e-2      &   2.95e-2    &  5.02e-2       &  2.57e-2         &   3.10e-2         &  4.87e-2     \\
\bottomrule
\end{tabular}}
\end{table}

\Jie{
The computational time for the portfolio optimization problem in Section~\ref{Sec:portfolio} is reported in Table~\ref{Table:time:portfolio}.
We observe that the computational time of $1$- or $2$-SDRO model increases mildly as the problem input size increases.
Also, SDRO models do not have the smallest computational time in general. 
The reason is that in this example, other DRO models have tractable finite-dimensional conic programming formulations so that off-the-shelf software can solve them efficiently.
In contrast, Sinkhorn DRO models do not have special reformulation, but they can still be solved in a reasonable amount of time.
\begin{table}[!ht]
  \centering
  \Jie{
     \caption{
     Average computational time (in seconds) per problem instance for portfolio optimization problem.
}
\label{Table:time:portfolio}
{\footnotesize
\begin{tabular}{l>{\centering}p{1.5cm}>{\centering}p{1.5cm}>{\centering}p{1.5cm}>{\centering}p{1.5cm}>{\centering}p{1.5cm}>{\centering\arraybackslash}p{1.5cm}}
\toprule
$(n, d)$ Values   & SAA & KL-DRO & 1-WDRO & 1-SDRO & 2-WDRO & 2-SDRO \\ \midrule
\rowcolor{gray!0}
$(30, 30)$     &     6.76e-03 & 1.42e-02 & 7.80e-03 & 4.91e-02 & 8.95e-03 & 5.00e-02 \\
$(50, 30)$      &    7.31e-03 & 1.84e-02 & 8.33e-03 & 1.87e-01 & 1.11e-02 & 5.88e-02 \\
$(100, 30)$      &   8.99e-03 & 2.95e-02 & 1.03e-02 & 2.78e-01 & 1.12e-02 & 6.00e-02 \\
$(150, 30)$     &    1.12e-02 & 4.14e-02 & 1.21e-02 & 2.80e-01 & 1.22e-02 & 6.95e-02 \\
$(200, 30)$   &      1.12e-02 & 5.66e-02 & 1.35e-02 & 2.99e-01 & 1.48e-02 & 7.67e-02 \\
$(400, 30)$ &        1.89e-02 & 6.45e-02 & 2.09e-02 & 2.99e-01 & 2.30e-02 & 1.62e-01 \\
$(100, 5)$ &         5.76e-03 & 1.46e-02 & 6.79e-03 & 1.05e-01 & 7.62e-03 & 5.40e-02 \\
$(100, 10)$ &        6.18e-03 & 1.70e-02 & 7.70e-03 & 1.08e-01 & 8.73e-03 & 5.55e-02 \\
$(100, 20)$ &        7.43e-03 & 1.82e-02 & 8.41e-03 & 1.12e-01 & 9.44e-03 & 5.58e-02 \\
$(100, 40)$ &        9.87e-03 & 3.25e-02 & 1.13e-02 & 1.16e-01 & 1.18e-02 & 5.70e-02 \\
$(100, 80)$ &        1.31e-02 & 6.48e-02 & 1.56e-02 & 1.19e-01 & 1.68e-02 & 5.72e-02 \\
$(100, 100)$ &       1.54e-02 & 7.00e-02 & 1.87e-02 & 1.22e-01 & 1.93e-02 & 5.73e-02 \\\bottomrule
\end{tabular}}}
\end{table}
}

\Jie{
The computational time of adversarial multi-class classification problem in Section~\ref{Sec:classification} is reported in Table~\ref{Table:time:classification}, with the basic statistics of classification datasets presented in Table~\ref{Table:sta:classification}.
The results indicate that Sinkhorn DRO models have shorter computational time than Wasserstein DRO models in general.
Note that we solve all baseline methods with stochastic algorithms. 
For large-scale datasets optimizing the log-sum-exp type loss for Sinkhorn DRO seems to be more efficient than solving the minimax game formulation for Wasserstein DRO.
\begin{table}[H]
  \centering
\Jie{
     \caption{
Basic statistics of adversarial multi-class logistic regression datasets.
}
\label{Table:sta:classification}{\footnotesize
\begin{tabular}{c>{\centering}p{3cm}>{\centering}p{2cm}>{\centering}p{2cm}>{\centering\arraybackslash}p{2cm}}
\toprule
   & MNIST & CIFAR-10 & tinyImageNet & STL-10\\ \midrule
\rowcolor{gray!00}
$
\begin{array}{c}
\mbox{Image Size}\\
\mbox{(before pre-processing)}
\end{array}    
$
&  784   &  3072  &  12288   &  27648   \\ \hline
$
\begin{array}{c}
\mbox{Feature Dimension}\\
\mbox{(after pre-processing)}
\end{array}    
$      &  512   & 512   &  512    &  512  \\ \hline
\rowcolor{gray!00}
$\#$ of classes      & 10 & 10 & 200 & 10           \\ \hline
Training Size     & 50000 & 50000 & 90000 & 5000     \\ \hline
\rowcolor{gray!00}
Testing Size     & 10000 & 10000 & 10000 & 8000     \\ \bottomrule
\end{tabular}}}
\end{table}
\Jie{
\begin{table}[H]
  \centering\Jie{
     \caption{
Average computational time (in seconds) per problem instance for adversarial multi-class logistic regression problem.
}
\label{Table:time:classification}{\footnotesize
\begin{tabular}{l>{\centering}p{1.5cm}>{\centering}p{1.5cm}>{\centering}p{1.5cm}>{\centering}p{1.5cm}>{\centering}p{1.5cm}>{\centering\arraybackslash}p{1.5cm}}
\toprule
Dataset   & SAA & KL-DRO & 1-WDRO & 1-SDRO & 2-WDRO & 2-SDRO \\ \midrule
MNIST     &   37.2  &  60.1   &   154    &   94.1     &  166      &   84.0    \\ \hline
CIFAR-10      &  31.6     &  51.7     &  133    &  98.3      &  140  &   80.6   \\ \hline
tinyImageNet      &   58.1  &  102       &  248  &    153   &  259     &  143    \\ \hline
STL-10     &  3.42	  &    5.15     &     13.5    &    10.1    &  14.2   & 8.61       \\ 
\bottomrule
\end{tabular}}}
\end{table}
}}

\Jie{
\subsection{Coefficient of Prescriptiveness for Different Parameter(s) Combination}
\label{Appendix:coeff:p:combination}
In this subsection, we report the coefficient of prescriptiveness for different parameter(s) combination on instances which are omitted in the main content.
Specifically, 
\begin{itemize}
    \item 
Fig.~\ref{fig:newsvendor:lognormal:merged} and \ref{fig:newsvendor:Gaussian:merged} correspond to the omitted experiment results in Section~\ref{Sec:newsvendor}.
    \item
Fig.~\ref{Fig:portfolio:hp:1} and \ref{Fig:portfolio:hp:2} correspond to the omitted experiment results in Section~\ref{Sec:portfolio}.
    \item
Fig.~\ref{fig:adversarial:merged} corresponds to the omitted experiment results in Section~\ref{Sec:classification}.
\end{itemize}
\begin{figure}[!ht]
    \centering
    \includegraphics[width=\linewidth]{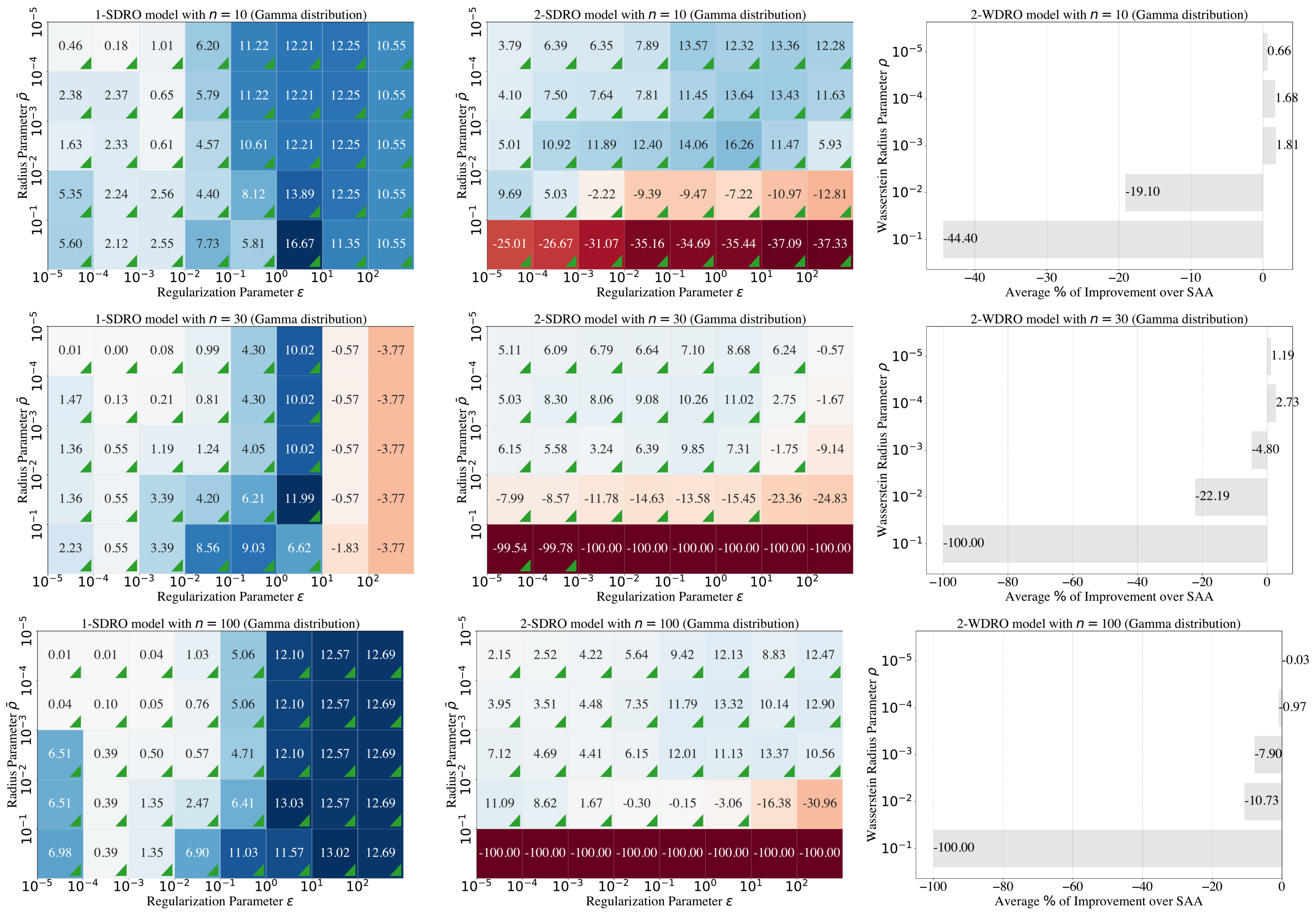}
    \caption{\Jie{Experiment results of the newsvendor model for gamma data distribution.
    Details of these subplots follow the same setup from Figure~\ref{fig:newsvendor:exponential:merged}.
    }}
    \label{fig:newsvendor:lognormal:merged}
\end{figure}
\begin{figure}[!ht]
    \centering
    \includegraphics[width=\linewidth]{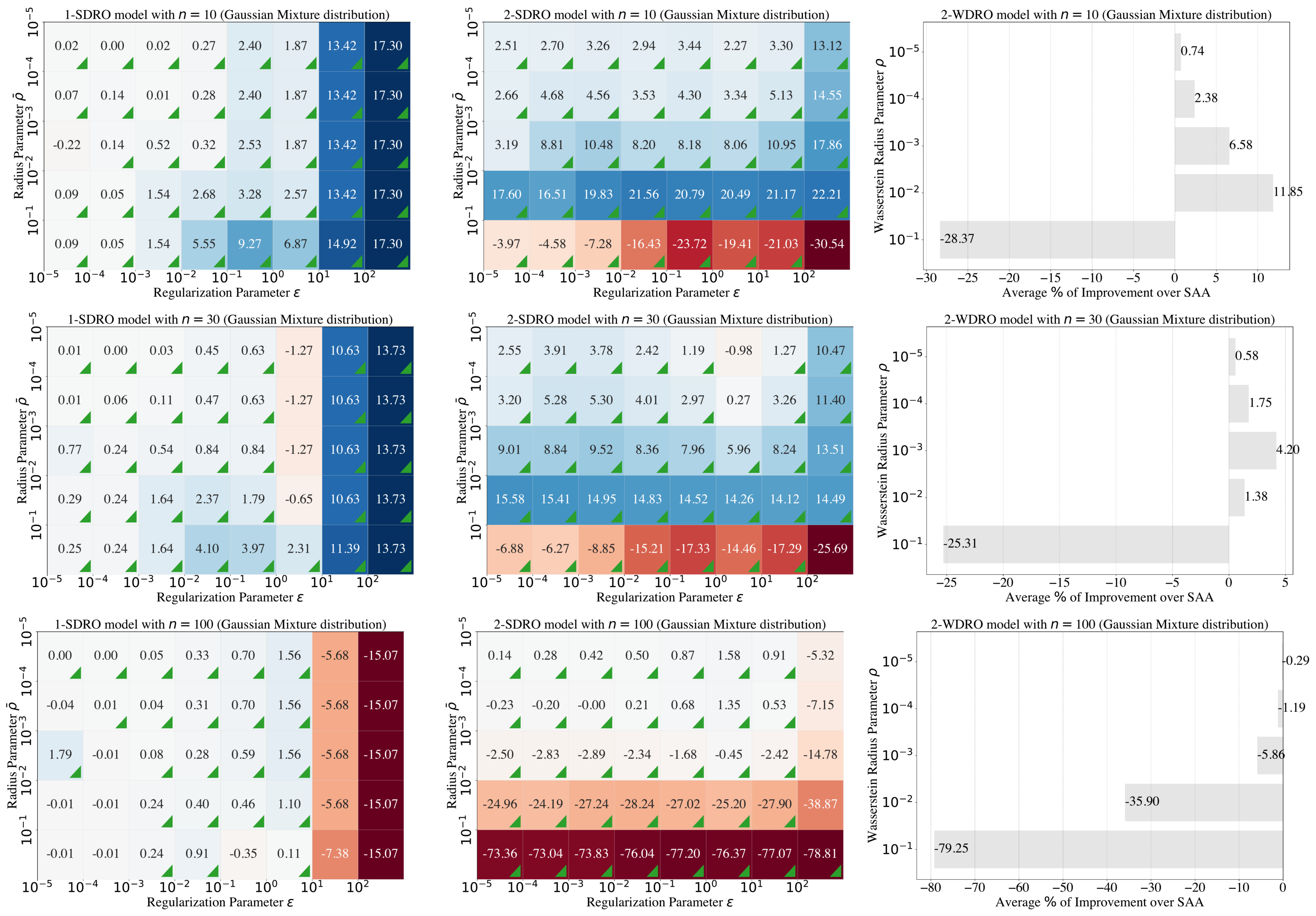}
    \caption{\Jie{Experiment results of the newsvendor model for the mixture of truncated normal distributions.
    Details of these subplots follow the same setup from Figure~\ref{fig:newsvendor:exponential:merged}.
    }}
    \label{fig:newsvendor:Gaussian:merged}
\end{figure}
\begin{figure}[!ht]
    \centering\Jie{
     \subfigure[$(n,d)=(50, 30)$]{\centering
     \includegraphics[width=0.45\textwidth]{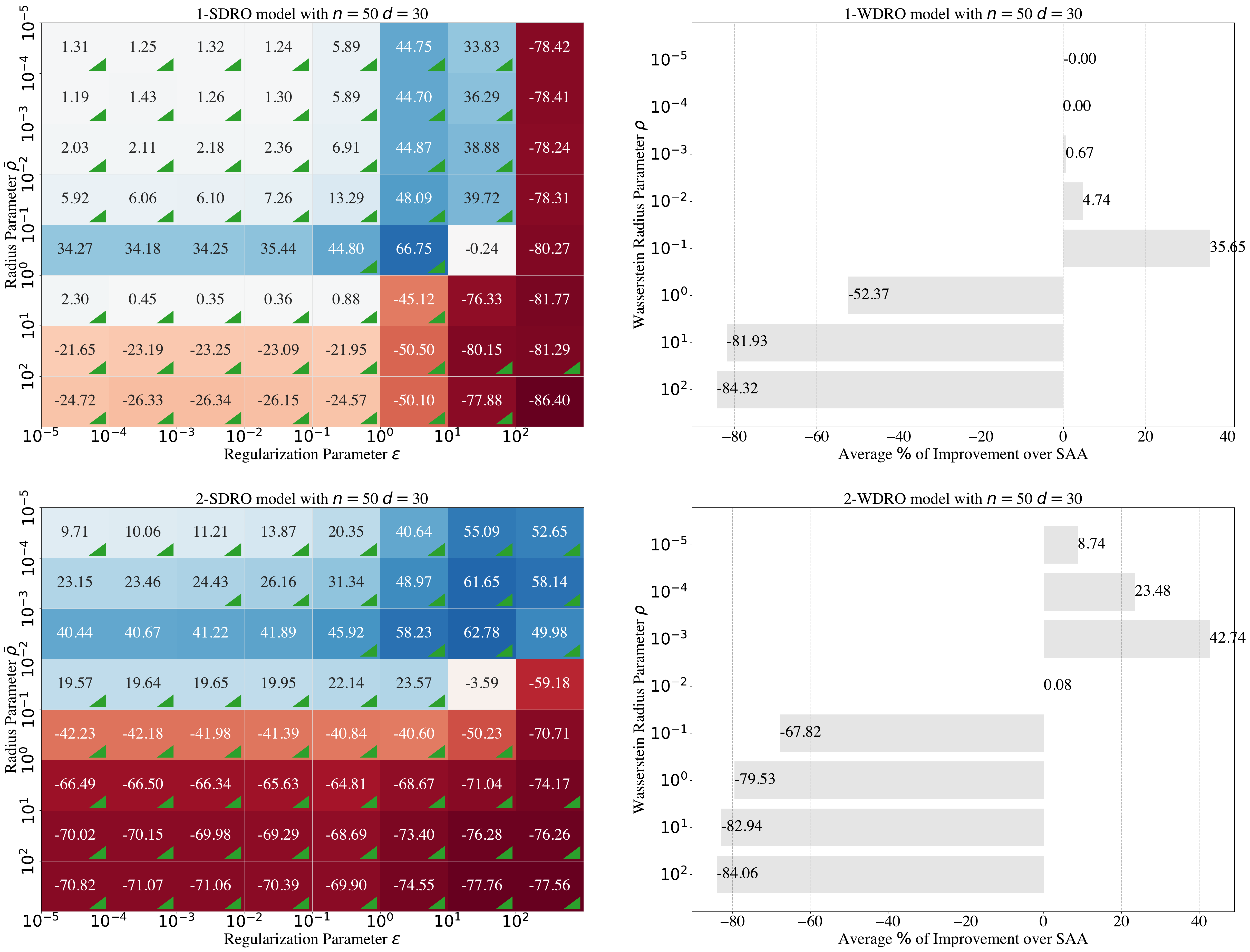}
    }
    \hfill
     \subfigure[$(n,d)=(100,30)$]{\centering
     \includegraphics[width=0.45\textwidth]{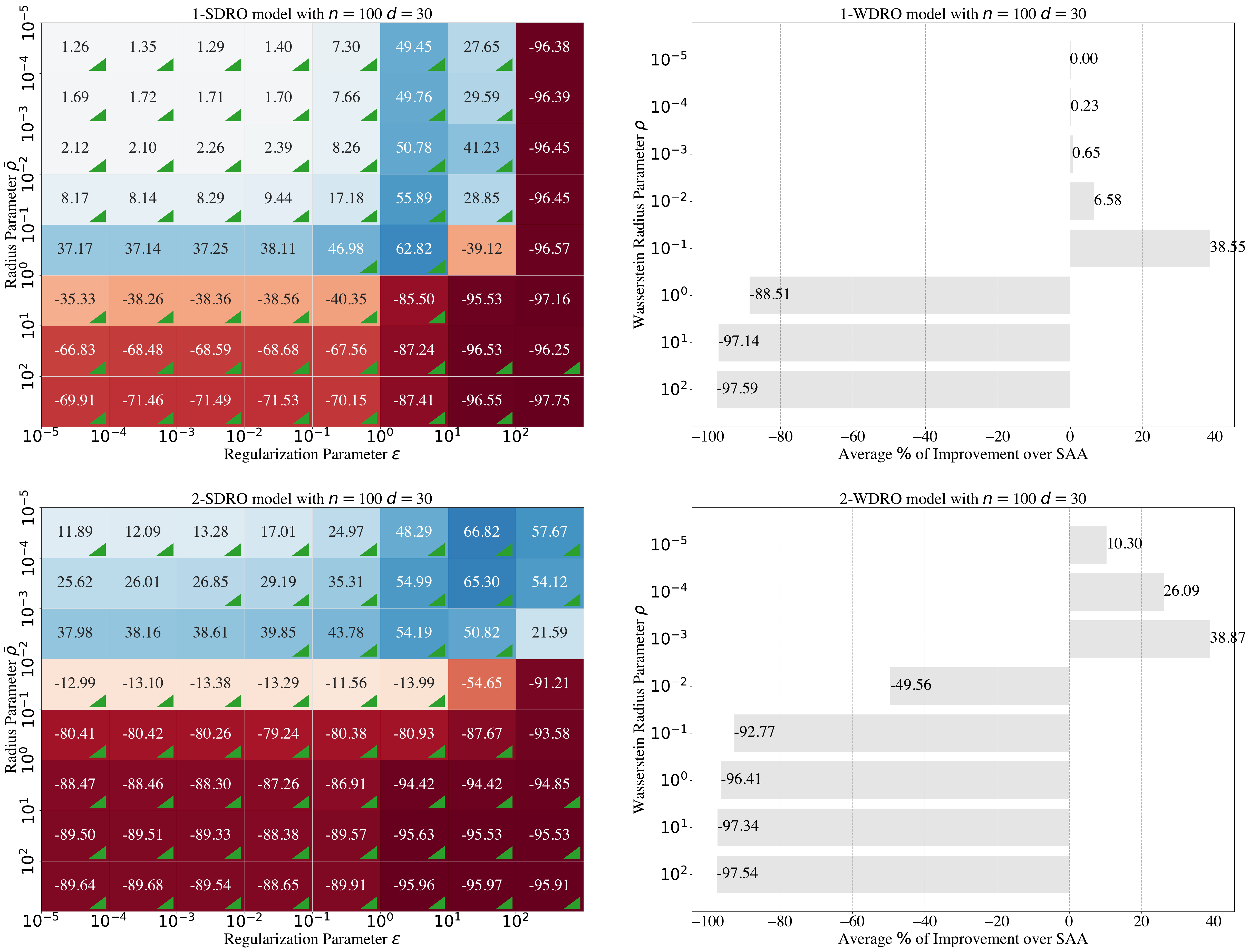}
    }\\
     \subfigure[$(n,d)=(150, 30)$]{\centering
     \includegraphics[width=0.45\textwidth]{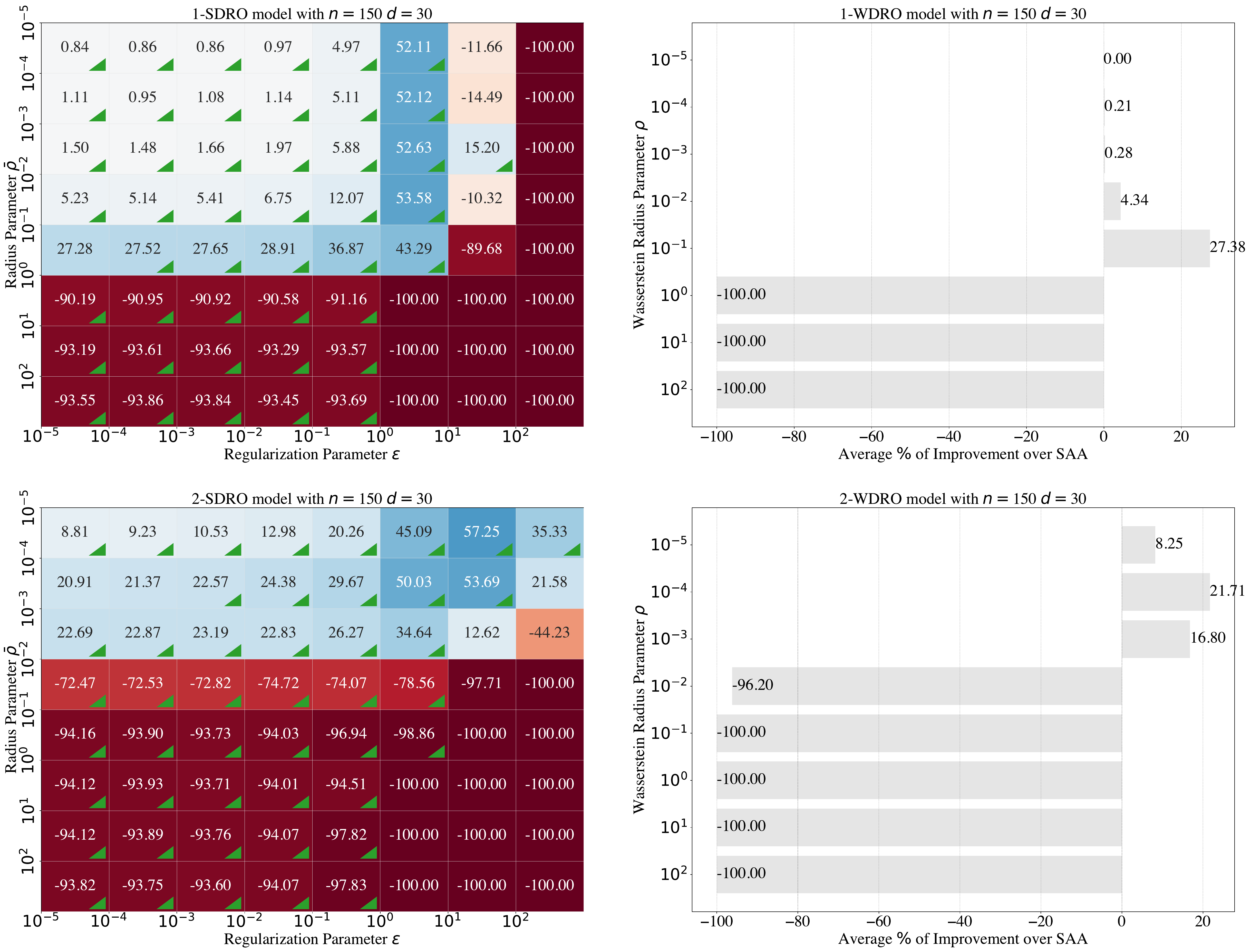}
    }
    \hfill
     \subfigure[$(n,d)=(200, 30)$]{\centering
     \includegraphics[width=0.45\textwidth]{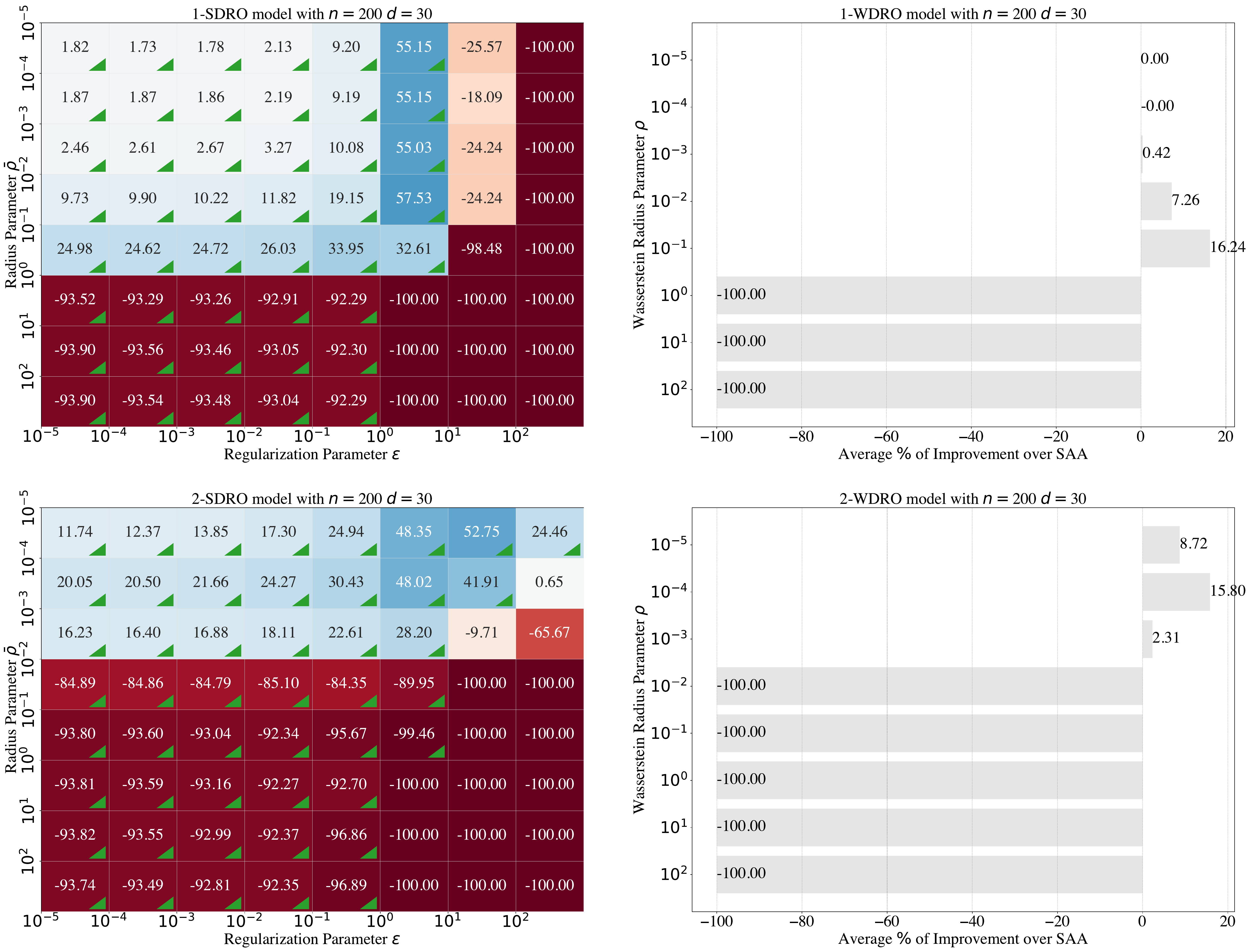}
    }
    \hfill
     \subfigure[$(n,d)=(400, 30)$]{\centering
     \includegraphics[width=0.45\textwidth]{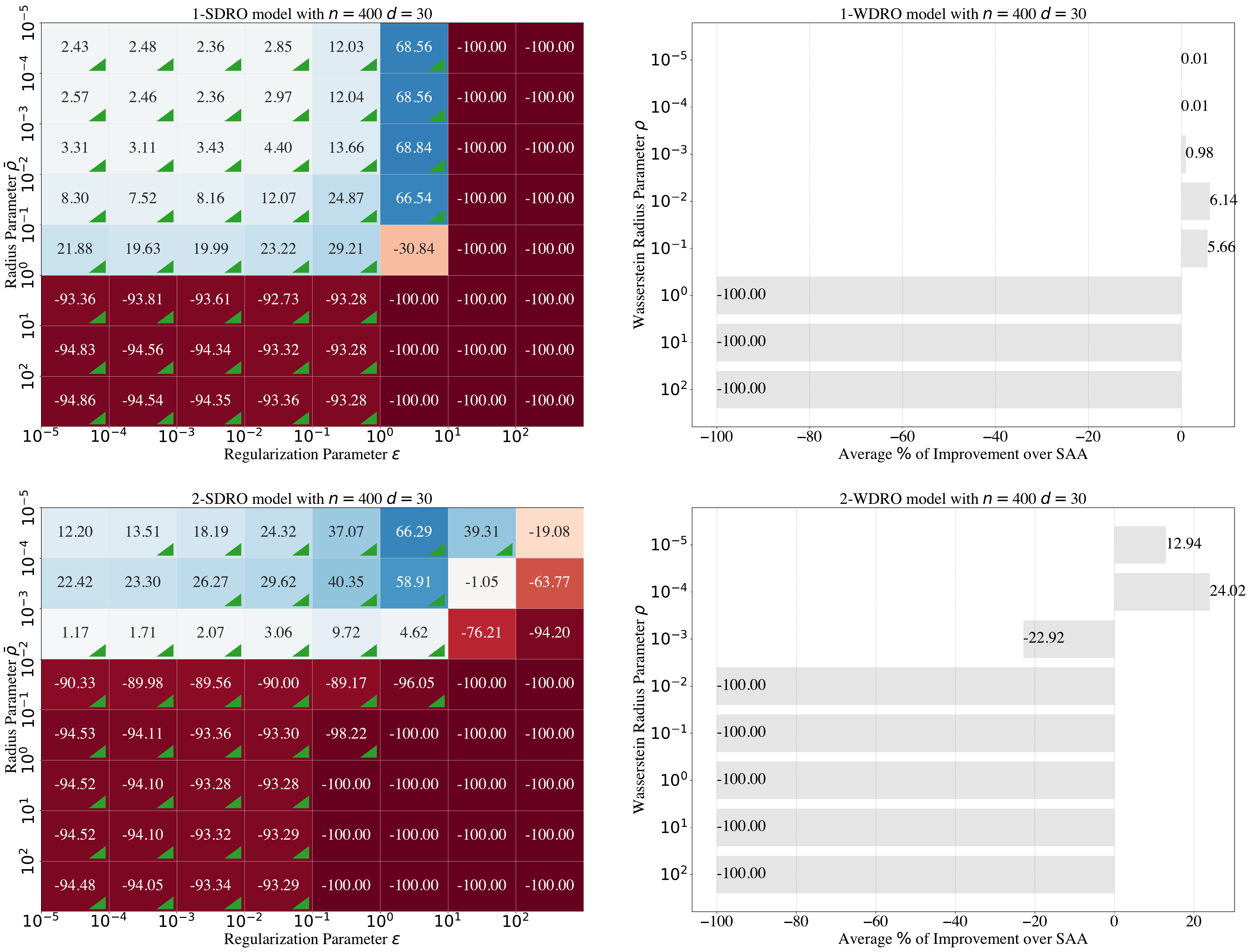}
    }
    \caption{
Additional experiment results of the portfolio optimization model for different data dimensions in heatmaps.
Here we fix the data dimension $d=30$ and vary the sample size $n\in\{50,100, 150, 200, 400\}$.
Details of these subplots follow the same setup from Fig.~\ref{fig:portfolio:heatmap:init}.
    }
    \label{Fig:portfolio:hp:1}
    }
\end{figure}
\begin{figure}[!ht]
    \centering\Jie{
     \subfigure[$(n,d)=(100,5)$]{\centering
     \includegraphics[width=0.45\textwidth]{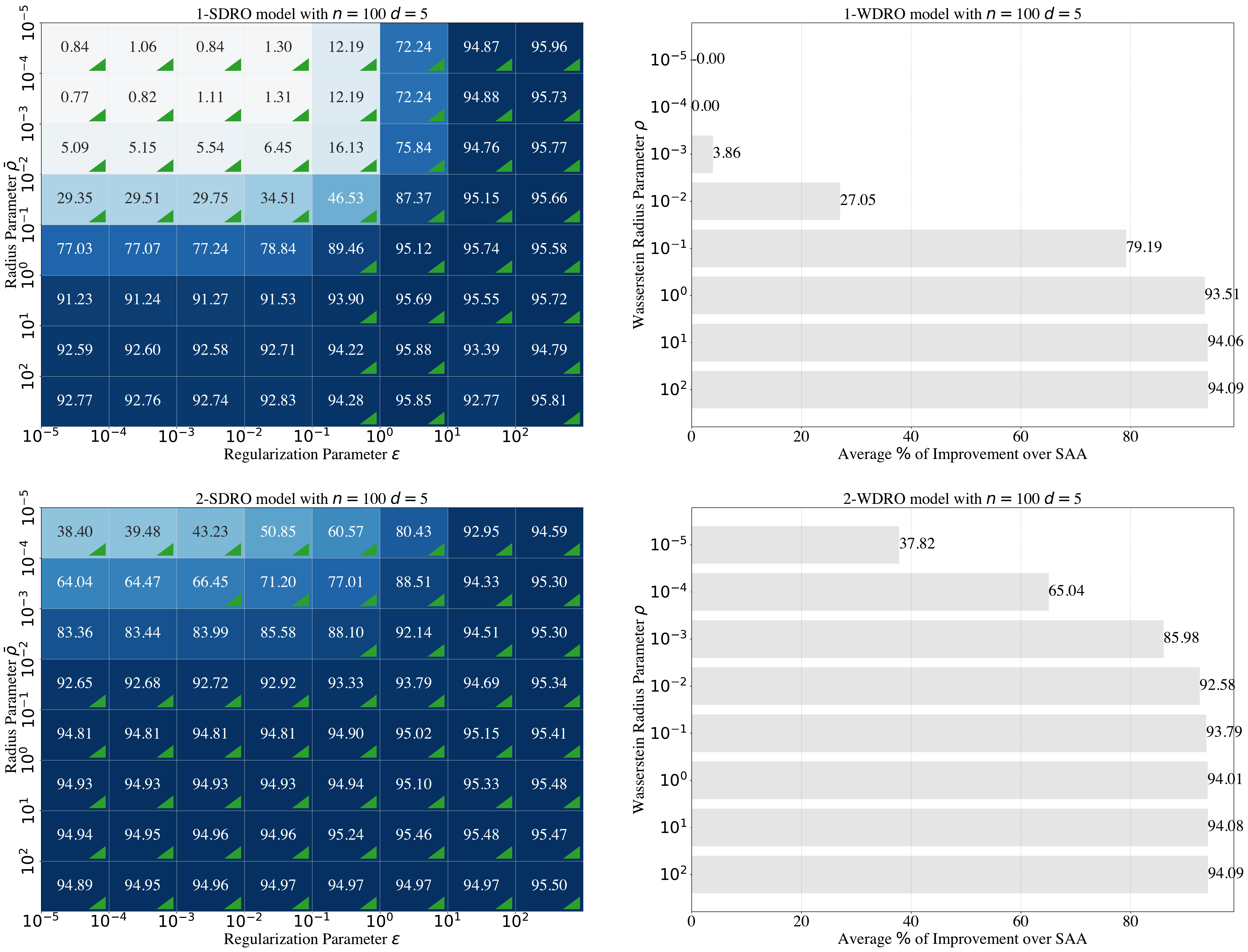}
    }
     \subfigure[$(n,d)=(100, 10)$]{\centering
     \includegraphics[width=0.45\textwidth]{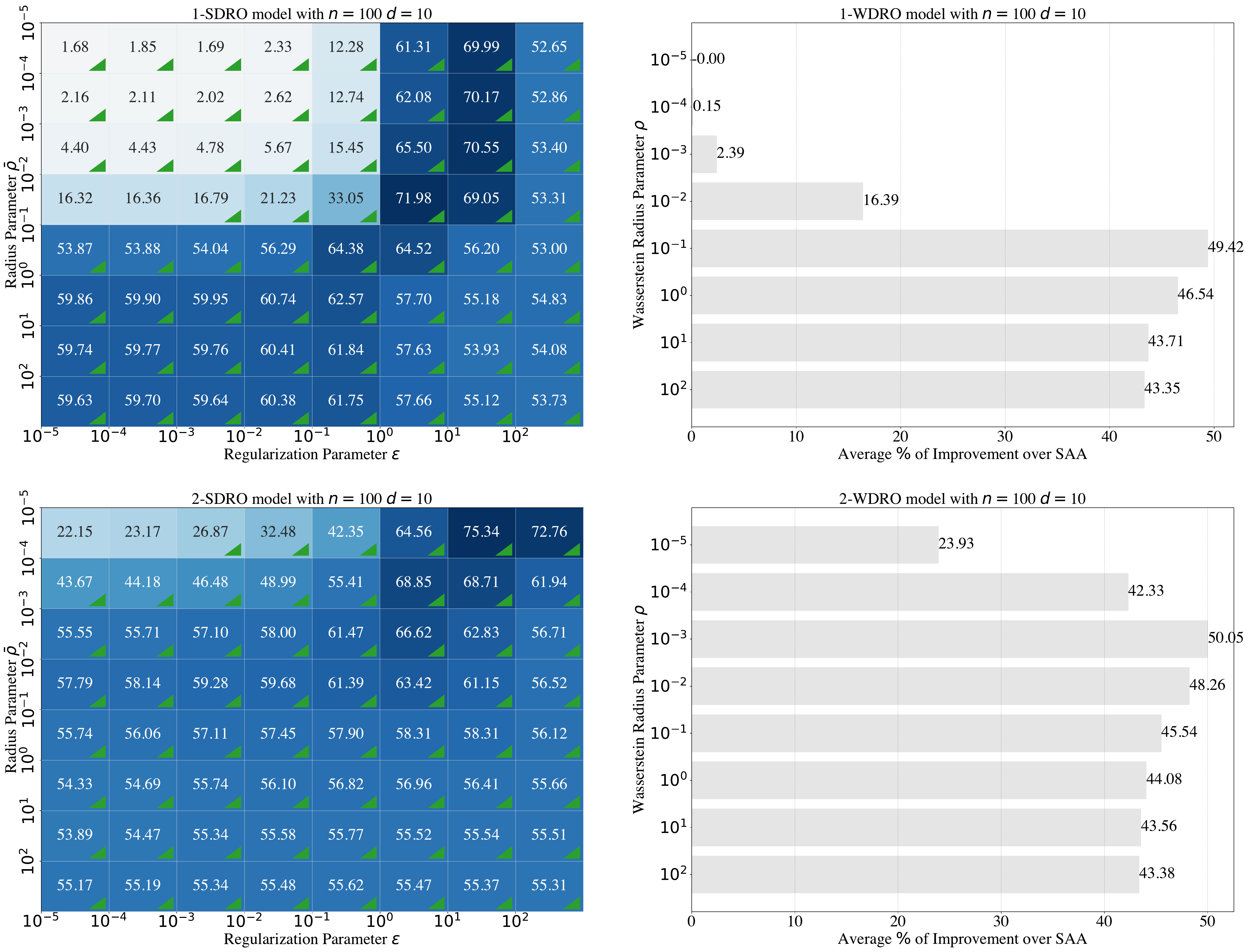}
    }
   \hfill
     \subfigure[$(n,d)=(100,20)$]{\centering
     \includegraphics[width=0.45\textwidth]{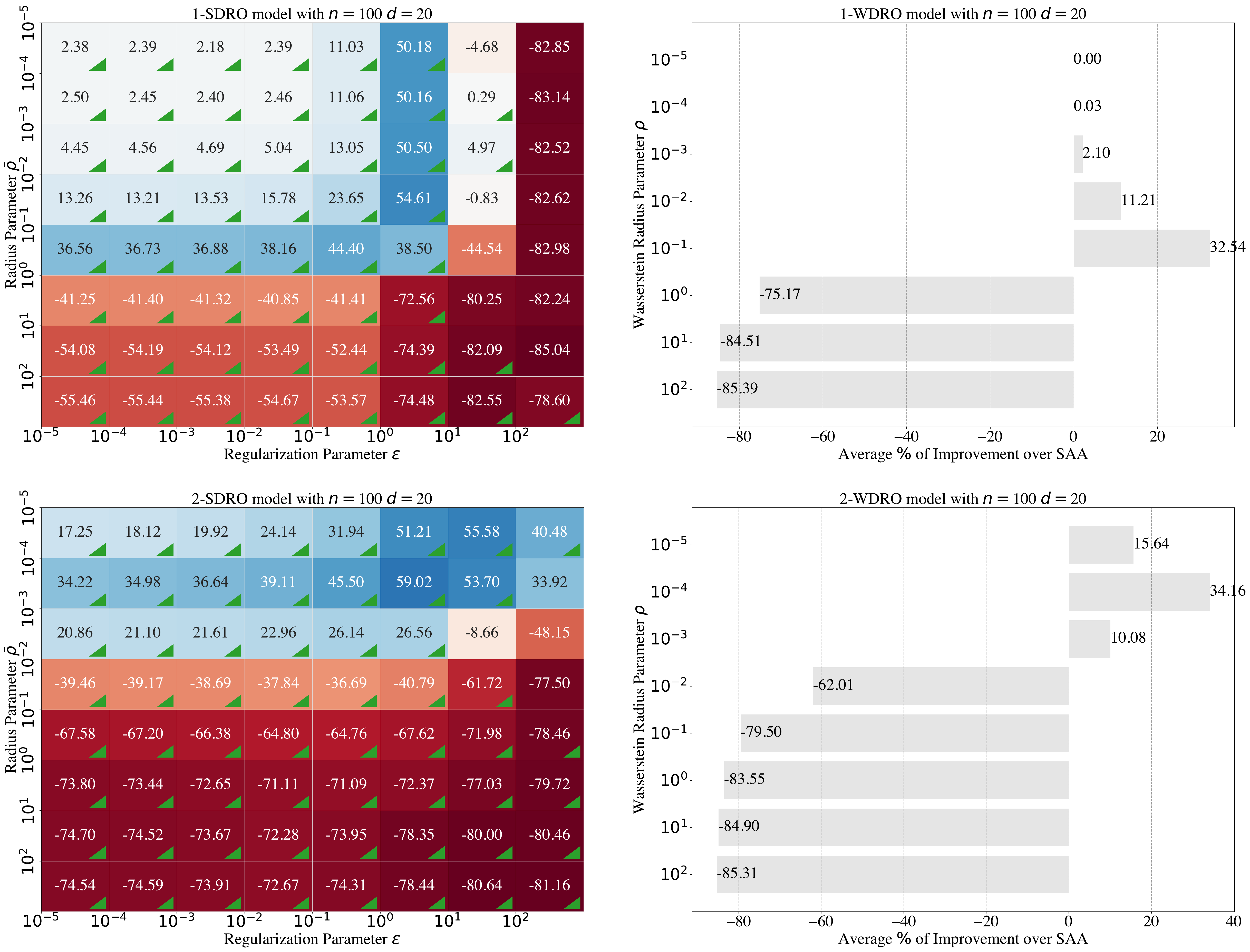}
    }
     \subfigure[$(n,d)=(100,40)$]{\centering
     \includegraphics[width=0.45\textwidth]{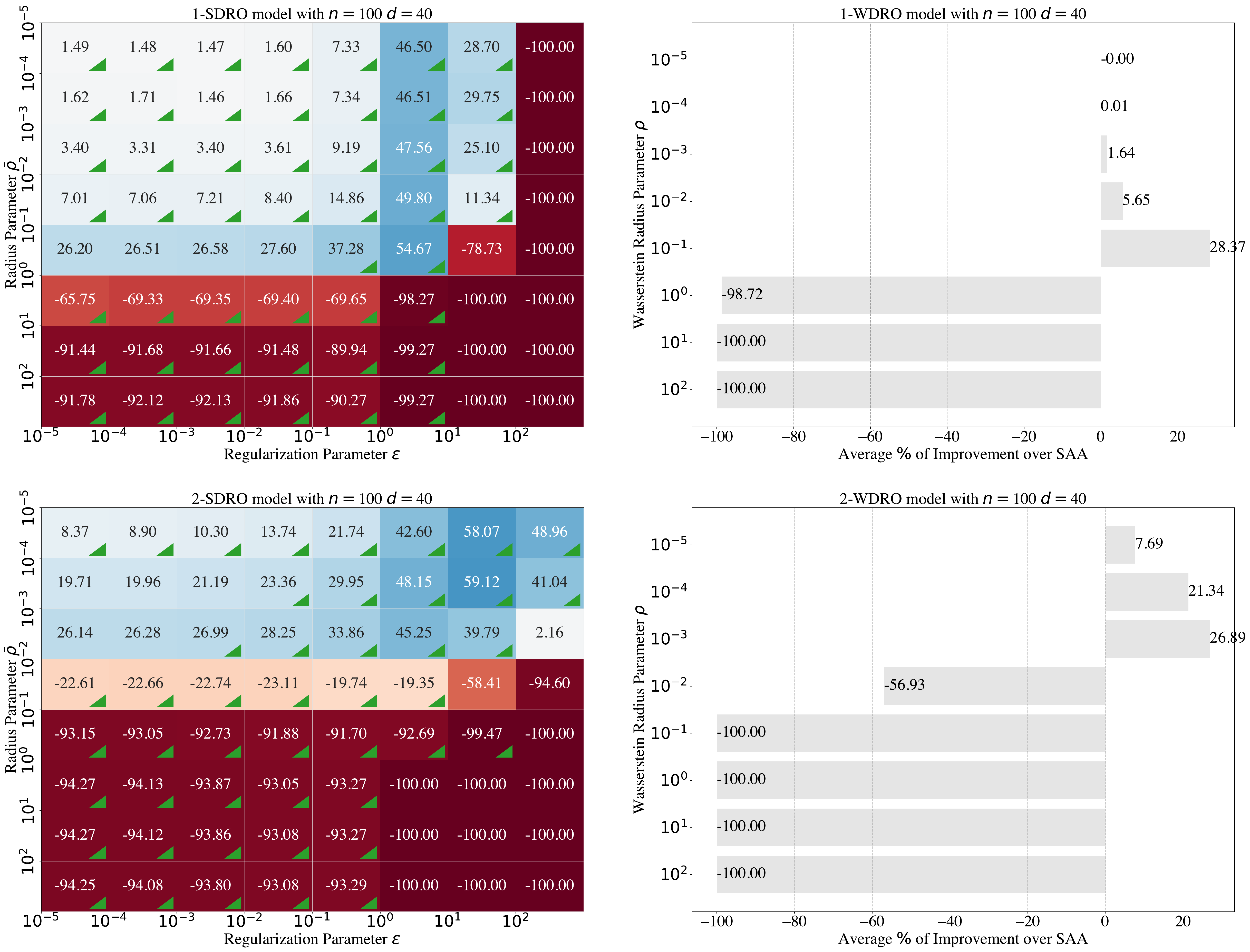}
    }
   \hfill
     \subfigure[$(n,d)=(100,80)$]{\centering
     \includegraphics[width=0.45\textwidth]{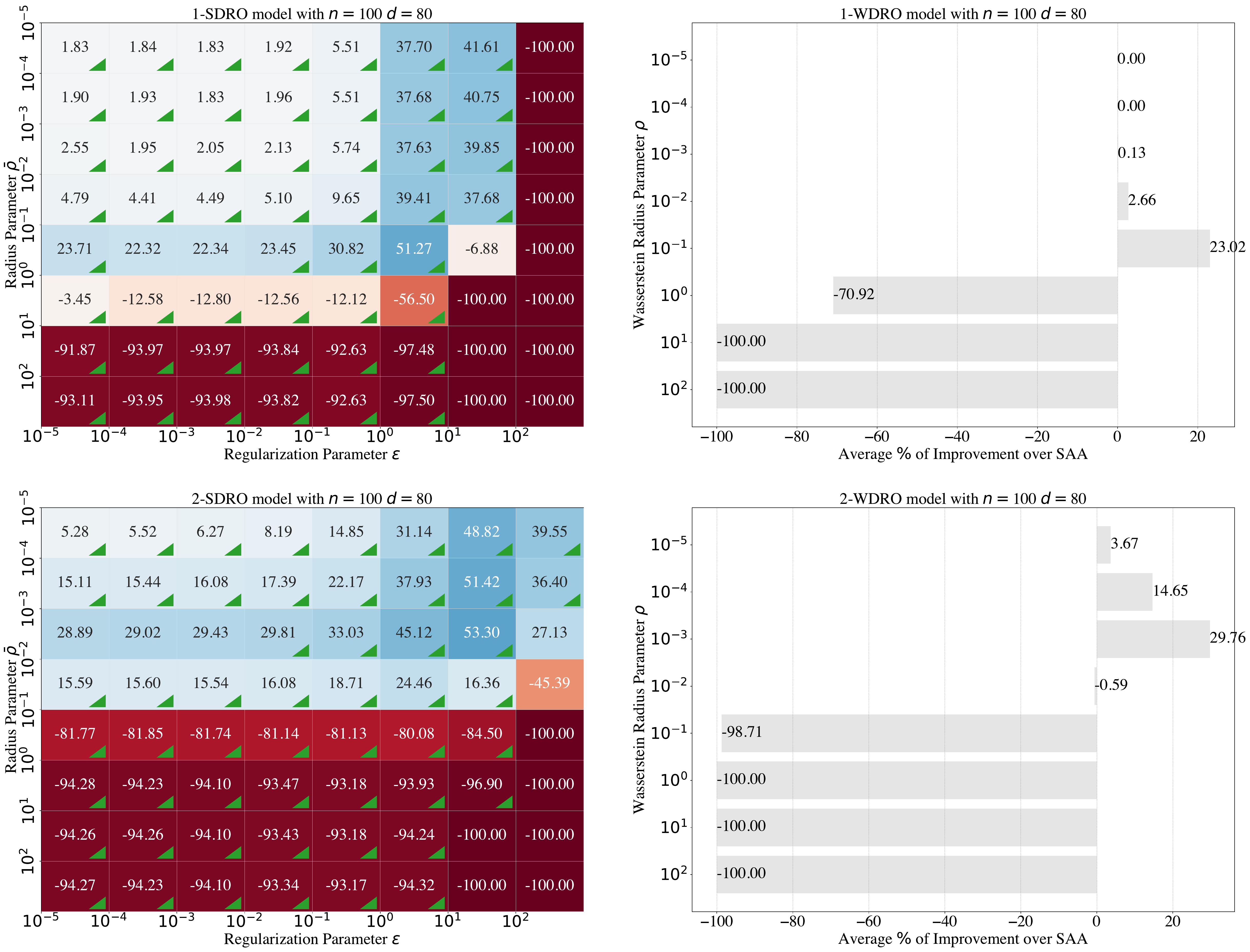}
    }
     \subfigure[$(n,d)=(100,100)$]{\centering
     \includegraphics[width=0.45\textwidth]{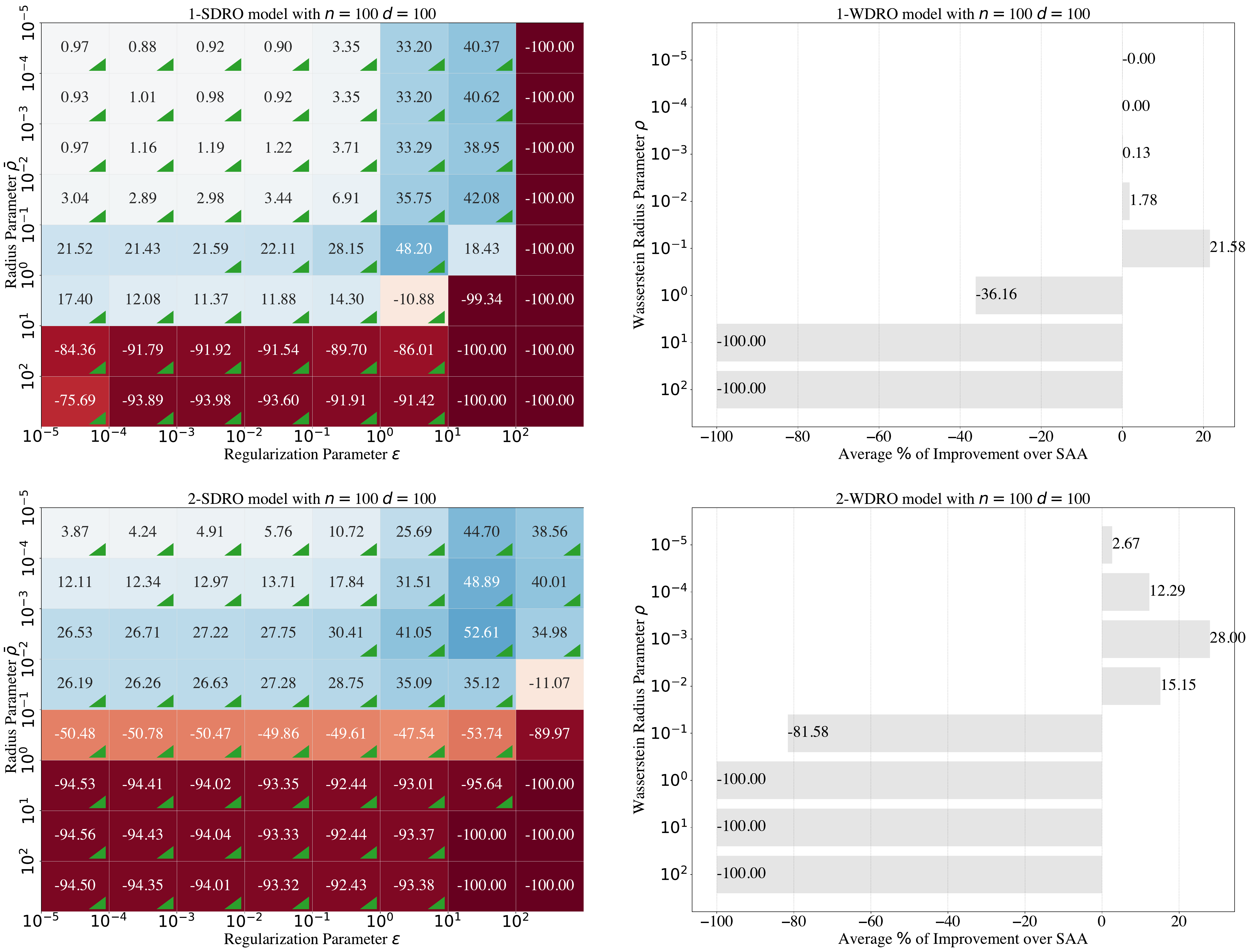}
    }
    \caption{
Additional experiment results of the portfolio optimization model for different data dimensions in heatmaps.
Here we fix the sample size $n=100$ and vary the data dimension $d\in\{5, 10, 20, 40, 80, 100\}$.
Details of these subplots follow the same setup from Fig.~\ref{fig:portfolio:heatmap:init}.
    } \label{Fig:portfolio:hp:2}
    }
\end{figure}
\begin{figure}[!ht]
    \centering
    \includegraphics[width=\linewidth]{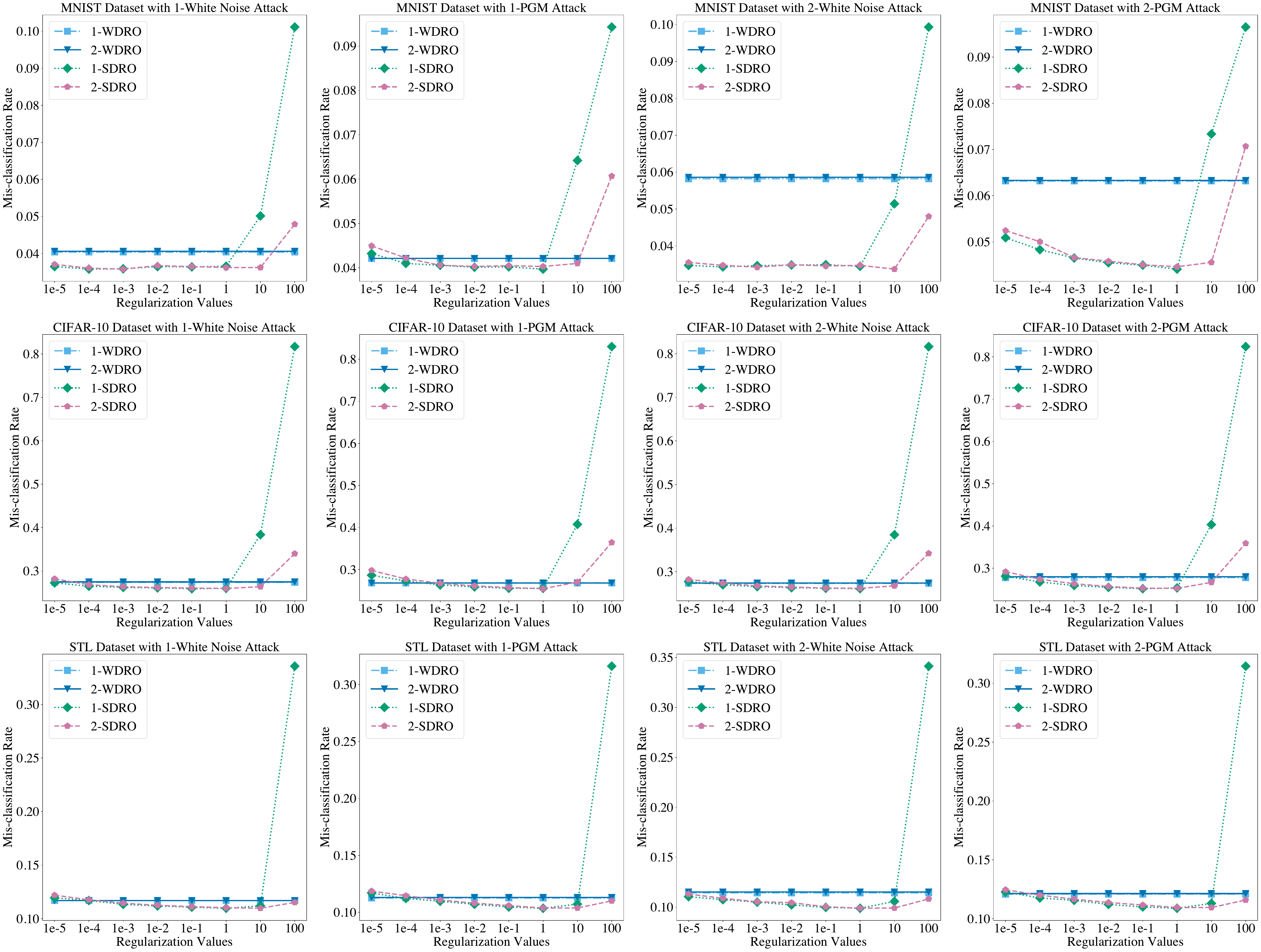}
    \caption{\Jie{Additional experiment results of the adversarial classification problem for different datasets and different types of perturbations.
    Details of these subplots follow the same setup from Fig.~\ref{fig:adversarial:ablation}.
    }}
    \label{fig:adversarial:merged}
\end{figure}
}

\clearpage
\section{Sufficient Condition for Condition~\ref{Assumption:light:tailed:f}}\label{app:sufficient condition}
  
\begin{proposition}\label{Proposition:sufficient:lambda:integrable}
Condition~\ref{Assumption:light:tailed:f} holds if there exists $p\ge1$ so that the following conditions are satisfied:
\begin{enumerate}
    \item
For any $x,y,z\in\cZ$, $c(x,y)\ge0$, and
$
(c(x,y))^{1/p}\le (c(x,z))^{1/p}+(c(z,y))^{1/p}.
$
    \item
The nominal distribution $\hP$ has a finite mean, denoted as $\ox$. 
Moreover, $\nu\{z:~0\le c(\ox,z)<\infty\}=1$ and $\Pr{}_{x\sim \hP}\{c(x,\ox)<\infty\}=1.$
    \item
    Assumption~\ref{Assumption:distance:measure:function}\ref{Assumption:distance:measure:function:3} holds, and 
there exists $\lambda>0$ such that 
$\bE_{z\sim\nu}\left[e^{f(z)/(\lambda\Reg)}e^{-2^{1-p}c(\ox,z)/\Reg}\right]<\infty.$
\end{enumerate}
\end{proposition}

We make some remarks for the sufficient conditions listed above.
The first condition can be satisfied by taking the transport cost as the $p$-th power of the metric defined on $\cZ$ for any $p\ge1$.
The second condition requires the nominal distribution $\hP$ is finite almost surely, e.g., it can be a subguassian distribution with respect to the transport cost $c$.
We first present an useful technical lemma before showing the proof of Proposition~\ref{Proposition:sufficient:lambda:integrable}.
\begin{lemma}\label{Lemma:well:defined:Q_x}
Under the first condition of Proposition~\ref{Proposition:sufficient:lambda:integrable}, for any $x\in\cZ$, it holds that
\[
\bE_{z\sim\nu}\left[
\int e^{-c(x,z)/\Reg}\right]\ge e^{-2^{p-1}c(x,\ox)/\Reg}
\bE_{z\sim\nu}\left[ e^{-2^{p-1}c(\ox,z)/\Reg}\right].
\]
\end{lemma}
\proof{Proof of Lemma~\ref{Lemma:well:defined:Q_x}.}
Based on the inequality $(a+b)^p\le 2^{p-1}(a^p+b^p)$, we can see that
\[
c(x,z)\le (c(y,z)^{1/p} + c(z,y)^{1/p})^{p}\le 2^{p-1}(c(y,z) + c(z,y)),\quad \forall x,y,z\in\cZ.
\]
Since $c(x,z)\le  2^{p-1}(c(\ox, z) + c(x,\ox))$, we can see that 
\begin{align*}
\bE_{z\sim\nu}\left[
\int e^{-c(x,z)/\Reg}\right]&\ge \exp\left( 
-2^{p-1}c(x,\ox)/\Reg
\right)\bE_{z\sim\nu}\left[ e^{-2^{p-1}c(\ox,z)/\Reg}\right].
\end{align*}
The proof is completed.
\QED\endproof

\proof{Proof of Proposition~\ref{Proposition:sufficient:lambda:integrable}.}
One can see that for any $x\in\mathrm{supp}\,\hP$, it holds that
\begin{align*}
&\mathbb{E}_{z\sim\bQ_{x,\Reg}}\left[ e^{f(z)/(\lambda\Reg)}\right]
=\bE_{z\sim \nu} \left[e^{f(z)/(\lambda\Reg)}\frac{e^{-c(x,z)/\Reg}}{\bE_{u\sim\nu}\left[e^{-c(x,u)/\Reg}\right]}\right]\\
\le&
\bE_{z\sim \nu} \left[e^{f(z)/(\lambda\Reg)}\frac{e^{-c(x,z)/\Reg}}{
\bE_{u\sim \nu}\left[e^{-2^{p-1}c(\ox,u)/\Reg}\right]
}\right]\le
\bE_{z\sim \nu} \left[e^{f(z)/(\lambda\Reg)}\frac{e^{-2^{1-p}c(\ox,z)/\Reg}e^{c(x,\ox)/\Reg}}{
\bE_{u\sim \nu}\left[e^{-2^{p-1}c(\ox,u)/\Reg}\right]
}\right]
\\
=&\frac{e^{c(x,\ox)(1 + 2^{p-1})/\Reg}}{\bE_{u\sim \nu}\left[e^{-2^{p-1}c(\ox,u)/\Reg}\right]}
\bE_{z\sim\nu}\left[e^{f(z)/(\lambda\Reg)}e^{-2^{1-p}c(\ox,z)/\Reg}\right],
\end{align*}
where the first inequality is based on the lower bound in Lemma~\ref{Lemma:well:defined:Q_x}, the second inequality is based on the triangular inequality $c(x,z)\ge 2^{1-p}c(\ox,z) - c(x,\ox)$.
Note that almost surely for all $x\in\mathrm{supp}\, \hP$, $c(x,\ox)<\infty$. 
Moreover, 
$0<
\bE_{z\sim\nu}\left[e^{-2^{p-1}c(\ox,z)/\Reg}\right]\le \bE_{z\sim\nu}\left[ e^{-c(\ox,z)/\Reg}\right]<\infty,$
where the lower bound is because $c(\ox,z)<\infty$ almost surely for all $z$, the upper bound is because $c(\ox,z)\ge0$ almost surely for all $z$.
Based on these observations, we have that 
\[
\mathbb{E}_{z\sim\bQ_{x,\Reg}}\left[ e^{f(z)/(\lambda\Reg)}\right]\le 
\frac{e^{c(x,\ox)(1 + 2^{p-1})/\Reg}}{\bE_{z\sim\nu}\left[e^{-2^{p-1}c(\ox,z)/\Reg}\right]}
\bE_{z\sim\nu}\left[e^{f(z)/(\lambda\Reg)}e^{-2^{1-p}c(\ox,z)/\Reg}\right]<\infty
\]
almost surely for all $x\sim\hP$.
\QED\endproof

\section{Proofs of Technical Results in Section~\ref{Sec:ent:duality} and \ref{Sec:example:reformulation}}\label{Sec:proof:sec:ent:duality}

\proof{Proof of Remark~\ref{Remark:Wasserstein:DRO}.}
Recall the dual objective function in \eqref{Eq:dual:SDRO:general} is
\[
v(\lambda;\Reg) = \lambda\rho +  \mathbb{E}_{x\sim\hP}\Big[ \lambda\Reg\; \log\mathbb{E}_{z\sim\nu}\big[
   e^{
  (f(z) - \lambda c(x,z))/(\lambda\Reg)}
  \big] \Big].
\]
We take limit for the second term in $v(\lambda; \Reg)$ to obtain:
\begin{align*}
&\lim_{\Reg\to0}\mathbb{E}_{x\sim\hP}\Big[ \lambda\Reg\; \log\mathbb{E}_{z\sim\nu}\big[
   e^{
  (f(z) - \lambda c(x,z))/(\lambda\Reg)}
  \big] \Big]
=
\mathbb{E}_{x\sim\hP}\left[ \lim_{\beta\to\infty}\frac{\lambda}{\beta} \log\mathbb{E}_{z\sim\nu}\big[
   e^{
  (f(z) - \lambda c(x,z))\beta/\lambda}
  \big] \right]\\
=&
\mathbb{E}_{x\sim\hP}\left[ \lim_{\beta\to\infty}
\lambda\nabla_{\beta}
\log\mathbb{E}_{z\sim\nu}\big[
   e^{
  (f(z) - \lambda c(x,z))\beta/\lambda}
  \big]\right]
=\mathbb{E}_{x\sim\hP}\left[ \lim_{\beta\to\infty}
\frac{
\mathbb{E}_{z\sim\nu}\big[
   e^{
  (f(z) - \lambda c(x,z))\beta/\lambda}
  \big[f(z) - \lambda c(x,z)\big]
  \big]
}{
\mathbb{E}_{z\sim\nu}\big[
   e^{
  (f(z) - \lambda c(x,z))\beta/\lambda}
  \big]
}
\right]\\
=&\mathbb{E}_{x\sim\hP}\Big[
\sup_{z\in\mathrm{supp}\,\nu}
\big\{f(z) - \lambda c(x,z)\big\}\Big].
\end{align*}
Particularly, when $\mathrm{supp}\,\nu=\cZ$, it holds that 
\[
\sup_{z\in\mathrm{supp}\,\nu}
\big\{f(z) - \lambda c(x,z)\big\}= \sup_{z\in\cZ}
\big\{f(z) - \lambda c(x,z)\big\}
\]
and in this case the dual objective function of the Sinkhorn DRO problem converges into that of the Wasserstein DRO problem.
\QED\endproof

\proof{Proof of Example~\ref{Example:DRO:LR}.}
In this example, the dual objective becomes
\begin{equation}
V_{\Dual} = \inf_{\lambda\ge0}~\left\{
\lambda\orho + \bE_{(a,b)\sim\hP}\left[
\lambda\Reg\log
\bE_{a'\sim\cN(a, \Reg I_d)}~\left[\exp\left(
\frac{(\theta\trans a' - b)^2}{\lambda\Reg}
\right)\right]\right]
\right\}.
\label{Eq:dual:LR}
\end{equation}
Specially, for any $a\in\bR^d, b\in\bR, \theta\in\bR^d$, it holds that
\begin{align*}
&\lambda\Reg\log\left( 
\bE_{a'\sim\cN(a, \Reg I_d)}~\exp\left(
\frac{(\theta\trans a' - b)^2}{\lambda\Reg}
\right)
\right)
=
\lambda\Reg\log\left( 
\bE_{\Delta_a\sim\cN(0, I_d)}~\exp\left(
\frac{\big[(\theta\trans a - b) + (\sqrt{\Reg}\theta)\trans \Delta_a\big]^2}{\lambda\Reg}
\right)
\right)\\
=&(\theta\trans a - b)^2 + \lambda\Reg\log\left(
\underbrace{
\bE_{\Delta_a\sim\cN(0, I_d)}~\exp\left( 
\frac{\Reg(\theta\trans\Delta_a)^2 - 2(b - \theta\trans a)\sqrt{\Reg}\theta\trans\Delta_a}{\lambda\Reg}
\right)}_{\text{(I)}}
\right).
\end{align*}
The term (I) can be simplified using the integral of exponential functions method:
\[
\text{(I)}=\left\{ 
\begin{aligned}
\det\left( 
I - \frac{2\theta\theta\trans}{\lambda}
\right)^{-1/2}\exp\left( 
2\frac{(\theta\trans a - b)^2}{\lambda^2\Reg}\theta\trans A^{-1}\theta
\right),&\quad\text{when $\|\theta\|_2^2<\frac{\lambda}{2}$,}\\
\infty,&\quad\text{otherwise},
\end{aligned}
\right.
\]
where the matrix $A = I - \frac{2\theta\theta\trans}{\lambda}$.
Finally, we obtain that if $\|\theta\|_2^2<\frac{\lambda}{2}$,
\begin{align*}
&\lambda\Reg\log\left( 
\bE_{a'\sim\cN(a, \Reg I_d)}~\exp\left(
\frac{(\theta\trans a' - b)^2}{\lambda\Reg}
\right)
\right)
=(\theta\trans a - b)^2 +\frac{(\theta\trans a - b)^2}{\frac{1}{2}\lambda\|\theta\|_2^{-2}-1} - \frac{\lambda\Reg}{2}\log\det\left( 
I - \frac{2\theta\theta\trans}{\lambda}
\right).%
\end{align*}
Substituting this expression into \eqref{Eq:dual:LR} gives the desired result.
\QED\endproof

\proof{Proof of Corollary~\ref{Corollary:conic:reformulation}.}
We now introduce the epi-graphical variables $s_i, i=1,\ldots,n$ to reformulate $V_{\Dual}$ as
\[
V_{\Dual}=\left\{
\begin{aligned}
\inf_{\lambda\ge0, s_i}&\quad \lambda\orho + \frac{1}{n}\sum_{i=1}^ns_i\\
\mbox{s.t.}&\quad \lambda\Reg\log
\mathbb{E}_{z\sim\bQ_{\hx_i,\Reg}}\left[ 
e^{f(z)/(\lambda\Reg)}
\right]
\le s_i, \forall i
\end{aligned}
\right.
\]
For fixed $i$, the $i$-th constraint can be reformulated as 
\[
\begin{aligned}
&\left\{
\exp\left(
\frac{s_i}{\lambda\Reg}
\right)
\ge 
\mathbb{E}_{z\sim\bQ_{\hx_i,\Reg}}\left[ 
e^{f(z)/(\lambda\Reg)}
\right]
\right\}
=
\left\{
1
\ge 
\mathbb{E}_{z\sim\bQ_{\hx_i,\Reg}}\left[ 
e^{(f(z)-s_i)/(\lambda\Reg)}
\right]
\right\}\\
=&
\left\{
\lambda\Reg
\ge 
\mathbb{E}_{z\sim\bQ_{\hx_i,\Reg}}\left[ 
\lambda\Reg e^{(f(z)-s_i)/(\lambda\Reg)}
\right]
\right\}\\
=&
\left\{
\lambda\Reg
\ge 
\sum_{\ell=1}^{L_{\max}}\bQ_{\hx_i,\Reg}(z_{\ell})a_{i,\ell}
\right\}
\bigcap
\left\{ 
a_{i,\ell}\ge \lambda\Reg\exp\left(\frac{f(z_{\ell}) - s_i}{\lambda\Reg}\right), \forall \ell
\right\},
\end{aligned}
\]
where the second constraint set can be formulated as $(\lambda\Reg, a_{i,\ell}, f(z_{\ell})-s_i)\in \mathcal{K}_{\exp}.$
Substituting this expression into $V_{\Dual}$ completes the proof.
\QED\endproof

\section{Proofs of Technical Results in Section~\ref{Proof:Theorem:strong:duality:main}}\label{Sec:appendix:proof}
We rely on the following technical lemma to derive our strong duality result.
\begin{lemma}{(\citep[Section~2.1]{hu2013kullback} or \citep{shapiro2017distributionally})}\label{Lemma:expression:vxlambda}
For fixed $\tau$ and a reference measure $\nu\in\cM(\cZ)$, consider the optimization problem
\begin{equation}\label{Eq:optimization:lemma:vlambda}
v(\tau)=\sup_{\bP\in\cP(\cZ)}~\left\{\mathbb{E}_{z\sim\bP}\left[ 
f(z) - \tau\log\left(
\frac{\diff\bP(z)}{\diff\nu(z)}
\right)
\right]\right\}.
\end{equation}
Suppose there exists a probability measure $\bQ\in\cP(\cZ)$ such that $\bQ\ll \nu$.
\begin{enumerate}
    \item
When $\tau=0$, 
\[
v(0)=\underset{\nu}{\esssup}(f)\triangleq\inf\{t\in\mathbb{R}:~
\nu\{f(z)>t\}=0
\}.
\]
    \item
When $\tau>0$ and 
\[
\mathbb{E}_{z\sim\nu}\left[ 
e^{f(z)/\tau}
\right]<\infty,
\]
it holds that
\[
v(\tau)=\tau\log\left(
\mathbb{E}_{z\sim \nu}\left[ 
e^{f(z)/\tau}
\right]
\right),
\]
and $\lim_{\tau\downarrow 0}v(\tau)=v(0)$.
The optimal solution in \eqref{Eq:optimization:lemma:vlambda} has the expression
\[
\diff\bP(z) = \frac{e^{f(z)/\tau}}{\bE_{u\sim\nu}\left[e^{f(u)/\tau}\right]}\diff\nu(z).
\]
\item
When $\tau>0$ and 
\[
\mathbb{E}_{z\sim\nu}\left[ 
e^{f(z)/\tau}
\right]=\infty,
\]
we have that $v(\tau)=\infty$.
\end{enumerate}
\end{lemma}

\Jie{
\begin{lemma}[Measurability of $v_x(\lambda)$]\label{Lemma:measure:vx}
Assume Assumptions~\ref{Assumption:distance:measure:function}\ref{Assumption:distance:measure:function:1}, \ref{Assumption:distance:measure:function}\ref{Assumption:distance:measure:function:2}, \ref{Assumption:distance:measure:function}\ref{Assumption:distance:measure:function:3} hold.
For fixed $\lambda\ge0$, define the function $v_x(\lambda):~\mathrm{supp}\,\hP\to\mathbb{R}\cup\{+\infty\}$ as
\[
v_x(\lambda) = \sup_{\gamma_x\in\cP(\cZ)}~
\left\{\mathbb{E}_{z\sim\gamma_x}\left[ 
f(z) - \lambda c(x,z) - \lambda\Reg\log\left(\frac{\diff \gamma_x(z)}{\diff\nu(z)}\right)
\right]\right\}.
\]
The function $v_x(\lambda)$ is measurable with respect to $x\sim\hP$ regardless of the choice of $\lambda\ge0$.
\end{lemma}
\proof{Proof of Lemma~\ref{Lemma:measure:vx}.}
When $\lambda=0$, by Lemma~\ref{Lemma:expression:vxlambda}, it holds that
\[
v_x(\lambda) = \esssup_{\nu}(f),
\]
which is a constant independent of $x$, which is clearly measurable.
When $\lambda>0$ and satisfies Condition~\ref{Assumption:light:tailed:f}, by Lemma~\ref{Lemma:expression:vxlambda}, it holds that
\[
v_x(\lambda) = 
\lambda\Reg\log\mathbb{E}_{z\sim\nu}\big[
   e^{
  (f(z) - \lambda c(x,z))/(\lambda\Reg)}
  \big]<\infty.
\]
As loss function $f$ and cost function $c$ are both measurable, by conditioning Lemma~\citep[Lemma~2.11]{kallenberg1997foundations}, $v_x(\lambda)$ is measurable.
When $\lambda>0$ such that the event 
\[
E=\left\{
x:~\mathbb{E}_{z\sim\bQ_{x,\Reg}}\left[ 
e^{f(z)/(\lambda\Reg)}
\right]=\infty
\right\}
=
\left\{
x:~
\mathbb{E}_{z\sim\nu}\big[
   e^{
  (f(z) - \lambda c(x,z))/(\lambda\Reg)}
  \big]=\infty
\right\}
\]
satisfies $\hP(E)>0$, by Lemma~\ref{Lemma:expression:vxlambda}, it holds that
\[
v_x(\lambda)=
\left\{ 
\begin{aligned}
\lambda\Reg\log\mathbb{E}_{z\sim\nu}\big[
   e^{
  (f(z) - \lambda c(x,z))/(\lambda\Reg)}
  \big]<\infty,&\quad \text{if }x\in E^c,\\
\infty,&\quad\text{if }x\in E.
\end{aligned}
\right.
\]
For fixed $\alpha\in\mathbb{R}$, the level set 
\[
\{x:~v_x(\lambda)\ge \alpha\} = \{x\in E^c:~v_x(\lambda)\ge \alpha\}\cup E
=\left\{x\in E^c:~\lambda\Reg\log\mathbb{E}_{z\sim\nu}\big[
   e^{
  (f(z) - \lambda c(x,z))/(\lambda\Reg)}
  \big]\ge \alpha\right\}\cup E,
\]
which is clearly a measurable set, and therefore $v_x(\lambda)$ is measurable.
The proof is completed.
\QED
\endproof
}

\proof{Proof of Lemma~\ref{Lemma:refor:Inf-W-E}.}
Recall from \eqref{Eq:V:primal:revision} that
\[
V_{\Primal}=\sup_{\{\gamma_x\}_{x\in\mathrm{supp}\,\hP}\subset \cP(\cZ)} \left\{
\bE_{x\sim\hP}\mathbb{E}_{z\sim\gamma_x}[f(z)]:\
\bE_{x\sim\hP}
\mathbb{E}_{z\sim \gamma_x}\left[ 
c(x, z) + \Reg\log\left(
\frac{\diff\gamma_x(z)}{\diff\nu(z)}
\right)
\right]
\le \rho
\right\}.
\]
Based on the change-of-measure identity $\log\left(\frac{\diff\gamma_x(z)}{\diff\nu(z)}\right)=\log\left(\frac{\diff\bQ_{x,\Reg}(z)}{\diff\nu(z)}\right) + \log\left(\frac{\diff\gamma_x(z)}{\diff\bQ_{x,\Reg}(z)}\right)$ and the expression of $\bQ_{x,\Reg}$, the constraint can be reformulated as 
\begin{equation*}
\bE_{x\sim\hP}
\mathbb{E}_{z\sim\gamma_x}\left[ 
c(x, z) + \Reg\log\left(\frac{e^{-c(x,z)/\Reg}}{\int e^{-c(x,u)/\Reg}\diff\nu(u)}\right) + \Reg\log\left(
\frac{\diff\gamma_x(z)}{\diff\bQ_{x,\Reg}(z)}
\right)
\right]
\le \rho.
\end{equation*}
Combining the first two terms within the expectation term and substituting the expression of $\orho$, it is equivalent to
\[
\Reg
\bE_{x\sim\hP}\mathbb{E}_{z\sim\gamma_x}\left[ 
\log\left( 
\frac{\diff\gamma_x(z)}{\diff\bQ_{x,\Reg}(z)}
\right)
\right]
\le \orho.
\]
In summary, the primal problem~\eqref{Inf-W-E} can be reformulated as a generalized KL-divergence DRO problem
\[
V_{\Primal}=\sup_{\{\gamma_x\}_{x\in\mathrm{supp}\,\hP}\subset \cP(\cZ)} \left\{\bE_{x\sim\hP}\mathbb{E}_{z\sim\gamma_x}[f(z)]:\
 \Reg
\bE_{x\sim\hP}\mathbb{E}_{z\sim\gamma_x}\left[ 
\log\left( 
\frac{\diff\gamma_x(z)}{\diff\bQ_{x,\Reg}(z)}
\right)
\right]
\le \orho
\right\}.\tag*{\QED}
\]
\endproof

In the remaining of this subsection, we provide the full proof of Theorem~\ref{Theorem:strong:duality}.
We first show that the dual minimizer exists.

\begin{lemma}[Existence of Dual Minimizer]\label{Lemma:finite:dual}
Suppose $\orho>0$ and Condition~\ref{Assumption:light:tailed:f} is satisfied, then the dual minimizer $\lambda^*$ exists, which either equals to $0$ or satisfies Condition~\ref{Assumption:light:tailed:f}.
\end{lemma}
\proof{Proof of Lemma~\ref{Lemma:finite:dual}.}
We first show that $\lambda^*<\infty$.
Denote by $v(\lambda)$ the objective function for the dual problem:
\[
v(\lambda) = 
\lambda\orho+\lambda\Reg\,\mathbb{E}_{x\sim\hP}\Big[\log
\mathbb{E}_{z\sim\bQ_{x,\Reg}}\big[ 
e^{f(z)/(\lambda\Reg)}
\big]\Big].
\]
The integrability condition for the dominated convergence theorem is satisfied, which implies
\begin{align*}
&\lim_{\lambda\to\infty}
\lambda\Reg\,\mathbb{E}_{x\sim\hP}\Big[\log
\mathbb{E}_{z\sim\bQ_{x,\Reg}}\big[ 
e^{f(z)/(\lambda\Reg)}
\big]\Big]
=
\mathbb{E}_{x\sim\hP}\left[ \lim_{\beta\to0}\frac{\Reg}{\beta} \log\mathbb{E}_{z\sim\bQ_{x,\Reg}}\big[ 
e^{\beta f(z)/\Reg}
\big] \right]\\
=&\mathbb{E}_{x\sim\hP}\left[ \lim_{\beta\to0}
\Reg\nabla_{\beta}
\log\mathbb{E}_{z\sim\bQ_{x,\Reg}}\big[ 
e^{\beta f(z)/\Reg}
\big]\right]
=\bE_{x\sim\hP}\left[\lim_{\beta\to0}
\frac{
\bE_{z\sim \bQ_{x,\Reg}}\left[f(z) e^{\beta f(z)/\Reg}\right]
}{\mathbb{E}_{z\sim\bQ_{x,\Reg}}\left[ 
e^{\beta f(z)/\Reg}
\right]}
\right]\\
=&
\bE_{x\sim\hP}
\mathbb{E}_{z\sim \bQ_{x,\Reg}}[f(z)],
\end{align*}
where the first equality follows from the change-of-variable technique with $\beta=1/\lambda$, 
the second equality follows from the definition of derivative,
the third and the last equality follows from the dominated convergence theorem.
As a consequence, as long as $\orho>0$, we have
$\lim_{\lambda\to\infty}v(\lambda)=\infty.$
We can take $\lambda$ satisfying Condition~\ref{Assumption:light:tailed:f} and then $v(\lambda)<\infty$.
This, toegther with the fact that $v(\cdot)$ is continuous, guarantees the existence of the dual minimizer.
Hence $\lambda^*<\infty$, which implies that either $\lambda^*=0$ or $\lambda^*$ satisfies Condition~\ref{Assumption:light:tailed:f}.
\QED\endproof

Next, we establish first-order optimality condition for cases $\lambda^*>0$ or $\lambda^*=0$, corresponding to whether the Sinkhorn distance constraint in \eqref{Inf-W-E} is binding or not.
Lemma~\ref{Lemma:necessary:lambda:0} below presents a necessary and sufficient condition for the dual minimizer $\lambda^*=0$, corresponding to the case where the Sinkhorn distance constraint in \eqref{Inf-W-E} is not binding.

\begin{lemma}[Necessary and Sufficient Condition for $\lambda^*=0$]\label{Lemma:necessary:lambda:0}
Suppose $\orho>0$ and Condition~\ref{Assumption:light:tailed:f} is satisfied, then the dual minimizer $\lambda^*=0$ if and only if all the following conditions hold:
\begin{enumerate}
\item\label{Lemma:necessary:lambda:0:I}
$\esssup_\nu
~f\triangleq \inf\{t:\,\nu\{f(z)>t\}=0\}<\infty$.
    \item\label{Lemma:necessary:lambda:0:II}
$\orho'=\orho+\Reg\bE_{x\sim \hP}\left[\log
\mathbb{E}_{z\sim\bQ_{x,\epsilon}}[1_A(z)]
\right]
\ge0$, where $A:=\{z:\,f(z)=\esssup_\nu~f\}$.
\end{enumerate}
\end{lemma}
Recall that we have the convention that  the dual objective evaluated at $\lambda=0$ equals $\esssup_\nu~f$. Thus Condition~\ref{Lemma:necessary:lambda:0:I} ensures that the dual objective function evaluated at the minimizer is finite.
When the minimizer $\lambda^*=0$, the Sinkhorn ball should be large enough to contain at least one distribution with objective value $\esssup_\nu~f$, and Condition~\ref{Lemma:necessary:lambda:0:II} characterizes the lower bound of $\orho$.
\proof{Proof of Lemma~\ref{Lemma:necessary:lambda:0}.}
Suppose the dual minimizer $\lambda^*=0$, then taking the limit of the dual objective function gives
\[
\lim_{\lambda\to0}~v(\lambda) = \bE_{x\sim\hP}\big[H^u(x)\big]<\infty,
\]
where $H^u(x) := \inf\{t:~\bQ_{x,\Reg}\{f(z)>t\} = 0\} \triangleq \underset{\bQ_{x,\Reg}}{\esssup}~f.$
For notational simplicity we take $H^u=\underset{\nu}{\esssup}~f$.
One can check that $H^u(x)\equiv H^u$ for any $x\in\mathrm{supp}\,\hP$: 
for any $t$ so that $\bQ_{x,\Reg}\{f(z)>t\}=0$, we have that
\[
\bE_{z\sim\nu}\Big[1\{f(z)>t\}e^{-c(x,z)/\Reg}\Big]=0,
\]
which, together with the fact that $\nu\{c(x,z)<\infty\}=1$ for fixed $x$, implies
\[
\bE_{z\sim\nu}[1\{f(z)>t\}]=0.
\]
On the contrary, for any $t$ so that $\nu\{f(z)>t\}=0$, we have that 
\[
0\le 
\bE_{z\sim\nu}\Big[1\{f(z)>t\}e^{-c(x,z)/\Reg}\Big]
\le \bE_{z\sim\nu}[1\{f(z)>t\}]=0,
\]
where the second inequality is because that $\nu\{c(x,z)\ge0\}=1$. 
As a consequence, $\bQ_{x,\Reg}\{f(z)>t\}=0$.
Hence we can assert that $H^u(x)=H^u$ for all $x\in\mathrm{supp}\,\hP$, which implies
\[
\lim_{\lambda\to0}~v(\lambda) = H^u < \infty.
\]

Then we show that almost surely for all $x$,
\[
\mathbb{E}_{z\sim\bQ_{x,\Reg}}[1_{A}(z)]>0,\quad\text{where } A=\{z:~f(z)=H^u\}.
\]
Denote by $D$ the collection of samples $x$ so that $\mathbb{E}_{z\sim\bQ_{x,\Reg}}[1_{A}(z)]=0$.
Assume the condition above does not hold, which means that $\hP\{D\}>0$.
For any $\tau>0$ and $x\in D$, there exists $H^l(x)<H^u$ such that 
\[
0<\fh_{x}:=\mathbb{E}_{z\sim\bQ_{x,\Reg}}[1_{B(x)}(z)]\le \tau, \quad \text{where }B(x)=\{z:~H^l(x)\le f(z)\le H^u\}.
\]
Define $H^{\mathrm{gap}}(x)=H^u-H^l(x)$, $\fh_{x}^c=1-\fh_{x}$.
Then we find that for $x\in D$,
\begin{align*}
v_{x}(\lambda)&=\lambda\Reg\log\left(
\mathbb{E}_{z\sim\bQ_{x,\Reg}}\left[e^{f(z)/(\lambda\Reg)}1_{B(x)}(z)\right]
+
\mathbb{E}_{z\sim\bQ_{x,\Reg}}\left[e^{f(z)/(\lambda\Reg)}1_{B(x)^c}(z)\right]
\right)\\
&\le H^u + \lambda\Reg\log\left(
\fh_x
+
e^{-H^{\mathrm{gap}}(x)/(\lambda\Reg)}
\fh_{x}^c
\right).
\end{align*}
Since $\hP\{D\}>0$, the dual objective function for $\lambda>0$ is upper bounded as
\begin{align*}
v(\lambda)&=\lambda\orho + \bE_{x\sim \hP}[v_x(\lambda)]\\
&\le H^u + \lambda\orho  + \lambda\Reg\bE_{x\sim\hP}\left[\log\left(
\fh_x
+
e^{-H^{\mathrm{gap}}(x)/(\lambda\Reg)}
\fh_{x}^c
\right)1_D(x)\right].
\end{align*}
We can see that 
\[
\lim_{\lambda\to0}~\lambda\orho  + \lambda\Reg\bE_{x\sim\hP}\left[\log\left(
\fh_x
+
e^{-H^{\mathrm{gap}}(x)/(\lambda\Reg)}
\fh_{x}^c
\right)1_D(x)\right]=0,
\]
and
\begin{align*}
&\lim_{\lambda\to0}~\nabla\left[ 
\lambda\orho  + \lambda\Reg\bE_{x\sim\hP}\left[\log\left(
\fh_x
+
e^{-H^{\mathrm{gap}}(x)/(\lambda\Reg)}
\fh_{x}^c
\right)1_D(x)\right]
\right]\\
=&\orho + \Reg\bE_{x\sim\hP}\left[\log\left( 
\fh_x
\right)1_D(x)\right]
\le\orho + \Reg\log(\tau)\hP\{D\}\le -\orho<0,
\end{align*}
where the second inequality is by taking the constant $\tau=\exp\left(-\frac{2\orho}{\Reg\hP\{D\}}\right)$.
Hence, there exists $\overline{\lambda}>0$ such that 
\begin{align*}
&v(\overline{\lambda})
\le 
H^u+
\overline{\lambda}\orho
+
\overline{\lambda}\Reg\bE_{x\sim\hP}\left[ \log\left(
\fh_x
+
e^{-H^{\mathrm{gap}}(x)/(\overline{\lambda}\Reg)}
\fh_{x}^c
\right)1_D(x)\right]<v(0),
\end{align*}
which contradicts to the optimality of $\lambda^*=0$.
As a result, almost surely for all $x$, we have that
\[
\mathbb{E}_{z\sim\bQ_{x,\Reg}}[1_{A}(z)]>0.
\]
To show the second condition, we re-write the dual objective function for $\lambda>0$ as 
\begin{align*}
&v(\lambda)=\lambda\orho + \lambda\Reg\bE_{x\sim\hP}\left[\log
\left(
\mathbb{E}_{z\sim\bQ_{x,\Reg}}[1_{A}(z)] + \mathbb{E}_{z\sim\bQ_{x,\Reg}}\left[ 
e^{[f(z)-H^u]/(\lambda\Reg)}1_{A^c}(z)
\right]
\right)
\right] + H^u.
\end{align*}
The gradient of $v(\lambda)$ becomes
\begin{align*}
\nabla v(\lambda)&=\orho + \Reg\bE_{x\sim\hP}\left[\log
\left(
\mathbb{E}_{z\sim\bQ_{x,\Reg}}[1_{A}(z)] + \mathbb{E}_{z\sim\bQ_{x,\Reg}}\left[ 
e^{[f(z)-H^u]/(\lambda\Reg)}1_{A^c}(z)
\right]
\right)
\right]\\
&\qquad + \bE_{x\sim\hP}\left[
\frac{\mathbb{E}_{z\sim\bQ_{x,\Reg}}\left[ 
e^{[f(z)-H^u]/(\lambda\Reg)}1_{A^c}(z)(H^u-f(z))/(\lambda)
\right]}{\mathbb{E}_{z\sim\bQ_{x,\Reg}}[1_{A}(z)] + \mathbb{E}_{z\sim\bQ_{x,\Reg}}\left[ 
e^{[f(z)-H^u]/(\lambda\Reg)}1_{A^c}(z)
\right]}
\right].
\end{align*}
We can see that $\lim_{\lambda\to\infty}\nabla v(\lambda)=\orho$.
Take
\[
v_{1,x}(\lambda) =  \mathbb{E}_{z\sim\bQ_{x,\Reg}}\left[
e^{[f(z)-H^u]/(\lambda\Reg)}1_{A^c}(z)
\right].
\]
Then $\lim_{\lambda\to0}v_{1,x}(\lambda)=0$ and $v_{1,x}(\lambda)\ge0$.
Take
\[
v_{2,x}(\lambda) =  
\frac{\mathbb{E}_{z\sim\bQ_{x,\Reg}}\left[ 
e^{[f(z)-H^u]/(\lambda\Reg)}1_{A^c}(z)(H^u-f(z))/(\lambda)
\right]}{\mathbb{E}_{z\sim\bQ_{x,\Reg}}[1_{A}(z)] + \mathbb{E}_{z\sim\bQ_{x,\Reg}}\left[ 
e^{[f(z)-H^u]/(\lambda\Reg)}1_{A^c}(z)
\right]}.
\]
Then $\lim_{\lambda\to0}v_{2,x}(\lambda)=0$ and $v_{2,x}(\lambda)\ge0$.
It follows that
\[
\lim_{\lambda\to0}\nabla v(\lambda ) = 
\orho+\Reg\bE_{x\sim\hP}\left[\log
\mathbb{E}_{z\sim\bQ_{x,\Reg}}[1_{A}(z)]
\right]=\orho'.
\]
Hence, if the last condition is violated, based on the mean value theorem, we can find $\overline{\lambda}>0$ so that $\nabla v(\overline{\lambda})=0$, which contradicts to the optimality of $\lambda^*=0$.

Now we show the converse direction.
For any $\lambda>0$, we find that 
\[
\nabla v(\lambda) = \orho + \Reg\bE_{x\sim\hP}\left[\log
\left(
\mathbb{E}_{z\sim\bQ_{x,\Reg}}[1_{A}(z)] + v_{1,x}(\lambda)
\right)
\right] + \bE_{x\sim\hP}[v_{2,x}(\lambda)].
\]
For fixed $x$, when $\mathbb{E}_{\bQ_{x,\Reg}}[1_{A}]=1$, we can see that $v_{1,x}(\lambda)=v_{2,x}(\lambda)=0$, then
\[
\orho + \Reg\left[\log
\left(
\mathbb{E}_{z\sim\bQ_{x,\Reg}}[1_{A}(z)] + v_{1,x}(\lambda)
\right)
\right] + v_{2,x}(\lambda) = \orho>0.
\]
When $\mathbb{E}_{z\sim\bQ_{x,\Reg}}[1_{A}(z)]\in(0,1)$, we can see that $v_{1,x}(\lambda)>0, v_{2,x}(\lambda)>0$.
Then
\[
\orho + \Reg\left[\log
\left(
\mathbb{E}_{z\sim\bQ_{x,\Reg}}[1_{A}(z)] + v_{1,x}(\lambda)
\right)
\right] + v_{2,x}(\lambda) > \orho + \Reg\log(\mathbb{E}_{z\sim\bQ_{x,\Reg}}[1_{A}(z)]) = \orho'\ge0.
\]
Therefore, $\nabla v(\lambda)>0$ for any $\lambda>0$.
By the convexity of $v(\lambda)$, the dual minimizer $\lambda^*=0$.
\QED\endproof

\Jie{
\proof{Proof of Lemma~\ref{Lemma:lambda:positive}.}
Recall that $v(\lambda)$ denotes the objective function for the dual
problem.
The optimality condition can be derived by taking $\nabla_{\lambda}~v(\lambda)\mid_{\lambda=\lambda^*}=0$.
To show the uniqueness of $\lambda^*$, we find that
\begin{align*}
&\nabla^2_{\lambda}v(\lambda)\\=
&\frac{1}{\lambda^3\Reg}\bE_{x\sim\hP}\left[ 
\Big( 
\bE_{z\sim\bQ_{x,\Reg}}[e^{f(z)/(\lambda\Reg)}]
\Big)^{-2}\cdot 
\Big( 
\bE_{z\sim\bQ_{x,\Reg}}[e^{f(z)/(\lambda\Reg)}f^2(z)]\bE_{z\sim\bQ_{x,\Reg}}[e^{f(z)/(\lambda\Reg)}]
-
\big\{\bE_{z\sim\bQ_{x,\Reg}}[e^{f(z)/(\lambda\Reg)}f(z)]\big\}^2
\Big)
\right].
\end{align*}
It can be shown by the Cauchy-Schwarz inequality that $\nabla^2_{\lambda}v(\lambda)\ge0$ for any $\lambda>0$, and the equality holds if and only if $f(\cdot)$ is a constant.
If it is the case, the dual objective $v(\lambda)$ has the unique minimizer $\lambda^*=0$, which contradicts to our assumption.
Hence, strict convexity holds for the dual objective and it implies the uniquess of $\lambda^*$.
\QED\endproof
}

\proof{Proof of Theorem~\ref{Theorem:strong:duality}.}
Recall the feasibility result in Theorem~\ref{Theorem:strong:duality}\ref{Theorem:strong:duality:I} can be easily shown by considering the reformulation of $V_{\Primal}$ in Lemma~\ref{Lemma:refor:Inf-W-E} and the non-negativity of KL-divergence.
When $\orho=0$, one can see that 
\begin{align*}
V_{\Dual}&=\inf_{\lambda\ge0}~\left\{\lambda\Reg\,\mathbb{E}_{x\sim\hP}\Big[\log
\mathbb{E}_{z\sim\bQ_{x,\Reg}}\big[ 
e^{f(z)/(\lambda\Reg)}
\big]\Big]\right\}\\
&\le \lim_{\lambda\to\infty}~\lambda\Reg\,\mathbb{E}_{x\sim\hP}\Big[\log
\mathbb{E}_{z\sim\bQ_{x,\Reg}}\big[ 
e^{f(z)/(\lambda\Reg)}
\big]\Big]
=\bE_{x\sim\hP}\bE_{z\sim \bQ_{x,\Reg}}[f(z)]=V.
\end{align*}
Therefore, the strong duality result holds in this case.
Theorem~\ref{Theorem:strong:duality}\ref{Theorem:strong:duality:IV} can be shown by Lemma~\ref{Lemma:necessary:lambda:0}.
It remains to show the strong duality result for $\orho>0$, which can be further separated to two cases: Condition~\ref{Assumption:light:tailed:f} holds or not.
\begin{itemize}
    \item 
When Condition~\ref{Assumption:light:tailed:f} holds, by Lemma~\ref{Lemma:finite:dual}, the dual minimizer $\lambda^*$ exists. 
The proof for $\lambda^*>0$ can be found in main context.
When $\lambda^*=0$, the optimality condition in Lemma~\ref{Lemma:necessary:lambda:0} holds.
We construct the primal (approximate) solution $\bP_\ast=\proj_{2\#}\gamma_*$, where $\gamma_*$ satisfies
\[
\diff\gamma_*(x,z) = \diff\gamma^x_*(z)\diff\hP(x),\quad
\text{where }\diff\gamma^x_*(y) =
\begin{cases}
0,&\quad\text{if }z\notin A,\\
\frac{e^{-c(x,z)/\Reg}\diff\nu(z)}{\bE_{u\sim\nu}\left[e^{-c(x,u)/\Reg}1_{A}\right]},&\quad\text{if }z\in A.
\end{cases}
\]
We can verify easily that the primal solution is feasible based on the optimality condition $\orho'\ge0$ in Lemma~\ref{Lemma:necessary:lambda:0}.
Moreover, we can check that the primal optimal value is lower bounded by the dual optimal value:
\begin{align*}
V_{\Primal}\ge&
\bE_{(x,z)\sim\gamma_{\ast}}[f(z)]=
\bE_{x\sim\hP}\bE_{z\sim \gamma^x_*}[f(z)]
=
\bE_{x\sim\hP}\bE_{z\sim \gamma^x_*}\left[\underset{\nu}{\esssup}~f\right]
=\underset{\nu}{\esssup}~f=V_{\Dual},
\end{align*}
where the second equality is because that $z\in A$ so that $f(z)=\underset{\nu}{\esssup}~f$.
This, together with the weak duality result, completes the proof in this part. 
\item
When Condition~\ref{Assumption:light:tailed:f} does not hold, we consider a sequence of real numbers $\{R_j\}_j$ such that $R_j\to\infty$ and take the objective function $f_j(z) = f(z)1\{f(z)\le R_j\}$.
Hence, there exists $\lambda>0$ satisfying $\Pr{}_{x\sim\hP}\left\{
x:~\mathbb{E}_{\bQ_{x,\Reg}}\left[ 
e^{f_j(z)/(\lambda\Reg)}
\right]=\infty
\right\}=0$.
According to the necessary condition in Lemma~\ref{Lemma:necessary:lambda:0}, the corresponding dual minimizer $\lambda_j^*>0$ for sufficiently large index $j$.
Then we can apply the duality result in the first part of Theorem~\ref{Theorem:strong:duality}\ref{Theorem:strong:duality:III} to show that for sufficiently large $j$, it holds that
\[
\sup_{\mathbb{P}\in\mathbb{B}_{\rho,\Reg}(\hP)}~
\left\{\mathbb{E}_{z\sim \mathbb{P}}[f_j(z)]\right\}\ge 
\lambda^*_j\orho+\lambda^*_j\Reg\bE_{x\sim\hP}\left[\log
\mathbb{E}_{z\sim\bQ_{x,\Reg}}\Big[ 
e^{f_j(z)/(\lambda\Reg)}
\Big]\right].
\]
Taking $j\to\infty$ both sides implies that $V_{\Primal}=\infty$.
\QED
\end{itemize}
\endproof

\section{Proof of Theorem~\ref{Theorem:complexity:BSMD} in Section~\ref{Sec:inner:step:1}}
\label{Appendix:convergence:analysis}
\Jie{
In this section, we omit the dependence of $\lambda$ when defining objective or subgradient terms, e.g., we write $F(\theta)$ for $F(\theta;\lambda)$.
}
We first present some preliminaries that can be useful for developing the proof result in Section~\ref{Sec:inner:step:1}.
As any two norms on a finite-dimensional vector space are equivalent, we impose the following assumption throughout Section~\ref{Sec:first:order} without loss of generality:
\begin{assumption}\label{assum:norm:bound}
There exists $\mathfrak{c}$ and $\mathfrak{d}$ such that 
$\mathfrak{c}
\|\cdot\|_2
\le
\|\cdot\|\le \mathfrak{d}\|\cdot\|_2.$
\end{assumption}
By Assumption~\ref{assum:norm:bound}, we obtain the bound regarding the dual norm $\|\cdot\|_*$:
\[
\mathfrak{d}^{-1}\cdot\|\cdot\|_2\le \|\cdot\|_*\le \mathfrak{c}^{-1}\cdot\|\cdot\|_2.
\]
The complexity result of our proposed gradient estimators is summarized below.

\begin{remark}[Complexity of Gradient Estimators~{\citep[Appendix~B]{hu2021biasvar}}]\label{Proposition:computation:cost}
\Jie{
To generate the SG estimator $v^{\text{SG}}(\theta)$,
one needs to generate $1$ sample from $\hP$ and $2^L$ samples from $\bQ_{x,\Reg}$ for some $x\in\mathrm{supp}\, \hP$. 
To generate the RT-MLMC estimator $v^{\text{RT-MLMC}}(\theta)$,
one needs to generate $1$ sample from $\hP$ and the required (expected) number of samples from $\bQ_{x,\Reg}$ for some $x\in\mathrm{supp}\, \hP$ equals $\frac{L}{2 - 2^{-L}} = \cO(L)$.
}
\QEG
\end{remark}

Next, we present some basic properties regarding the approximation function $F^{\ell}(\theta)$ defined in \eqref{Eq:stat:robust:formula:approximation} in Lemma~\ref{Facts:about:F}, which can be used to show Theorem~\ref{Theorem:complexity:BSMD}.
Recall that we defined the constant $K_{\lambda,\Reg,B}=B/(\lambda\Reg)$.
\begin{lemma}
\label{Facts:about:F}
\begin{enumerate}
    \item\label{Facts:about:F:I}
Under Assumption~\ref{Assumption:throughout:loss}\ref{Assumption:throughout:loss:bound}, it holds that
\[
\big| 
F^{\ell}(\theta) - F(\theta)
\big|\le \lambda\Reg \exp\left(2K_{\lambda,\Reg,B}\right)\cdot 2^{-(\ell+1)},\qquad \forall \theta\in\Theta.
\]
    \item\label{Facts:about:F:II}
Under Assumption~\ref{Assumption:throughout:loss}\ref{Assumption:throughout:loss:bound} and \ref{Assumption:throughout:loss}\ref{Assumption:throughout:loss:lip}, it holds that
\[
\big\|
\nabla F^{\ell}(\theta) - \nabla F(\theta)
\big\|_2^2\le L_f^2\exp\left(4K_{\lambda,\Reg,B}\right)\cdot 2^{-\ell},\qquad \forall \theta\in\Theta.
\]
\item\label{Facts:about:F:IV}
Under Assumption~\ref{Assumption:throughout:loss}\ref{Assumption:throughout:loss:lip}, it holds that
\[
\bE\left[\big\|
g^{\ell}(\theta, \zeta^{\ell})
\big\|^2_2\right]\le L_f^2,\qquad \forall \theta\in\Theta.
\]
Additionally when Assumption~\ref{Assumption:throughout:loss}\ref{Assumption:throughout:loss:bound} holds, it holds that
\[
\bE\left[\big\|
G^{\ell}(\theta, \zeta^{\ell})
\big\|^2_2\right]\le L_f^2\exp\left(4K_{\lambda,\Reg,B}\right)\cdot 2^{-\ell},\qquad \forall \theta\in\Theta.
\]
\end{enumerate}
\end{lemma}

\Jie{
\proof{Proof of Lemma~\ref{Facts:about:F}.}
Recall that \eqref{Eq:F:lambda} is a special CSO problem in \eqref{Eq:CSO:speific}, by taking  $H^1(\boldsymbol{\cdot})=\lambda\Reg\log(\boldsymbol{\cdot})$ and $H^2(\boldsymbol{\cdot},z)=\exp(f_{\boldsymbol{\cdot}}(z) / (\lambda\Reg))$.
Under the assumptions stated in Lemma~\ref{Facts:about:F}, it can be shown that $H^2(\boldsymbol{\cdot},z)$ is $\exp(K_{\lambda,\Reg,B})$-uniformly bounded, $\exp(K_{\lambda,\Reg,B})L_f/(\lambda\Reg)$-Lipschitz continuous. The function $H^1(\boldsymbol{\cdot})$ has the domain set $[1, \exp(K_{\lambda,\Reg,B})]$, and is therefore $\lambda\Reg$-Lipschitz continuous and $\lambda\Reg$-smooth.
Thus, the desired results hold by applying \citep[Lemma~3.1]{hu2020sample} and \citep[Proposition~4.1]{hu2021biasvar}.
\endproof
}

\subsection{Proof of Theorem~\ref{Theorem:complexity:BSMD}}
\Jie{
We first study the convergence guarantees for solving a generic nonsmooth convex optimization $\min\limits_{\theta\in\Theta}~F(\theta)$.
Let $\overline{F}(\theta)$ denote its approximation, with the approximation bias $\Delta_F$ satisfying
\[
|\overline{F}(\theta) - F(\theta)|\le \Delta_F,\quad \forall \theta\in\Theta.
\]
Denote by $\nabla \overline{F}(\theta)$ a subgradient of $\overline{F}$ at $\theta$.
Suppose for a given $\theta$, the subgradient estimate of $F(\theta)$, denoted as $v(\theta)$, satisfies
\[
\bE[v(\theta)] = \nabla \overline{F}(\theta),\quad
\bE\left[ 
\|v(\theta)\|_*^2
\right]\le M_*^2.
\]
Let $\overline{\theta}^* \in \argmin\limits_{\theta\in\Theta}~\overline{F}(\theta)$ and $\theta^*\in \argmin\limits_{\theta\in\Theta}~{F}(\theta)$.
We then establish the following result.
\begin{lemma}[BSMD for Nonsmooth Convex Optimization]\label{Lemma:SGD:scheme:convex:nonsmooth}
Under the assumptions stated above and with the initial guess $\theta_0\in\Theta$, 
consider the BSMD algorithm that generates the following iteration:
\[
\theta_{t+1} = \prox_{\theta_t}\big( h v(\theta_t)\big),\quad \theta_0\in\Theta,\quad t=0,\ldots,T-1,
\]
where the stepsize parameter $h = \sqrt{\frac{2\kappa D_{\omega}(\theta_0, \overline{\theta}^*)}{TM_*^2}}$.
Let the estimated optimal solution generated by BSMD algorithm be 
$\widehat{\theta} = \frac{1}{T}\sum_{t=1}^T\theta_t$.
Then, the suboptimality gap satisfies:
\[
\bE\big[F(\widehat{\theta}) - F(\theta^*)\big]\le 2\Delta_F + M_*\sqrt{\frac{2D_{\omega}(\theta_0, \overline{\theta}^*)}{\kappa T}}.
\]
\end{lemma}
\begin{remark}
If the approximation bias is zero (i.e., $\Delta_F=0$), the BSMD algorithm reduces to the standard SMD studied in \citep{nemirovski2009robust}.
By \citep[Section~2.3]{nemirovski2009robust}, the suboptimality gap in Lemma~\ref{Lemma:SGD:scheme:convex:nonsmooth} is bounded by $M_*\sqrt{\frac{2D_{\omega}(\theta_0, \overline{\theta}^*)}{\kappa T}}$.
For the case where $\Delta_F>0$, the proof of Lemma~\ref{Lemma:SGD:scheme:convex:nonsmooth} follows from the decomposition argument similar to \citep[Eq.~(9)]{hu2021biasvar}.
However, our result generalizes to the BSMD algorithm with (potentially) nonsmooth loss functions, whereas \citep{hu2021biasvar} focuses only on the SGD algorithm for unconstrained optimization with smooth loss functions.
 \QEG
\end{remark}
Now we are ready to show complexity results for BSMD using SG and RT-MLMC estimators.
Both estimators rely on the same approximation function $F^L(\theta)$ defined in \eqref{Eq:stat:robust:formula:approximation}.
By Lemma~\ref{Facts:about:F}\ref{Facts:about:F:I}, 
$\Delta_F = \lambda\Reg\exp\left(2K_{\lambda,\Reg,B}\right)\cdot 2^{-(L+1)}.$
We now analyze each estimator separately.
\\\noindent
{\bf SG.}
It can be shown from the first part of Lemma~\ref{Facts:about:F}\ref{Facts:about:F:IV} that 
$\bE\left[ 
\|v^{\mathrm{SG}}(\theta)\|_*^2
\right]\le
\big( 
M_*^{\mathrm{SG}}
\big)^2:=
 \mathfrak{c}^{-2}L_f^2.$
To obtain $\delta$-optimal solution for SG estimator, by Lemma~\ref{Lemma:SGD:scheme:convex:nonsmooth}, it suffices to ensure
\[
2\Delta_F\le \frac{\delta}{2},\quad
M_*^{\mathrm{SG}}\sqrt{\frac{2D_{\omega}(\theta_0, \overline{\theta}^*)}{\kappa T}}\le\frac{\delta}{2}.
\]
To satisfy these conditions, we specify the following hyper-parameters:
\[
L=\left\lceil\frac{1}{\log2}\left[\log\frac{2\lambda\Reg \exp(2K_{\lambda,\Reg,B})}{\delta}\right]\right\rceil,\quad 
T=\left\lceil\frac{8L_f^2D_{\omega}(\theta_0, \overline{\theta}^*)}{\kappa\mathfrak{c}^{2}\delta^2}\right\rceil,\quad 
h=\sqrt{\frac{2\kappa\mathfrak{c}^2D_{\omega}(\theta_0, \overline{\theta}^*)}{TL_f^2}}.
\]
\noindent
{\bf RT-MLMC.}
By the second part of Lemma~\ref{Facts:about:F}\ref{Facts:about:F:IV} and basic calculation, we find
\begin{align*}
&\bE\left[\big\|
v^{\mathrm{RT-MLMC}}(\theta)
\big\|^2_*\right]
\le  \mathfrak{c}^{-2}\bE\left[\big\|
v^{\mathrm{RT-MLMC}}(\theta)
\big\|^2_2\right]
=
\sum_{\ell=0}^L
\frac{1}{p_{\ell}}
\bE\left[\big\|
G^{\ell}(\theta, \zeta_1^{\ell})
\big\|^2_2\right]
\\
\le&\big( 
M_*^{\mathrm{RT-MLMC}}
\big)^2:= 
2(L+1)L_f^2\exp(4K_{\lambda,\Reg,B}).
\end{align*}
Similar to the case of SG, we ensure
\[
2\Delta_F\le \frac{\delta}{2},\quad
M_*^{\mathrm{RT-MLMC}}\sqrt{\frac{2D_{\omega}(\theta_0, \overline{\theta}^*)}{\kappa T}}\le\frac{\delta}{2}.
\]
To satisfy these conditions, we select the following hyper-parameters:
\begin{multline*}
L=\left\lceil\frac{1}{\log2}\left[\log\frac{2\lambda\Reg\exp(2K_{\lambda,\Reg,B})}{\delta}\right]\right\rceil,\quad 
T=\left\lceil\frac{16(L+1)L_f^2D_{\omega}(\theta_0, \overline{\theta}^*)\exp(4K_{\lambda,\Reg,B})}{\kappa\mathfrak{c}^2\delta^2}\right\rceil,\\ 
h=\sqrt{\frac{2\kappa D_{\omega}(\theta_0, \overline{\theta}^*)}{T\big( 
M_*^{\mathrm{RT-MLMC}}
\big)^2}}.
\end{multline*}
By Remark~\ref{Proposition:computation:cost}, when running BSMD with SG estimator, the sample complexity from $\hP$ equals $\cO(T)$ and that from $\bQ_{x,\Reg}$ equals $\cO(T2^L)$; when running BSMD with RT-MLMC estimator, the sample complexity from $\hP$ equals $\cO(T)$ and that from $\bQ_{x,\Reg}$ equals $\cO(TL)$.
Substituting the expressions of $T, L$ gives the desired result.
}

\section{Proofs of Technical Results in Section~\ref{Sec:inner:iteration:step:2}}
\label{proof:Sec:inner:iteration:step:2}
\Jie{
We first provide two technical lemmas that can be useful to show the main results in Section~\ref{Sec:inner:iteration:step:2}.
\begin{lemma}\label{Lemma:stat:A}
Under Assumption~\ref{Assumption:throughout:loss}\ref{Assumption:throughout:loss:bound}, it holds that $\bE[(A^{\ell}(\theta, \zeta^{\ell};\lambda))^2]\le 
\lambda^2\Reg^2\exp(2K_{\lambda,\Reg,B})\cdot 2^{-\ell}.$
\end{lemma}
\proof{Proof of Lemma~\ref{Lemma:stat:A}.}
The proof follows the similar procedure from \citep[Proposition~4.1]{hu2021biasvar}.
\endproof
\begin{lemma}[Complexity of RT-MLMC-based Objective Estimator]\label{Proposition:MLMC:obj}
Let error probability $\alpha\in(0,1)$ and accuracy level $\delta>0$.
Assume  Assumption~\ref{Assumption:throughout:loss}\ref{Assumption:throughout:loss:bound} holds and 
specify 
\begin{equation}
L = \left\lceil\frac{1}{\log2}\left[\log\frac{\lambda\Reg\exp(2K_{\lambda,\Reg,B})}{\delta}\right]\right\rceil,\quad 
m'=\cO(1)\frac{\lambda^2\Reg^2\exp(2K_{\lambda,\Reg,B})(L+1)}{\delta^2}\cdot\log\frac{2}{\alpha}.\label{Eq:lemma:MLMC:obj}
\end{equation}
Then, the RT-MLMC estimator~\eqref{Eq:RTMLMC:est} has an accuracy error $\delta$ with probability at least $1-\alpha$.
Its sample complexity from $\hP$ equals $\cO(m')=
\tO(
\lambda^2\Reg^2K_{\lambda,\Reg,B}\exp(2K_{\lambda,\Reg,B})\cdot \delta^{-2}
)
$ and that from $\bQ_{x,\Reg}$ equals $\cO(m'\cdot L) = \tO(
\lambda^2\Reg^2K_{\lambda,\Reg,B}^2\exp(2K_{\lambda,\Reg,B})\cdot \delta^{-2}
)$.
Here $\tO(\cdot)$ hides constants linearly depending on $(\log\frac{\lambda\Reg}{\delta})^2$ and $\log\frac{1}{\alpha}$.
\end{lemma}
\proof{Proof of Lemma~\ref{Proposition:MLMC:obj}.}
We first specify $L$ as in \eqref{Eq:lemma:MLMC:obj} such that $|F^{L}(\theta;\lambda) - F(\theta;\lambda)|\le\frac{\delta}{2}$.
The RT-MLMC estimator~\eqref{Eq:RTMLMC:est} satisfies that 
\begin{align*}
\bE[\widehat{F}(\theta;\lambda)] &= F^L(\theta;\lambda),\\
\Var\left(
\widehat{F}(\theta;\lambda)
\right)&\le \frac{1}{m'}\sum_{\ell=0}^L\frac{1}{p_{\ell}}\bE[(A^{\ell}(\theta, \zeta^{\ell}))^2]
\le
\frac{1}{m'}
\lambda^2\Reg^2\exp(2K_{\lambda,\Reg,B})\cdot(L+1).
\end{align*}
Consequently, there exists $\delta'>0$ such that
\begin{align*}
&\text{Pr}\left\{
|F(\theta;\lambda) - \widehat{F}(\theta;\lambda)|>\delta
\right\}
\le\text{Pr}\left\{
|F^L(\theta;\lambda) - \widehat{F}(\theta;\lambda)|>\frac{\delta}{2}
\right\}\\
\le&2\exp\left( 
-\frac{\delta^2}{4(\delta' + 2)\Var\left(
\widehat{F}(\theta;\lambda)
\right)}
\right)\le 
2\exp\left( 
-\frac{\delta^2m'}{4(\delta'+2)\lambda^2\Reg^2\exp(2K_{\lambda,\Reg,B})(L+1)}
\right),
\end{align*}
where the second inequality is based on the Cramer's large deviation theorem~\citep{kleywegt2002sample}, and the last inequality is by the upper bound on $\Var\left(
\widehat{F}(\theta;\lambda)
\right)$.
To make the desired coverage probability, we take $m'$ as in \eqref{Eq:lemma:MLMC:obj}.
The complexity results are derived by standard calculation similar to Remark~\ref{Proposition:computation:cost}.
\QED
\endproof
}

In the following, we provide the proof of Proposition~\ref{Proposition:com:estimate:optval}.

\Jie{
\proof{Proof of Proposition~\ref{Proposition:com:estimate:optval}.}
Denote by $\theta^*=\argmin\limits_{\theta\in\Theta}~F(\theta;\lambda)$.
The goal is to choose hyper-parameters such that 
\[
\text{Pr}\bigg\{ 
\left|\min_{i\in[m]}~\hF(\widehat{\theta}_i; \lambda) - F(\theta^*;\lambda)\right|\le \delta
\bigg\}\ge 1 - \eta.
\]
On the one hand, 
\[
\min_{i\in[m]}~\hF(\widehat{\theta}_i; \lambda) - F(\theta^*;\lambda)\le 
\min_{i\in[m]}~F(\widehat{\theta}_i; \lambda) - F(\theta^*;\lambda)
+
\max_{i\in[m]}~|F(\widehat{\theta}_i; \lambda) - \hF(\widehat{\theta}_i; \lambda)|.
\]
On the other hand,
\[
F(\theta^*;\lambda) - \min_{i\in[m]}~\hF(\widehat{\theta}_i; \lambda)\le 
F(\theta^*;\lambda) - \min_{i\in[m]}~F(\widehat{\theta}_i; \lambda) +\max_{i\in[m]}~|F(\widehat{\theta}_i; \lambda) - \hF(\widehat{\theta}_i; \lambda)|
\le 
\max_{i\in[m]}~|F(\widehat{\theta}_i; \lambda) - \hF(\widehat{\theta}_i; \lambda)|.
\]
Based on those two inequalities, it suffices to choose hyper-parameters such that
\begin{equation}\label{Eq:VF:bound}
\text{Pr}\bigg\{ 
\max_{i\in[m]}~|F(\widehat{\theta}_i; \lambda) - \hF(\widehat{\theta}_i; \lambda)|\le \frac{\delta}{2}
\bigg\}\ge 1 - \frac{\eta}{2}
\end{equation}
and
\begin{equation}\label{Eq:FF:bound}
\text{Pr}\bigg\{ 
\min_{i\in[m]}~F(\widehat{\theta}_i; \lambda) - F(\theta^*;\lambda)\le \frac{\delta}{2}
\bigg\}\ge 1 - \frac{\eta}{2}.
\end{equation}
To ensure the relation~\eqref{Eq:VF:bound}, it suffices to apply Lemma~\ref{Proposition:MLMC:obj} with error probability $\frac{\eta}{2m}$ and accuracy level $\delta/2$.
It implies that the sample complexity from $\hP$ at Step~3 of Algorithm~\ref{Alg:inexact:obj} for each independent repetition is $\tO(
\lambda^2\Reg^2K_{\lambda,\Reg,B}\exp(2K_{\lambda,\Reg,B})\cdot \delta^{-2}
)$, and that from $\bQ_{x,\Reg}$ is $\tO(
\lambda^2\Reg^2K_{\lambda,\Reg,B}^2\exp(2K_{\lambda,\Reg,B})\cdot \delta^{-2})
$.
To ensure the relation~\eqref{Eq:FF:bound}, it suffices to take
\[
\text{Pr}\bigg\{ 
F(\widehat{\theta}_i;\lambda) - F(\theta^*;\lambda)\le \frac{\delta}{2}
\bigg\}\ge 1 - \left(\frac{\eta}{2}\right)^{1/m},\quad \forall i\in[m].
\]
By Markov's inequality, it suffices to ensure 
\begin{equation}
\bE[F(\widehat{\theta}_i;\lambda) - F(\theta^*;\lambda)]\le \frac{\delta}{2}\left(\frac{\eta}{2}\right)^{1/m},\quad \forall i\in[m].
\label{Eq:relation:final}
\end{equation}
By Theorem~\ref{Theorem:complexity:BSMD}\ref{Theorem:complexity:BSMD:smooth} with accuracy level $\frac{\delta}{2}\left(\frac{\eta}{2}\right)^{1/m}$, the sample complexity from $\hP$ at Step~2 of Algorithm~\ref{Alg:inexact:obj} for each independent repetition is $\tO(K_{\lambda,\Reg,B}\exp(4K_{\lambda,\Reg,B})\cdot\delta^{-2}\eta^{-2/m})$, and that from $\bQ_{x,\Reg}$ is $\tO(K_{\lambda,\Reg,B}^2\exp(4K_{\lambda,\Reg,B})\cdot\delta^{-2}\eta^{-2/m})$.
Therefore, the sample complexity from $\hP$ of Algorithm~\ref{Alg:inexact:obj} is
\[
\begin{aligned}
&m\cdot\Big[ 
\tO(K_{\lambda,\Reg,B}\exp(4K_{\lambda,\Reg,B})\cdot\delta^{-2}\eta^{-2/m}) + \tO(
\lambda^2\Reg^2K_{\lambda,\Reg,B}\exp(2K_{\lambda,\Reg,B})\cdot \delta^{-2}
)
\Big]\\
=&\tO( H_{\lambda,\Reg,B} K_{\lambda,\Reg,B}\exp(2K_{\lambda,\Reg,B})\delta^{-2}\cdot m(1 + \eta^{-2/m}))
\end{aligned}
\]
and that from $\bQ_{x,\Reg}$ is
\[
\begin{aligned}
&m\cdot\Big[ 
\tO(K_{\lambda,\Reg,B}^2\exp(4K_{\lambda,\Reg,B})\cdot\delta^{-2}\eta^{-2/m}) + \tO(
\lambda^2\Reg^2K_{\lambda,\Reg,B}^2\exp(2K_{\lambda,\Reg,B})\cdot \delta^{-2})
\Big]\\
=&\tO( H_{\lambda,\Reg,B} K_{\lambda,\Reg,B}^2\exp(2K_{\lambda,\Reg,B})\delta^{-2}\cdot m(1 + \eta^{-2/m})).
\end{aligned}
\]
In the above deviation, we defined the constant $H_{\lambda,\Reg,B} = \max(\exp(2K_{\lambda,\Reg,B}), \lambda^2\Reg^2)$.
Hence, it suffices to specify $m$ such that $\tO(m(1 + \eta^{-2/m}))$ is minimized. One valid choice is $m=\lceil\log_2\frac{2}{\eta}\rceil$, which leads to the desired complexity bounds.
\QED\endproof
}

Finally, we show the proof of Theorem~\ref{Proposition:complexity:bisection}.
A key technique is the following complexity result on bisection search with inexact oracles.
\Jie{
\begin{lemma}[Complexity for Noisy Bisection]\label{Lemma:noisy:bisection:complexity}
Let the accuracy level $\delta>0$, and $\Psi:~\bR \to \bR$ be a $L_{\Psi}$-Lipschitz continuous and convex function defined on the interval $[\lambda_l,\lambda_u]$.
Assume there exists an oracle $\widehat{\Psi}:~\bR\to\bR$ such that $|\widehat{\Psi}(\lambda) - \Psi(\lambda)|\le \delta, \forall \lambda$.
Let us run Algorithm~\ref{Alg:outer:bisection} for $T'=\lceil \log_2\big(\frac{L_{\Psi}(\lambda_u - \lambda_l)}{\delta}\big)\rceil$ iterations, then with at most $3 + 2T'$ calls to $\widehat{\Psi}$, Algorithm~\ref{Alg:outer:bisection} outputs $\widehat{\lambda}$ so that
\[
\Psi(\widehat{\lambda}) - \min_{\lambda\in [\lambda_l,\lambda_u]}~\Psi(\lambda)\le 4{\delta}.
\]
\end{lemma}
\proof{Proof of Lemma~\ref{Lemma:noisy:bisection:complexity}.}
The proof is straightforward by following~\citep[Lemma 33]{cohen2016geometric}
\QED
\endproof
}

\Jie{
\proof{Proof of Theorem~\ref{Proposition:complexity:bisection}.}
It can be verified that $\Psi$ is a convex function with a subgradient 
\begin{align*}
\frac{\partial}{\partial\lambda}\Psi(\lambda)
&=\orho + \bE_{x\sim\hP}\left[\Reg\log\bE_{z\sim\bQ_{x,\Reg}}\Big[e^{f_{\theta^*_{\lambda}}(z)/(\lambda\Reg)}\Big]\right] - \bE_{x\sim\hP}\left[ 
\frac{\bE_{z\sim\bQ_{x,\Reg}}\left[e^{f_{\theta^*_{\lambda}}(z)/ (\lambda\Reg)}f_{\theta^*_{\lambda}}(z)\right]}{\lambda\bE_{z\sim\bQ_{x,\Reg}}\left[e^{f_{\theta^*_{\lambda}}(z)/(\lambda\Reg)}\right]}
\right],
\end{align*}
where $\theta^*_{\lambda}\in\argmin\limits_{\theta\in\Theta}~F(\theta;\lambda)$.
By Assumption~\ref{Assumption:throughout:loss} and $\lambda\in[\lambda_l,\lambda_u]$, this subgradient vector is bounded: 
\[
\left| 
\frac{\partial}{\partial\lambda}\Psi(\lambda)
\right|\le L_{\Psi}:=\orho + \frac{B}{\lambda_{l}}\left[1 + \exp({K_{\lambda_l, \Reg, B})}\right].
\]
In summary, $\Psi(\lambda)$ is a $L_{\Psi}$-Lipschitz and convex function defined on $[\lambda_{l}, \lambda_u]$.
Applying Lemma~\ref{Lemma:noisy:bisection:complexity} with accuracy level $\delta/4$ together with the union bound, we are able to find the optimal multiplier up to accuracy $\delta$ with probability at least $1-\eta$ by calling the oracle $\widehat{\Psi}$ for 
$3 + 2\left\lceil \log_2\big(\frac{4L_{\Psi}(\lambda_u - \lambda_l)}{\delta}\big)\right\rceil$
times.\QED\endproof
}

\begin{lemma}\label{Lemma:last}
Under Assumption~\ref{Assumption:throughout:loss}, the optimal multiplier $\lambda^*$ to \eqref{Eq:reformulate:statistical:learning} satisfies $\lambda^*\le \frac{B}{\orho}$.
\end{lemma}

\proof{Proof.}
It can be verified that
\begin{align*}
0&=\left. \frac{\partial}{\partial\lambda}\Psi(\lambda)\right|_{\lambda=\lambda^*}\\
&=\orho + \bE_{x\sim\hP}\left[\Reg\log\bE_{z\sim\bQ_{x,\Reg}}\Big[e^{f_{\theta^*_{\lambda^*}}(z)/(\lambda\Reg)}\Big]\right] - \bE_{x\sim\hP}\left[ 
\frac{\bE_{z\sim\bQ_{x,\Reg}}\left[e^{f_{\theta^*_{\lambda^*}}(z)/ (\lambda\Reg)}f_{\theta^*_{\lambda^*}}(z)\right]}{\lambda^*\bE_{z\sim\bQ_{x,\Reg}}\left[e^{f_{\theta^*_{\lambda^*}}(z)/(\lambda\Reg)}\right]}
\right]\\
&\ge \orho - \bE_{x\sim\hP}\left[ 
\frac{\bE_{z\sim\bQ_{x,\Reg}}\left[e^{f_{\theta^*_{\lambda^*}}(z)/ (\lambda\Reg)}f_{\theta^*_{\lambda^*}}(z)\right]}{\lambda^*\bE_{z\sim\bQ_{x,\Reg}}\left[e^{f_{\theta^*_{\lambda^*}}(z)/(\lambda\Reg)}\right]}
\right]\\
&\ge \orho - \frac{B}{\lambda^*},
\end{align*}
where the two inequalities is based on the fact that $0\le f_{\theta}(z)\le B$.
The desired result holds directly.
\QED
\endproof

\end{document}